\documentclass{memo-l}

\usepackage{amsmath,amssymb,amscd,amsfonts}
\usepackage[all,cmtip]{xy}
\usepackage{graphicx}
\usepackage{tikz-cd}

\newcommand{\A}{\mathbb{A}}
\newcommand{\bK}{\mathbb{K}}
\newcommand{\Q}{\mathbb{Q}}
\newcommand{\R}{\mathbb{R}}
\newcommand{\Z}{\mathbb{Z}}

\newcommand{\cA}{\mathcal{A}}
\newcommand{\cB}{\mathcal{B}}
\newcommand{\cI}{\mathcal{I}}
\newcommand{\cO}{\mathcal{O}}
\newcommand{\cP}{\mathcal{P}}
\newcommand{\cT}{\mathcal{T}}
\newcommand{\cBT}{\mathcal{BT}}

\newcommand{\BM}{\mathrm{BM}}
\newcommand{\cell}{\mathrm{cell}}
\newcommand{\cl}{\mathrm{cl}}
\newcommand{\Coker}{\mathrm{Coker}}
\newcommand{\cores}{\mathrm{cor}}
\newcommand{\diag}{\mathrm{diag}}
\newcommand{\Ext}{\mathrm{Ext}}
\newcommand{\fin}{\mathrm{fin}}
\newcommand{\GL}{\mathrm{GL}}
\newcommand{\Hom}{\mathrm{Hom}}
\newcommand{\id}{\mathrm{id}}
\newcommand{\Ker}{\mathrm{Ker}}
\newcommand{\Lat}{\mathrm{Lat}}
\newcommand{\Latbar}{\overline{\mathrm{Lat}}}
\newcommand{\Map}{\mathrm{Map}}
\newcommand{\MS}{\mathrm{MS}}
\newcommand{\nd}{\mathrm{nd}}

\newcommand{\PGL}{\mathrm{PGL}}
\newcommand{\sgn}{\mathrm{sgn}}
\newcommand{\SL}{\mathrm{SL}}
\newcommand{\SO}{\mathrm{SO}}
\newcommand{\Tor}{\mathrm{Tor}}
\newcommand{\tot}{\mathrm{tot}}
\newcommand{\tr}{\mathrm{tr}}

\newcommand{\bsl}{\backslash}
\newcommand{\eps}{\epsilon}
\newcommand{\pmone}{{\{\pm 1\}}}

\newcommand{\inj}{\hookrightarrow}
\newcommand{\surj}{\twoheadrightarrow}
\newcommand{\resp}{resp.\ }
\newcommand{\xto}[1]{\xrightarrow{#1}}
\newcommand{\wt}[1]{\widetilde{#1}}
\newcommand{\wh}[1]{\widehat{#1}}

\newcommand{\C}{\mathbb{C}}
\newcommand{\Ind}{\mathrm{Ind}}
\newcommand{\St}{\mathrm{St}}
\newcommand{\cC}{\mathcal{C}}
\newcommand{\cH}{\mathcal{H}}
\newcommand{\Spec}{\mathrm{Spec}\,}
\newcommand{\cF}{\mathcal{F}}
\newcommand{\rank}{\mathrm{rank}}
\newcommand{\cL}{\mathcal{L}}
\newcommand{\Supp}{\mathrm{Supp}\,}
\newcommand{\cD}{\mathcal{D}}
\newcommand{\Mat}{\mathrm{Mat}}
\newcommand{\divi}{\mathrm{div}}
\newcommand{\Flag}{\mathrm{Flag}}
\newcommand{\Vertex}{\mathrm{Vert}}

\newtheorem{theorem}{Theorem}[chapter]

\theoremstyle{definition}
\newtheorem{definition}[theorem]{Definition}

\theoremstyle{remark}
\newtheorem{remark}[theorem]{Remark}

\numberwithin{section}{chapter}
\numberwithin{equation}{chapter}

\newtheorem{thm}{Theorem}
\newtheorem{prop}[thm]{Proposition}
\newtheorem{lem}[thm]{Lemma}
\newtheorem{cor}[thm]{Corollary}
\makeindex

\usepackage{tikz}
\usetikzlibrary{shapes}

\newcommand{\grid}[1]{
\foreach \i [count=\row from 0, remember=\row as \lastrow (initially 0)] in {0,...,#1}{
    \foreach \j [count=\column from 0, remember=\column as \lastcolumn (initially 0)] in {0,...,\i}{
        \ifnum\row=0
            \node[tri](0-0){}; 
        \else
            \ifnum\column=0
                \node[tri, anchor=north](\row-0) at (\lastrow-0.corner 2){}; 
            \else
                \node[tri, anchor=north](\row-\column) at (\lastrow-\lastcolumn.corner 3){}; 
            \fi
        \fi}}
}

\begin{document}

\frontmatter

\title{Arithmetic quotients of the Bruhat-Tits building for projective general linear group in positive characteristic}


\author{Satoshi Kondo}
\address{Corresponding author: 
Satoshi Kondo \\
Middle East Technical University \\
Northern Cyprus Campus, Kalkanli \\
Guzelyurt, Mersin 10, Turkey;
Kavli Institute for the Physics and Mathematics of the Universe\\
University of Tokyo\\
5-1-5 Kashiwanoha\\
Kashiwa 277-8583\\ Japan
}
\email{satoshi.kondo@gmail.com}
\thanks{the first author was partially supported by World Premier International Research Center Initiative (WPI Initiative), MEXT, Japan}

\author{Seidai Yasuda}
\address{Seidai Yasuda \\ Department of Mathematics, Faculty of Science, Hokkaido University
Kita 10, Nishi 8, Kita-Ku, Sapporo, Hokkaido, 060-0810, Japan}
\curraddr{}
\email{sese@math.sci.hokudai.ac.jp}

\date{}




\begin{abstract}
Let $d \ge 1$.  We study a subspace of the space of 
automorphic forms 
of $\GL_d$ over a global field of positive characteristic
(or, a function field of a curve over a finite field).
We fix a place $\infty$ of $F$, and we consider
the subspace $\cA_{\St}$ 
consisting of automorphic
forms such that the local component at $\infty$
of the associated automorphic representation
is the Steinberg representation (to be made precise
in the text).

We have two results.

One theorem (Theorem~\ref{7_PROP1}) 
describes the constituents of $\cA_{\St}$
as automorphic representation and gives a multiplicity one type statement.   

For the other theorem (Theorem~\ref{lem:apartment}), 
we construct, using the geometry of the Bruhat-Tits building,
an analogue of modular symbols in $\cA_{\St}$ integrally
(that is, in the space of $\Z$-valued automorphic forms).
We show that the quotient is finite and give 
a bound on the exponent of this quotient.
\end{abstract}

\maketitle

\tableofcontents

{\Large{\bf{Acknowledgments}}}\\

During this research, the first author was supported as a Twenty-First Century COE Kyoto Mathematics Fellow and was partially supported by JSPS KAKENHI Grant number JP17740016, JP21654002, JP25800005 and by World Premier International Research Center Initiative (WPI Initiative), MEXT, Japan. The second author was partially supported by JSPS KAKENHI Grant number JP16244120, JP21540013, JP24540018, and JP15H03610.

The second author would like to thank Fumiharu Kato for letting him know the works
\cite{We2}, \cite{We1} of
Werner.

We thank Professor Kazuya Kato for always raising important questions.

\mainmatter
\chapter{Introduction}

\section{Modular forms and automorphic forms}
\label{sec:modauto}
Automorphic forms are fundamental objects of study in number theory.   We obtain some very basic results concerning automorphic forms (satisfying some condition at a fixed prime)
for $\GL_d$ (where $d$ is a positive integer) over the
function field of a curve over a finite field (or a global field
of positive characteristic).

Automorphic forms are defined for 
a reductive algebraic group $G$ 
over a global field $F$ 
(either a finite extension of $\Q$ or of $\mathbb{F}_q(t)$ for 
a finite field $\mathbb{F}_q$ 
of $q$ elements).   They are functions on the adele points $G(\A_F)$ which are 
invariant by the translation by $G(F)$, satisfying certain conditions at some 
place(s) of $F$.    In practice, for characteristic 0 fields $F$, 
they are realized as functions on some quotient 
of a symmetric space, satisfying certain (real analytic) conditions.
For example, (elliptic) modular forms are automorphic forms for $\GL_2$
over the rationals $\Q$.   They are functions on the upper half plane 
(or, equally, on the quotient $\mathrm{SL}_2(\R)/\mathrm{SO}_2(\R)$),
and certain conditions at the real place amounts to the modular 
condition of some weight.    

In this article, we study certain automorphic forms for $\GL_d$ for $d \ge 1$
over a global field of positive characteristic.    We fix a place $\infty$
of $F$.     The automorphic forms that we consider 
are functions on $\GL_d(\A_F)$, which are invariant under 
$\GL_d(F)$, 
satisfying a condition at the place $\infty$. 
The slogan here is that we consider ``automorphic 
forms whose associated automorphic representation
is the Steinberg representation at $\infty$".  
This will be made precise later.
As an analogue of the symmetric space in positive characteristic,
there are the Bruhat-Tits buildings.   We use the Bruhat-Tits building
of $\PGL_d$ over $F_\infty$,  where $F_\infty$ is the local field (completion) 
at the place $\infty$.      
It is a simplicial complex of dimension $d-1$ whose 
set of simplices are quotients of $\PGL_d(F_\infty)$.   
For example, the set of zero simplices is isomorphic to 
$\GL_d(F_\infty)/\GL_d(O_\infty)$ where $O_\infty \subset F_\infty$
is the ring of integers, and the set of $(d-1)$-dimensional 
simplices is isomorphic to $\GL_d(F_\infty)/\cI$ 
where $\cI \subset \GL_d(\cO_\infty)$ 
is the Iwahori subgroup consisting of those matrices that are
congruent to an upper triangular matrix modulo the maximal ideal
of $\cO_\infty$.
There are 
many interpretations of the simplices (e.g. in terms of $\cO_\infty$-lattices, 
of norms, and, for quotients, of vector bundles over the proper 
smooth curve $C$ whose function field is $F$)
which  we we will use.

For the dictionary between our function field setup and 
the setup for modular forms, see the table 
in Section~\ref{sec:dictionary}.

\section{Classical modular symbols}
In the study of modular forms, 
especially of weight 2 
(automorphic forms 
for $\GL_2$ over $\Q$), 
one useful tool is modular symbols, as invented by Shimura 
and Eichler and developed by Manin
\cite{Manin1}.      
To introduce them, recall the following 
geometric setup for modular forms.
We consider some arithmetic subgroup 
$\Gamma \subset \mathrm{SL}_2(\Z)$.
Main examples include congruence subgroups
such as $\Gamma_1(N)$ for a positive integer $N$.   
Then modular forms 
appear as 1-forms on the analytic 
space $\Gamma \backslash \cH$ or on $\Gamma\backslash\cH^*$
where $\cH$ is the upper half space and 
$\cH^*=\cH \cup \mathbb{P}^1(\Q)$.
The quotients have algebraic models 
defined over some number field,
and the $\Gamma\backslash\cH^*$ is 
a smooth compactification of $\Gamma\backslash\cH$.
The set $\mathbb{P}^1(\Q)$ or 
the quotient $\Gamma\backslash\mathbb{P}^1(\Q)$
is called the set of cusps.
The Eichler-Shimura isomorphism states that
1-forms on $\Gamma\backslash\cH^*$
are exactly the cusp forms of weight 2
(automorphic forms with some condition
near the boundary in a compactification;
non-cusp forms are easier to study in that,
for one thing, 
they come from a smaller algebraic group).

A modular symbol in this geometric setup 
is a path 
from $a$ to $b$ in $\cH^*$ 
for cusps $a,b$, or its class 
in the relative homology
$H_1(\Gamma\backslash\cH^*, \{\text{cusps}\}; \C)$.
The relation with modular forms 
is given by integration, integrating a 
1-form (a modular form) 
on the path from $a$ to $b$.
We may regard modular symbols as elements
in the dual of cusp forms.
One main theorem is that 
the modular symbols generate the 
dual space of cusp forms.
This property enables us to 
study cusp forms using modular 
symbols, which are amenable to computation.

We refer to Manin's fairly recent introductory article 
\cite{Manin} for more information.   The reader is also referred to
\cite{FKSms} for applications in Iwasawa theory.

\section{Higher dimensional modular symbols of Ash and Rudolph}
\label{sec:AR ms intro}
Ash and Rudolph consider 
higher dimensional modular symbols in \cite{AR}.    
The automorphic forms they treat are for $\GL_d$
over the rationals $\Q$, and they are functions 
on the symmetric space 
$X=\mathrm{SL}_d(\R)/\mathrm{SO}_d(\R)$.    
They use the Borel-Serre 
bordification (a partial compacitification) $\overline{X}$
of $X$ (as a generalization of $\cH^*$).
Let $\Gamma$ be a torsion free arithmetic subgroup
and set 
$M=\Gamma\backslash X$ and $\overline{M}=\Gamma\backslash\overline{X}$.
Then their modular symbols are elements 
of the relative homology group 
$H_{d-1}(M, \overline{M}; \Z)$.
By Poincar\'e duality,
this group is dual to 
$H^N(M;\Z)=H^N(\Gamma;\Z)$
where $N=d(d-1)/2$,
which is closely related to the space 
of cusp forms.

To construct elements in the relative homology, 
they introduce universal modular symbols
as elements of the homology groups 
$H_{d-2}(T_d, \Z)$ 
of the Tits building $T_d$
(the simplicial complex associated with the poset of flags
in $\Q^d$).
(The Solomon-Tits theorem says that $T_d$ is 
homotopy equivalent to the bouquet of $(d-2)$-spheres.)
This $T_d$ is homotopy equivalent to 
$\overline{X} \setminus X$
and the projection map to the quotient
give elements in $H_{d-1}(M, \overline{M}; \Z)$.

Given an ordered
basis $q_1, \dots, q_d$ of $\Q^d$,
they construct a map from a $(d-2)$-sphere
to $T_d$.
They call the homology class of 
this sphere a universal modular symbol
and denote it 
by 
$[q_1, \dots, q_d]$.
It is known that universal modular symbols 
generate the homology group
$H_{d-2}(T_d, \Z)$.

Their main theorem 
is that the ordered bases
$q_1, \dots, q_d$ such that the matrix 
$(q_1, \dots, q_d)$ regarded as a matrix 
lying in $\mathrm{SL}_d$
(not $\mathrm{GL}_d$)
span the homology group
$H_{d-1}(\overline{M}, M;\Z)$
when $\Gamma$ is torsion free.
If $\Gamma$ is not torsion free,
there always exists a torsion free 
subgroup $\Gamma' \subset \Gamma$
of finite index.   Their result in this case
is that the homology group divided by
the space of modular symbols is killed
by the index $[\Gamma:\Gamma']$.

They also give 
an algorithm in expressing a general modular symbol
in terms of such ``unimodular" modular symbols.
We do not consider the analogue of this theorem.
One reason is that their proof make 
use of the fact that $\Z$ is a 
Euclidean domain.   The analogue of $\Z$, namely $A$, 
in our setup need not be a Euclidean domain.

\section{Our result on modular symbols}
\label{sec:result ms}
We take $F$ to be a global field of positive characteristic
and fix a place $\infty$ as above.

The space of automorphic forms we consider can be 
defined as follows.    
Let $C_\C=\Hom(\GL_d(F) \backslash \GL_d(\A), \C)$
be the space of $\C$-valued functions 
on the adele points that are invariant by the action
of $\GL_d(F)$.    
We let $C_\C^\infty=\bigcup_\bK C_\C^\bK$
where $\bK$ runs over the open compact subgroups
of $\GL_d(\A)$ where the superscript $()^\bK$ 
means the invariants.    This is the space 
of smooth vectors.
Let $\St_\C$ denote the Steinberg 
representation of $\GL_d(F_\infty)$
and $\cI \subset \GL_d(O_\infty)$
be the Iwahori subgroup.   
It is known that the Iwahori fixed part
$\St_{d,\C}^\cI$ is one dimensional.
Take a non-zero vector $v \in \St_{d,\C}^\cI$.   
Then we define 
\[
\cA_{\St_\C,\C}=\mathrm{Image}(\Hom_{\GL_d(F_\infty)}
(\St_{d,\C}, C_\C^\infty) \to
C_\C^\infty)
\]
where the arrow is the evaluation at $v$.
This is the space of automorphic forms 
with Steinberg at infinity, as mentioned in Section~\ref{sec:modauto},
and is the main object of our study.

We use the Bruhat-Tits building $\cBT_\bullet$
of $\PGL_d(F_\infty)$ as an analogue of 
the symmetric space $\mathrm{SL}_d(\R)/\mathrm{SO}_d(\R)$.
We set 
$X_{\bK, \bullet}
=\GL_d(F) \backslash \GL_d(\A^\infty) \times
\cBT_\bullet/\bK$
for the ring of finite adeles $\A^\infty$  and 
for a compact open subgroup $\bK \subset \GL_d(\A^\infty)$.
(This is the analogue of the double coset 
description of 
the $\C$-valued points of a Shimura variety.)
We show that there is a canonical isomorphism
\[
\varinjlim_\bK 
H_{d-1}^{\BM}
(X_{\bK, \bullet}, \C) 
\cong \cA_{\St_\C, \C}
\]
and study the geometry of $X_{\bK, \bullet}$.

Let us assume for simplicity that 
$X_{\bK, \bullet}=\Gamma \backslash \cBT_\bullet$ 
for some arithmetic subgroup $\Gamma \subset \GL_d(F)$.
(In general, it is a finite disjoint union of such quotients.   
This is analogous to that 
$\coprod_i \Gamma_i\backslash\cH$,
where $\Gamma_i$ are some arithmetic groups, is 
isomorphic to the $\C$-valued points of a modular curve.)

Let us give the definition of our modular symbols.
Our definition of modular symbols is analogous to
the paths in $\cH^*$.
By the definition of Bruhat-Tits buildings,
$\cBT_\bullet$ is the union of subsimplicial
complexes called apartments.   The apartments are 
indexed by the set of $F_\infty$-bases
of $F_\infty^{\oplus d}$.
Let $A_{q_1,\dots, q_d, \bullet} \subset \cBT_\bullet$ 
denote the apartement corresponding to 
the basis $q_1, \dots, q_d$.
We have a map
\[
A_{q_1, \dots, q_d, \bullet} \subset \cBT_\bullet
\to \Gamma \backslash \cBT_\bullet.
\]
For an arithmetic subgroup $\Gamma$
and an $F$-basis (that is, a basis of $F^{\oplus d}$
regarded as a basis of $F_\infty^{\oplus d}$),
this map is locally finite (in the sense that the inverse 
image of a simplex is a finite set).
The image of the 
fundamental class of the apartment 
$H_{d-1}^{\BM}(A_{q_1, \dots, q_d, \bullet}, \Z)$
in $H_{d-1}^\BM(\Gamma \backslash \cBT_\bullet, \Z)$
is well defined in this case, where $H^\BM$ means
the Borel-Moore homology.  Our modular symbols are 
defined to be the elements of this form as $q_1, \dots, q_d$ 
runs over all $F$-bases.  

Our modular symbols are a priori 
different from the modular symbols coming from
universal modular symbols (the analogue of those of Ash and Rudolph).
We prove that they actually coincide.

Our main theorem computes a bound of the index 
of 
the subgroup generated by 
modular symbols 
inside 
$H_{d-1}^\BM(\Gamma \backslash \cBT_\bullet, \Z)$.
We have a uniform bound which is independent of the choice of $\Gamma$.
The prime-to-$p$ (the characteristic of $F$)
part depends only on the base field, and 
divides the order of $\GL_d(\mathbb{F}_{q'})$ for 
some explicitly given $q'$.    The exponent of the $p$-part 
is given explicitly in terms of $d$.

\section{Outline of proof for universal modular symbols}
\label{sec:ums outline}
In Chapter~\ref{ch:pf for ums}, we prove our main theorem on 
modular symbols but for universal modular symbols.
We give an outline of the proof in this section.
\subsection{}
Let $\Gamma \subset \GL_d(F)$ be an arithmetic subgroup.
We construct (as described in Section~\ref{sec:result ms} 
above) modular 
symbols in the Borel-Moore homology
$H^\BM(\Gamma \backslash \cBT_\bullet, \Z)$
as the (fundamental) classes 
of the apartments corresponding to 
$F$-bases.
The goal is to describe the size of cokernel
of the (injective) map the ($\Z$-module of) 
modular symbols to the Borel-Moore homology
of $\Gamma \backslash \cBT_\bullet$.

To achieve this goal, we use the universal 
modular symbols of Ash and Rudolph.   That is,
we construct a map
\[
H_{d-2}(T_d) \to H_{d-1}^\BM(\Gamma \backslash \cBT_\bullet)
\]
from the homology of the Tits building $T_d$.
This homology group is the space generated by universal
modular symbols.
We compute the cokernel of this map.
Then we show that the image coincides with 
the space of (our) modular symbols.

\subsection{}
The analogue of the map 
\[
H_{d-2}(T_d) \to H_{d-1}^\BM(\Gamma \backslash \cBT_\bullet)
\]
appears in Ash and Rudolph as described briefly in Section~\ref{sec:AR ms intro}, 
however, the construction
is different.   We believe our method is simpler and might apply 
to their case as well.   There is a remark in Section~\ref{sec:compare AR}.

First, we express the Borel-Moore homology as an inverse limit:
\[
H_{d-1}^\BM(\Gamma \backslash \cBT_\bullet)
\cong \varprojlim_\alpha H_{d-1}(\Gamma \backslash \cBT_\bullet,
\Gamma \backslash \cBT_\bullet^{(\alpha)}).
\]
Here, $\alpha$ runs over 
positive real numbers,
and $\cBT_\bullet^{(\alpha)}$ is a subsimplicial complex
consisting of ``more unstable than $\alpha$" vector bundles.

\subsection{}
The following is the key sequence in our proof:
\[
\begin{array}{ll}
H_{d-2}(T_d)
\cong
H_{d-1}(\cBT_\bullet, \cBT_\bullet^{(\alpha)})
\twoheadrightarrow
H_0(\Gamma, H_{d-1}(\cBT_\bullet, \cBT_\bullet^{(\alpha)}))
\\
\cong
H_{d-1}^\Gamma(\cBT_\bullet, \cBT_\bullet^{(\alpha)})
\to 
H_{d-1}(\Gamma \backslash \cBT_\bullet,
\Gamma \backslash \cBT_\bullet^{(\alpha)}).
\end{array}
\]
The second map is the canonical surjection to 
coinvariants.
Let us describe the other three maps.

\subsection{}
The first map (isomorphism) is obtained by 
following the method of Grayson's (see \cite[Cor 4.2]{Gra}).
Let us explain this in this section.

Recall that each 0-simplex of $\Gamma \backslash \cBT_\bullet$
can be interpreted as a locally free sheaf (vector bundle)
on the curve $C$ whose function field is $F$.
It can be seen from the work of Grayson that 
the semi-stable ones lie in the ``middle", whereas
the unstable ones are closer to the ``boundary".
The picture to have in mind is in Serre's book 
(near \cite[p.106, II.2, Thm 9]{Trees}), 
where the quotient of the 
building of dimension 1 (tree) is discussed in detail.
Only the finite graph in the middle consist of semi-stable ones
and  the halflines, after a few steps,  consist of unstable ones only.

We define a subsimplicial complexes 
$\cBT_\bullet^{(\alpha)}$
by using the Harder-Narasimhan function which measures
how unstable a vector bundle is.    In particular,
all the simplices of $\cBT_\bullet^{(\alpha)}$ 
for sufficiently large $\alpha$
correspond to unstable ones.
(If $d=2$, $\cBT_\bullet^{(\alpha)}$
consists of the halflines which 
become shorter as $\alpha$
grows bigger.)

Recall on the other hand that a simplex of 
the Tits building $T_d$ corresponds to a flag in 
$F^d$.    Suppose that a $0$-simplex of $\cBT_\bullet^{(\alpha)}$
corresponds to an unstable vector bundle. 
Then there is a nontrivial Harder-Narasimhan filtration,
and by taking the generic fiber, we obtain a 
filtration or a flag of $F^{\oplus d}$.
This is how the two spaces $\cBT_\bullet^{(\alpha)}$
and $T_d$ are related.  
Using Grayson's method, we see that they are 
homotopy equivalent.   
Then, using that $\cBT_\bullet$ is contractible,
we obtain the first isomorphism.

\subsection{}
Let us look at the third map which is an isomorphism.
There is a Lyndon-Hochshild-Serre spectral sequence
(see Section~\ref{sec:LHS ss})
for a pair of spaces converging to the equivariant homology:
\[
E^2_{p,q}=
H_p(\Gamma, H_q(\cBT_\bullet, \cBT_\bullet^{(\alpha)}; \Z)) 
\Rightarrow 
H_{p+q}^\Gamma(\cBT_\bullet, \cBT_\bullet^{(\alpha)}; \Z).
\]
We can compute the $E_2$-page.
The Solomon-Tits theorem says that the 
homotopy type of $T_d$ 
is a bouquet of $(d-2)$-spheres.
This means that many of the $E_2$-terms are 
zero and 
the relevant terms are $E^2_{0,t}$
which are the homology groups of 
the relative space
$(\cBT_\bullet, \cBT_\bullet^{(\alpha)})$.
Hence we obtain the third map.

\subsection{}
For the fourth map, we use the following 
spectral sequence
\[
E_{p,q}^1=
\bigoplus_{\sigma \in \Sigma_p}
H_q(\Gamma_\sigma, \chi_\sigma)
\Rightarrow
H_{p+q}^\Gamma(\cBT_\bullet, \cBT_\bullet^{(\alpha)}; \Z)
\]
where $\Sigma_p$
is the set of $p$-simplices of 
$\cBT_\bullet \setminus \cBT_\bullet^{(\alpha)}$,
 $\Gamma_\sigma$ is the stabilizer group,
and $\chi_\sigma$ is the representation
associated with orientation.

Because it is a first quadrant spectral sequence
and the $E_1$-terms $E_{s,t}^1$ 
vanish for $t>d$,
we obtain the fourth map as the composition
of (injective) differential maps.
Thus the cokernel of the composition is 
equal to the cokernel of the fourth map.
To compute it, we bound the order of the $E^1$
terms, or the order of the stabilizer groups.

\subsection{Torsion in arithmetic subgroups}
It is common in characteristic zero case (see e.g. Borel-Serre, Ash-Rudolph)
to assume that the arithmetic subgroup is 
without torsion.   In that case, all the stabilizer groups 
are trival since finite.   Then the sought cokernel 
turns out to be trivial.   The general case is reduced
to the torsion free case because given an arithmetic subgroup
there always exists a torsion free 
arithmetic subgroup of finite index.    The exponent
of the cokernel 
in this case is killed by this index.

However, in positive characteristic, an arithmetic subgroup
always has a nontrivial $p$-torsion subgroup so we do not 
expect to do as well as in the characteristic zero case.
Instead, we will see that given an arithmetic subgroup $\Gamma$, 
there always exists 
an arithmetic subgroup $\Gamma' \subset \Gamma$
of finite index
which is $p'$-torsion free 
(each torsion element has order prime to $p$).

We can bound the sought cokernel in the $p'$-torsion free
case in terms of some power of $p$ depending only 
on the dimension $d$ 
by inspecting the shape of $p$-torsion subgroups
of an arithmetic subgroup.   
We can also bound the index $[\Gamma: \Gamma']$
in terms of $d$ and $q$.   We arrive at the bound of the 
cokernel in general in terms of $d$, $p$, and $q$.

\section{Outline for comparing modular symbols}
\label{sec:outline comparison}
After computing the cokernel for the image 
of the universal modular symbols,
there remains the task of comparing the image
with the space of our modular symbols.

The strategy is summarized in the following diagram (in Section~\ref{sec:CD comparison}):
\\
\begin{tikzcd}
H_{d-1}^\BM(\Gamma \backslash \cBT_\bullet)   \arrow[d,"\sim", "(1)"']                                                       &      H_{d-1}^\BM(A_\bullet)  \arrow[dd, "\sim"]  \arrow[l, "(10)"] 
\\
\varprojlim_\alpha H_{d-1}(\Gamma \backslash \cBT_\bullet, \Gamma\backslash \cBT_\bullet^{(\alpha)})    &    \\
\varprojlim_\alpha H_{d-1}(\cBT_\bullet, \cBT_\bullet^{(\alpha)})  \arrow[u, "(2)"']  \arrow[d, "\sim", "(3)"'] &   \varprojlim_\alpha H_{d-1}(A_\bullet, A_\bullet^{(\alpha)})  \arrow[l] \arrow[d, "\sim", "(6)"']
\\
H_{d-1}(\cBT_\bullet, \cBT_\bullet^{(\alpha)}) \arrow[d, "\sim", "(4)"']                                                   &    H_{d-1}(A_\bullet, A_\bullet^{(\alpha)})  \arrow[l]  \arrow[d, "\sim", "(7)"'] 
\\
H_{d-2}(\cBT_\bullet^{(\alpha)})   \arrow[d, "\sim", "(5)"']                                                                     &    H_{d-2}(A_\bullet^{(\alpha)})  \arrow[d, "\sim", "(8)"'],  \arrow[l]                      
\\
H_{d-2}(T_d)                                                                                                                                  &  H_{d-2}(|\cP'(B)|)                    \arrow[l, "(9)"]      
\\
                                                                                                                                                  &  H_{d-2}(S^{d-2})                     \arrow[u,"\sim"']
\end{tikzcd}

We first take an ordered basis $v_1,\dots,v_d$ of $F^{\oplus d}$.
We then have an embedding
$\varphi_{v_1, \dots, v_d}:
A_\bullet \to \cBT_\bullet$.
(The image of this was denoted $A_{v_1, \dots, v_d, \bullet}$ earlier.)
This gives the top arrow as the pushforward
by this embedding composed with the projection map
$\cBT_\bullet \to \Gamma \backslash \cBT_\bullet$.
We know (show) that the Borel-Moore homology
$H_{d-1}^\BM(A_\bullet)$
of an apartment is isomorphic to $\Z$
(the fundamental class is a generator).
The image of $1 \in \Z$ by this pushforward is the 
definition of our modular symbol.

The image of $1 \in \Z$ at the right bottom 
via the horizontal map
followed by the left vertical map 
is the definition of the Ash-Rudolph modular symbol
(or the modular symbol coming from the universal
modular symbol).

To prove that the two modular symbols coincide, 
we add the right vertical column and show 
that the diagram commutes.
The commutativity of the squares except for the 
last one is not difficult.

We follow Grayson \cite{Gra} for the 
construction of the map (5).
We need to look at the proof (which uses \cite{Quillen})
and see that it comes from a zigzag of morphisms of posets.
We interpret the map (8) in a similar manner
in terms of posets and 
then the commutativity of the last square follows.

\section{On the structure of the space of automorphic forms}
We study automorphic forms 
$\cA_{\St, \C}$ with Steinberg at infinity
and the intersection $\cA_{\St, \mathrm{cusp}, \C}$ 
with the cusp forms. 
  
Using known results, we verify (Proposition~\ref{prop:66_3})
that 
$\cA_{\St, \mathrm{cusp}, \C}$ is the direct sum of
irreducible cuspidal automorphic representations
that are cuspidal and whose local representation
at infinity is the Steinberg representation.
We also obtain (Theorem~\ref{7_PROP1}) 
a similar description for 
the space $\cA_{\St, \C}$.
While the space of $L^2$ automorphic forms
is well studied, 
not all subrepresentations of $\cA_{\St, \C}$ 
are $L^2$ and the $L^2$ methods do not apply.

Let us give some ideas and the outline of proof.
We first show that the space of automorphic 
forms is isomorphic to the Borel-Moore homology
of (a finite union of) the quotient $X_{\bK, \bullet}$ 
of 
the Bruhat-Tits building by an arithmetic subgroup.
We take the dual and work with the cohomology
with compact support.

Now, if there existed some good compactification $\overline{X}$
of $X=X_{\bK, \bullet}$
then we would be looking at an exact sequence
\[
H^{d-2}(\overline{X}\setminus X)
\to
H^{d-1}_c(X)
\to
H^{d-1}(\overline{X}).
\]
(The analogy in the case of modular forms is 
$X=Y_0(N), \overline{X}=X_0(N), \overline{X}\setminus X=\{cusps\}$,
and $d=2$.   Then $H^1(\overline{X})$ is (roughly) the space of cusp forms,
and the remaining task is to compute 
$H^0(\overline{X}\setminus X)$ as an 
automorphic representation.)
While the compactifications given in \cite{FKS} might be helpful,
we do not use them.

There are two steps.
The first step is to regard the cohomology as the limit 
as the boundary becomes smaller (again using the spaces 
$\cBT_\bullet^{(\alpha)}$).   
The second step is to express the limit as induced representation
using a covering spectral sequence.

\section{Remarks}
We give miscellaneous remarks.   

\subsection{Automorphic forms in positive characteristic}
We study automorphic forms in positive characteristic.
Because there are no archimedean primes, 
there is no difficulty coming from complex or real analysis.
In terms of automorphic forms, this means that 
all the smooth automorphic 
forms are admissible automorphic forms (we refer to Cogdell's Lectures 3,4 
in \cite{CKM} for the definitions).
Then the theory of automorphic representations is equivalent
to the (simpler) theory of representatations of Hecke algebras.
For example, decomposition into irreducible automorphic
representations is simpler in positive characteristic.

This also means that it is meaningful to consider the automorphic
forms with values in $\Z$, and this gives one natural $\Z$-structure 
($\Z$-submodule which spans the space of ($\C$-valued)
automorphic forms).
In the modular form case, there is the cohomology with coefficients
in $\Z$ of the 
topological space $\Gamma \backslash \cH$, which spans the space of modular
forms.   However, this does not have an interpretation as $\Z$-valued functions.

One simplification occurs.   The Borel-Moore homology, or the relative
homology, in the Ash-Rudolph case is not readily related to the 
space of automorphic forms, but its Poincar\'e dual is.
The pairing, in terms of automorphic forms, involves integration
and hence periods of automorphic forms.    For example, 
in the modular form case, the pairing involves values 
such as $\int_0^{i\infty} f(z)dz$ where $f$ is a modular form.

In our case, the Borel-Moore homology is equal to the space of 
automorphic forms, hence our results are stated directly in terms
of automorphic forms, with no reference to any pairing.    
However, in our application \cite{KY:Zeta elements}, 
we do take a pairing of modular symbols and cusp forms,
much analogous to the integration.    The second author does 
not know what the natural formulation is.

\subsection{}
The space $\cA_{\St, \C}$ 
of automorphic forms 
with Steinberg at infinity 
arises as the \'etale 
cohomology of Drinfeld modular
varieties.   Drinfeld modular varieties 
may be regarded as a Shimura variety 
for $\PGL_d$ over $F$.   (Shimura varieties
are defined only in characteristic zero.   Note
also that there is no Shimura variety for $\GL_d$
if $d \ge 3$.)
Thus by studying the cohomology,
Laumon \cite{Laumon1}, \cite{Laumon2}
obtains a result in 
the Langlands program for $GL_d$
over a global field in positive characteristic.

One drawback is that, via this method,
one can treat only those automorphic
representations whose local representation 
is the Steinberg representation at the fixed prime 
infinity.

\subsection{Our old preprint}
This article may be regarded as a revised version of 
our preprint \cite{KY:preprint}.   
The main result of loc.cit. 
is that the space of $\Q$-valued cusp forms with 
Steinberg at infinity 
is contained in the space generated by modular symbols
tensored with $\Q$.
The result in this article implies that.

We took a very different approach there.
In \cite{KY:preprint}, 
we used the Werner compactification
\cite{We1}, \cite{We2}
and used the duality twice.
One key observation there
was that the Werner compactification
was best suited 
when studying the group cohomology
of arithmetic subgroups (not the 
homology of the quotient spaces),
which in turn is closely related to the 
space of cusp forms.    We no longer 
use this observation here.

\section{Comparing with the Ash-Rudolph method}
\label{sec:compare AR}
Let us compare our method with that of Ash and Rudolph.

\subsection{}
Ash and Rudolph \cite[p.245]{AR} 
use the following sequence:
\[
H_{d-2}(T_d)
\cong H_{d-2}(\partial \overline X)
\cong H_{d-1}(\overline X, \partial \overline X)
\xto{\pi_*}
H_{d-1}(M, \partial M).
\]
Here 
$X=\SL_n(\R)/\SO_n(\R)$,
$\overline{X}$ is the Borel-Serre bordification,
$M=\Gamma \backslash \overline{X}$
is a manifold with boundary for some arithmetic subgroup $\Gamma$, 
and $M\setminus \partial M=\Gamma \backslash X$.
The map $\pi_*$ is the canonical projection, which is shown to be surjective
if $\Gamma$ is torsion free.   
In the proof of their result \cite[Prop 3.2]{AR}, however,
they use Poincar\'e duality and show
$H_c^N(\overline{X}) \to H_c^N(M)$ is 
injective where $N=d(d-1)/2$.   Then they translate these groups
into group cohomology of $\Gamma$
using that the action is free and the contractibility of spaces
$X$ and $\overline{X}$.
They apply the Borel-Serre duality (of group cohomology) 
to prove the claimed injectivity.

\subsection{}
Let us point out some differences with our method.

There is a difficulty in using a (partial) compactification of 
$\cBT_\bullet$.
There exists a compactification of Borel-Serre of $\cBT_\bullet$
(which is in the same spirit as the Borel-Serre bordification of $X$).
(There are also other compactifications.   See \cite{FKS} for 
overview.)
One difference is that an arithmetic subgroup does not act 
freely on (the boundary of) the compactification, 
even when restricted to some subgroup.    
Thus there is some difficulty in connecting the 
homology groups of spaces to group cohomology
(as was done by Ash and Rudolph).   
We also did not find the analogue 
of the Borel-Serre duality in the literature.   
There may also be some difficulty 
arising from the fact that
an arithmetic subgroup is usually not torsion free.

These considerations suggest that their method may
not apply directly to our case.   
On the other hand, we believe that our method applies
to the case treated by Ash and Rudolph.
Our method is more straightforward in that
we do not use the two dualities.

\section{Dictionary}
\label{sec:dictionary}
The following table is a dictionary between our function field 
setup (the right column) and the classical 
modular forms setup (the left column).   See also Sections~\ref{sec:modauto}.

\begin{center}
    \begin{tabular}{| p{2cm} | p{45mm} | p{45mm} |}
    \hline
base field         &         $\Q$                     &                       $F$ (a global field in positive characteristic; a function field of a curve over a finite field)   \\
\hline
place             & the real place $\infty$       &  a fixed place $\infty$   \\
\hline
integers        &    $\Z$                           &   $A$ (integral at all but $\infty$)  \\
\hline
completion     &    $\R$                           & $F_\infty$,  ($\cO_\infty$ the ring of integers)  \\ 
\hline
rank, dimension &                 $d=2$                    &                       $d\ge 1$            \\
\hline
algebraic group  &   $\mathrm{SL}_2$ over $\Q$                &                   $\PGL_d$ over $F$                                          \\
\hline
symmetric space &    $\mathfrak{H}=\mathrm{SL}_2(\R)/\mathrm{SO}_2(\R)$ the upper half plane  (real manifold)  &  $\cB\cT_\bullet$, $\cB\cT_0=\PGL_d(F_\infty)/\PGL_d(\cO_\infty)$  the Bruhat-Tits building (simplicial complex)       \\
\hline
modular symbol & geodesic from $a$ to $b$ ($a,b \in \mathbb{P}^1(\Q)$)   &   apartment corresponding to a basis of $F^{\oplus d}$  \\
\hline  
arithmetic subgroup  & $\Gamma \subset \mathrm{SL}_d(\Z)$           & $\Gamma \subset \GL_d(F)$  \\
\hline
automorphic form  &  modular form of weight 2 level $\Gamma$     &  $\Gamma$-invariant harmonic function (cochain) on $\cB\cT_{d-1}$  \\
\hline 
automorphic representation &  Of $\GL_2(\A_\Q)$, discrete series (for weight 2) at $\infty$  & Of $\GL_d(\A_F)$, Steinberg at $\infty$   \\
\\
\hline
homology  &   $H^{\BM}_1(\Gamma \backslash \cH, \Z)$   &   $H^\BM_{d-1}(\Gamma \backslash\cB\cT_\bullet, \Z)$     \\
\hline
cusp  &    $\mathbb{P}^1(\Q)$                 &   Several choices                    \\
\hline
(partial) compactification  &    $\cH \cup \mathbb{P}^1(\Q)$   & Several choices \\
\hline
    \end{tabular}
\end{center}

\chapter{Simplicial complexes and their (co)homology}

\label{sec:2}
The material of this chapter (except for the remark in Section~\ref{sec:cellular})
already appeared in Sections 3 and 5 of \cite{KY:Zeta elements}.
We collected them for the convenience of readers.

We define generalized simplicial complexes in Section~\ref{sec:simplicial complex},
define their 4 (co)homology theories (homology, cohomology, Borel-Moore
homology, cohomology with compact support), and 
mention the universal coefficient theorem and geometric realization.

Later, we will consider quotients of the Bruhat-Tits building by
arithmetic subgroups (to be defined).   While the Bruhat-Tits building is canonically a simplicial complex, the arithmetic quotient is in general not a simplicial 
complex.   This issue was addressed by Prasad in 
a paper of Harder \cite[p.140, Bemerkung]{Harder2}.   
We introduce this notion because 
the quotients are naturally (generlized) simplicial complexes.
An example of a generalized simplicial complex which is not a strict simplicial complex consists of two vertices with two edges between the vertices.

When defining the 4 (co)homology theories of simplicial complexes,
usually we fix an orientation of each simplex and then construct
suitable complexes of abelian groups and compute the (co)homology groups.
We give a slightly different complex,
where we do not fix a choice of orientation of each simplex.
This is because the Bruhat-Tits building is not canonically oriented
(even though they look like they are canonically oriented).

We end this chapter with a remark in Section~\ref{sec:cellular}.

\section{Generalized simplicial complexes}
\label{sec:simplicial complex}
\subsection{}
Let us recall the notion of (abstract) simplicial complex.
\begin{definition}
A (strict) simplicial complex is a pair $(Y_0,\Delta)$ of a set $Y_0$
and a set $\Delta$ of finite subsets of $Y_0$ which satisfies
the following conditions:
\begin{itemize}
\item If $S \in \Delta$ and $T\subset S$, then $T \in \Delta$.
\item If $v \in Y_0$, then $\{v \} \in \Delta$.
\end{itemize}
\end{definition}
In this paper we call a simplicial complex in the sense above
a strict simplicial complex, and use the terminology 
``simplicial complex" in a little broader sense as defined below.
We recall that in a strict (abstract) simplicial complex,
as recalled above, each simplex is uniquely determined
by the set of its vertices.

\subsection{}
The definition of (generalized) simplicial complex is as follows.
For a set $S$, let $\cP^\fin(S)$ denote the
category whose objects are the non-empty finite subsets of
$S$ and whose morphisms are the inclusions.
\begin{definition}
A simplicial complex is a pair $(Y_0,F)$ of a set $Y_0$
and a presheaf $F$ of sets on $\cP^\fin(Y_0)$
such that $F(\{\sigma \}) = \{\sigma\}$ holds for every
$\sigma \in Y_0$. 
\end{definition}

\subsection{}
The definition above is too abstract in practice.   We now give a working definition
which is equivalent to the definition above.
\begin{definition}
A  (generalized) simplicial complex is a collection 
$Y_\bullet = (Y_i)_{i \ge 0}$ of sets 
indexed by non-negative integers, equipped with
the following additional data:
\begin{itemize}
\item a subset $V(\sigma) \subset Y_0$ 
with cardinality $i+1$, for each $i \ge 0$ and
for each $\sigma \in Y_i$
(we call $V(\sigma)$ the set of vertices of $\sigma$),
and
\item an element in $Y_j$, for each $i \ge j \ge 0$, 
for each $\sigma \in Y_i$,
and for each subset $V' \subset V(\sigma)$ with cardinality
$j+1$ (we denote this element in $Y_j$ by the symbol
$\sigma \times_{V(\sigma)} V'$ and call it the face of
$\sigma$ corresponding to $V'$)
\end{itemize}
which satisfy the following conditions:
\begin{itemize}
\item For each $\sigma \in Y_0$, the equality
$V(\sigma) = \{\sigma\}$ holds,
\item For each $i \ge 0$, for each $\sigma \in Y_i$,
and for each non-empty subset $V' \subset V(\sigma)$,
the equality $V(\sigma \times_{V(\sigma)} V') = V'$ holds.
\item For each $i \ge 0$ and for each $\sigma \in Y_i$,
the equality $\sigma \times_{V(\sigma)} V(\sigma) = \sigma$
holds, and
\item For each $i \ge 0$, for each $\sigma \in Y_i$,
and for each non-empty subsets $V', V'' \subset V(\sigma)$
with $V'' \subset V'$, the equality 
$(\sigma \times_{V(\sigma)} V')\times_{V'} V'' = 
\sigma \times_{V(\sigma)} V''$ holds.
\end{itemize}
Let us call the element of the form 
$\sigma\times_{V(\sigma)} V'$ for $j$ and $V'$
as above, 
the $j$-dimensional face of $\sigma$ corresponding to $V'$.
We remark here that the symbol $\times_{V(\sigma)}$
does not mean a fiber product in any way.
\end{definition}

\subsection{}
This equivalence of the two definitions is explicitly
described as follows.  Suppose we are 
given a simplicial complex $Y_\bullet$.
Then the corresponding $F$ is the presheaf 
which associates, to a non-empty finite subset $V \subset Y_0$
with cardinality $i+1$, the set of elements $\sigma \in Y_i$
satisfying $V(\sigma)=V$.

This alternative definition of a simplicial complex is much 
closer to the definition of a strict simplicial complex.

\subsection{}
Any strict simplicial complex gives a simplicial 
complex in the sense above in the following way.
Let $(Y_0,\Delta)$ be a strict simplicial complex.
We identify $Y_0$ with the set of subsets of $Y_0$
with cardinality $1$.
For $i \ge 1$ let $Y_i$ denote the set of the elements
in $\Delta$ which has cardinality $i+1$ as a subset of $Y_0$.
For $i \ge 1$ and for $\sigma \in Y_i$, 
we set $V(\sigma)= \sigma$ regarded
as a subset of $Y_0$. 
%
For a non-empty subset $V \subset V(\sigma)$,
of cardinality $i'+1$, we set 
$\sigma \times_{V(\sigma)} V = V$ regarded
as an element of $Y_{i'}$. 
Then it is easily checked that 
the collection $Y_\bullet = (Y_i)_{i \ge 0}$ 
together with the assignments $\sigma \mapsto V(\sigma)$
and $(\sigma, V) \mapsto \sigma \times_{V(\sigma)} V$
forms a simplicial complex.

\subsection{}
An example of a (generalized) simplicial complex which is not 
a strict simplicial complex consists of two vertices with two edges 
between the two vertices.

\subsection{}
The Bruhat-Tits building (for PGL as introduced in Definition~\ref{def:BT}) 
is a strict simplicial complex.
However, its arithmetic quotients are generally merely a (generalized) simplicial complex.

\subsection{}
Let $(Y_0, F)$ and $(Z_0, G)$ be simplicial complexes.
\begin{definition}
A morphism $(Y_0, F) \to (Z_0, G)$ is a map of sets 
$f_0: Y_0 \to Z_0$ and a morphism of presheaves
$F \to G \circ P_{\mathrm{fin}}(f_0)$
on $P_{\mathrm{fin}}(f_0)$.   Here $P_{\mathrm{fin}}(f_0): P_{\mathrm{fin}}(Y_0) 
\to P_{\mathrm{fin}}(Z_0)$ is the functor induced by $f_0$.
\end{definition}

\subsection{}
The definition above in terms of the working definition
is as follows.
Let $Y_\bullet$ and $Z_\bullet$ be simplicial complexes.
We define a map (morphism) from $Y_\bullet$ to $Z_\bullet$ to
be a collection 
$f=(f_i)_{i \ge 0}$ of maps $f_i : Y_i \to Z_i$ of sets 
which satisfies the following
conditions:
\begin{itemize}
\item for any
$i \ge 0 $ and for any $\sigma \in Y_i$, 
the restriction of $f_0$ to $V(\sigma)$ is injective
and the image of $f|_{V(\sigma)}$ is equal to the set
$V(f_i(\sigma))$, and 
\item for any
$i \ge j \ge 0$, for any $\sigma \in Y_i$, and for any
non-empty subset $V' \subset V(\sigma)$ with cardinality $j+1$
we have $f_j(\sigma \times_{V(\sigma)} V')
= f_i(\sigma) \times_{V(f_i(\sigma))} f_0(V')$.
\end{itemize}

\section{Orientation}
Usually the (co)homology groups of $Y_\bullet$ are defined 
to be the (co)homology groups of a complex $C_\bullet$ 
whose component in degree $i$ is the free abelian group 
generated by the $i$-simplices of $Y_\bullet$. 
For a precise definition of the boundary homomorphism
of the complex $C_\bullet$, 
we need to choose an orientation of each simplex. 
In this paper we adopt an alternative, 
equivalent definition of homology groups which 
does not require any choice of orientations. 
The latter definition seems a little complicated at first glance, 
however it will soon turn out to be a better way 
for describing the (co)homology of the arithmetic quotients 
Bruhat-Tits building, which seems to have no canonical, 
good choice of orientations.

\subsection{}
\label{sec:orientation}
We introduce the notion of orientation of a simplex.
It will be a $\{\pm 1\}$-torsor (where $\{\pm 1\}$
is the abelian group $\Z/2\Z$)
associated with each simplex.

Let $Y_\bullet$ be a simplicial complex
and let $i \ge 0$ be a non-negative integer.
For an $i$-simplex $\sigma \in Y_i$,
we let $T(\sigma)$ denote the set of all bijections
from the finite set $\{1,\ldots, i+1 \}$ of cardinality
$i+1$ to the set $V(\sigma)$ of vertices of $\sigma$.
The symmetric group $S_{i+1}$ acts on the set 
$\{1,\ldots, i+1 \}$ from the left 
and hence on the set $T(\sigma)$ from the right.
Through this action the set $T(\sigma)$ is
a right $S_{i+1}$-torsor.

\begin{definition}
We define the set $O(\sigma)$ of orientations of $\sigma$
to be the $\pmone$-torsor 
$O(\sigma) = T(\sigma) \times_{S_{i+1},\sgn} \pmone$ 
which is the push-forward 
of the $S_{i+1}$-torsor $T(\sigma)$ 
with respect to the signature character 
$\sgn: S_{i+1} \to \pmone$.
\end{definition}

We give the basic properties which we 
need in order to consider (co)homology.

\subsection{}
When $i \ge 1$, the $\pmone$-torsor $O(\sigma)$ is
isomorphic, as a set, to the quotient 
$T(\sigma)/A_{i+1}$ of $T(\sigma)$
by the action of the alternating group 
$A_{i+1} = \Ker\, \sgn \subset S_{i+1}$. 
When $i=0$, the $\pmone$-torsor $O(\sigma)$ 
is isomorphic to the product
$O(\sigma) = T(\sigma) \times \pmone$, on which the group $\pmone$
acts via its natural action on the second factor.

\subsection{}
Let $i \ge 1$ and let $\sigma \in Y_i$.
For $v \in V(\sigma)$ let $\sigma_v$ denote the
$(i-1)$-simplex 
$\sigma_v 
= \sigma \times_{V(\sigma)} (V(\sigma) \setminus \{v\})$.
There is a canonical map
$s_v : O(\sigma) \to O(\sigma_v)$ of 
$\pmone$-torsors defined as follows.
Let $\nu \in O(\sigma)$ and take a lift
$\wt{\nu}:\{1,\ldots,i+1\} \xto{\cong} V(\sigma)$
of $\nu$ in $T(\sigma)$. Let 
$\wt{\iota}_v : \{1,\ldots,i\} \inj \{1,\ldots,i+1\}$
denote the unique order-preserving injection
whose image is equal to $\{1,\ldots,i+1\} \setminus 
\{\wt{\nu}^{-1}(v)\}$. It follows from the
definition of $\wt{\iota}_v$ that the composite
$\wt{\nu} \circ \wt{\iota}_v: \{1,\ldots,i\} 
\to V(\sigma)$ induces a bijection
$\wt{\nu}_v : \{1,\ldots,i\}
\xto{\cong} V(\sigma) \setminus \{v\}
= V(\sigma_v)$. We regard $\wt{\nu}_v$
as an element in $T(\sigma_v)$. 
We define $s_v : O(\sigma) \to O(\sigma_v)$
to be the map which sends $\nu \in O(\sigma)$
to $(-1)^{\wt{\nu}^{-1}(v)}$ times
the class of $\wt{\nu}_v$. It is easy to check that
the map $s_v$ is well-defined.

\subsection{}
Let $i \ge 2$ and $\sigma \in Y_i$.
Let $v, v' \in V(\sigma)$ with $v \neq v'$.
We have $(\sigma_v)_{v'} = (\sigma_{v'})_v$.
Let us consider the two composites
$s_{v'} \circ s_v :
O(\sigma) \to O((\sigma_v)_{v'})$ and
$s_v \circ s_{v'} :
O(\sigma) \to O((\sigma_{v'})_v)$.
It is easy to check that the equality
\begin{equation} \label{formula1}
s_{v'} \circ s_v (\nu) 
= (-1) \cdot s_v \circ s_{v'} (\nu)
\end{equation}
holds for every $\nu \in O(\sigma)$.

\section{Cohomology and homology}
\label{sec:def homology}

\begin{definition}
We say that a simplicial complex
$Y_\bullet$ is locally finite if for any $i \ge 0$
and for any $\tau \in Y_i$, there exist only
finitely many $\sigma \in Y_{i+1}$ such that
$\tau$ is a face of $\sigma$.
\end{definition}
We give the four notions of homology or cohomology
for a locally finite (generalized) simplicial complex.

\subsection{}
Let $Y_\bullet$ be a simplicial complex 
(\resp a locally finite simplicial complex).
For an integer $i\ge 0$, we set 
$Y_i'=\coprod_{\sigma \in Y_i} O(\sigma)$ 
and regard it as a $\pmone$-set.
Given an abelian group $M$, 
we regard the abelian groups
$\bigoplus_{\nu \in Y_i'} M$
and $\prod_{\nu \in Y_i'}M$
as $\pmone$-modules in such a way that
the component at $\eps\cdot \nu$ of
$\eps \cdot (m_\nu)$ is equal to 
$\eps m_\nu$ for $\eps \in \pmone$
and for $\nu \in Y_i'$.

\subsection{}
For $i \ge 1$, we let 
$\wt{\partial}_{i,\oplus}: \bigoplus_{\nu \in Y_i'}M
\to \bigoplus_{\nu \in Y_{i-1}'}M$ 
(\resp $\wt{\partial}_{i,\prod}: \prod_{\nu \in Y_i'}M
\to \prod_{\nu \in Y_{i-1}'}M$)
denote the homomorphism of abelian groups 
that sends
$m = (m_\nu)_{\nu \in Y_i'}$ to the element
$\wt{\partial}_i(m)$ whose coordinate 
at $\nu' \in O(\sigma') \subset Y_{i-1}'$ is equal to
\begin{eqnarray}\label{eqn:boundary}
\wt{\partial}_i(m)_{\nu'} =
\sum_{(v,\sigma,\nu)} m_\nu.
\end{eqnarray}
Here in the sum on the right hand side
$(v,\sigma,\nu)$ runs over the triples
of $v \in Y_0 \setminus V(\sigma')$,
an element $\sigma \in Y_i$, and $\nu \in O(\sigma)$
which satisfy $V(\sigma) = V(\sigma') \amalg \{v\}$ and
$s_v(\nu) = \nu'$.

Note that the sum on 
the right hand side is a finite sum
for $\wt{\partial}_{i,\oplus}$ by definition.
One can see also that the sum is a finite sum
in the case of $\wt{\partial}_{i,\prod}$ 
using the locally finiteness of $Y_\bullet$.

Each of $\wt{\partial}_{i,\oplus}$ 
and $\wt{\partial}_{i,\prod}$ is a
homomorphism of $\pmone$-modules. 
Hence it induces a homomorphism
$\partial_{i,\oplus} : (\bigoplus_{\nu \in Y_i'}M)_\pmone
\to (\bigoplus_{\nu \in Y_{i-1}'}M)_\pmone$
(\resp $\partial_{i,\prod} : (\prod_{\nu \in Y_i'}M)_\pmone
\to (\prod_{\nu \in Y_{i-1}'}M)_\pmone$) of abelian groups,
where the subscript $\pmone$ means the coinvariants.

\subsection{}
It follows from the formula (\ref{formula1}) 
and the definition
of $\partial_{i,\oplus}$ and $\partial_{i,\prod}$ that
the family of abelian groups
$((\bigoplus_{\nu \in Y_i'}M)_\pmone)_{i\ge 0}$
(resp.
$((\prod_{\nu \in Y_i'}M)_\pmone)_{i\ge 0}$)
indexed by the integer $i \ge 0$, together with
the homomorphisms $\partial_{i,\oplus}$ 
(resp. $\partial_{i,\prod}$)
for $i \ge 1$,
forms a complex of abelian groups.
The homology groups of this complex are 
the homology groups $H_*(Y_\bullet, M)$
(resp. the Borel-Moore homology groups
$H_*^\BM(Y_\bullet, M)$)
of $Y_\bullet$ 
with coefficients in $M$.

\subsection{}
We note that there is a canonical map
$H_*(Y_\bullet, M) \to H_*^\BM(Y_\bullet, M)$
from homology to Borel-Moore homology
induced by the map 
of complexes
$((\bigoplus_{\nu \in Y_i'}M)_\pmone)_{i\ge 0}
\to ((\prod_{\nu \in Y_i'}M)_\pmone)_{i\ge 0}$
given by inclusion at each degree.

\subsection{}
The family of abelian groups 
$(\Map_{\pmone}(Y_i', M))_{i\ge 0}$
(resp.
$(\Map_{\mathrm{fin}, {\pmone}}(Y_i', M))_{i\ge 0}$
where the subscript $\mathrm{fin}$ means finite support)
of the $\pmone$-equivariant maps from $Y_i'$ to $M$
forms a complex of
abelian groups in a similar manner.
(One uses the locally finiteness of $Y_\bullet$
for the latter.)  
The cohomology groups of this complex
are the cohomology groups $H^*(Y_\bullet, M)$
(resp. the cohomology groups with compact support 
$H_c^*(Y_\bullet, M)$)
of $Y_\bullet$ 
with coefficients in $M$.
There is a canonical map from cohomology with compact support to 
cohomology.

\section{Universal coefficient theorem}
\label{univ_coeff}
It follows from the definitions that
the following universal coefficient
theorem holds.  
\subsection{}
For a simplicial complex $Y_\bullet$, 
there exist canonical short exact sequences
$$
0 \to H_i(Y_\bullet, \Z) \otimes M
\to H_i(Y_\bullet, M) \to 
\Tor_1^\Z (H_{i-1}(Y_\bullet, \Z),M) \to 0
$$
and
$$
0 \to \Ext^1_\Z(H_{i-1}(Y_\bullet, \Z),M)
\to H^i(Y_\bullet, M) \to 
\Hom_\Z (H_i(Y_\bullet, \Z),M) \to 0.
$$
for any abelian group $M$. 

\subsection{}
Similarly, for a locally finite simplicial complex $Y_\bullet$,
we have short exact sequences
$$
0 \to \Ext^1_\Z(H_c^{i+1}(Y_\bullet, \Z),M)
\to H^\BM_i(Y_\bullet,M)
\to \Hom_\Z (H_c^i(Y_\bullet, \Z),M) \to 0
$$
and
$$
0 \to H^i_c (Y_\bullet, \Z) \otimes M
\to H_c^i(Y_\bullet, M) \to 
\Tor_1^\Z (H_c^{i+1}(Y_\bullet, \Z),M) \to 0
$$
for any abelian group $M$. 
%

%
\section{Some induced maps}
\label{sec:quasifinite}
Let $f=(f_i)_{i \ge 0}:
Y_\bullet \to Z_\bullet$ be a map of
simplicial complexes. For each integer $i \ge 0$
and for each abelian group $M$,
the map $f$ canonically 
induces homomorphisms $f_* : H_i(Y_\bullet, M)
\to H_i(Z_\bullet,M)$ and
$f^* : H^i(Z_\bullet,M) \to H^i(Y_\bullet,M)$.
We say that the map $f$ is finite
if the subset $f_i^{-1}(\sigma)$ of $Y_i$ is
finite for any $i \ge 0$ and for any 
$\sigma \in Z_i$. If $Y_\bullet$ and $Z_\bullet$ are
locally finite, and if $f$ is finite, then
$f$ canonically 
induces the 
pushforward homomorphism $f_* : H^\BM_i(Y_\bullet, M)
\to H^\BM_i(Z_\bullet,M)$ and
the pullback homomorphism
$f^*: H_c^i(Z_\bullet,M) \to H_c^i(Y_\bullet,M)$.

\section{Geometric realization} \label{sec:CWcomplex}
 Let $Y_\bullet$ be a simplicial complex.
We associate a CW complex $|Y_\bullet|$ which we call
the geometric realization of $Y_\bullet$.

Let $I(Y_\bullet)$ denote the disjoint
union $\coprod_{i \ge 0} Y_i$. We define a partial order
on the set $I(Y_\bullet)$ by saying that $\tau \le \sigma$
if and only if $\tau$ is a face of $\sigma$.
For $\sigma \in I(Y_\bullet)$, we let $\Delta_\sigma$ denote
the set of maps $f:V(\sigma) \to \R_{\ge 0}$ satisfying
$\sum_{v \in V(\sigma)} f(v) =1$. We regard $\Delta_\sigma$
as a topological space whose topology is induced from that
of the real vector space $\Map(V(\sigma),\R)$.
If $\tau$ is a face of $\sigma$, we regard the space
$\Delta_\tau$ as the closed subspace of $\Delta_\sigma$
which consists of the maps $V(\sigma) \to \R_{\ge 0}$
whose support is contained in the subset 
$V(\tau) \subset V(\sigma)$.
We let $|Y_\bullet|$ denote the colimit
$\varinjlim_{\sigma \in I(Y_\bullet)} \Delta_\sigma$
in the category of topological spaces
and call it the geometric realization of $Y_\bullet$.
%
It follows from the definition that
the geometric realization $|Y_\bullet|$ has a canonical
structure of CW complex.

\section{Cellular versus singular}
\label{sec:cellular}
We give a remark on the use of the term 
``Borel-Moore homology''
in this paragraph.
Given a strict simplicial complex,
its cohomology, homology and cohomology with compact support
(for a locally finite strict simplicial complex)
are usually defined as above,
and called cellular (co)homology.  See for example \cite{Hatcher}.

On the other hand, there is also the singular (co)homology
and with compact support that are 
defined using the singular (co)chain complex.  
It is well-known that the cellular (co)homology groups
(with compact support) are isomorphic to the 
singular (co)homology groups (with compact support)
of the geometric realization.
The same proof applies to the generalized simplicial 
complexes and gives an isomorphism between the
cellular and the singular theories.

For the Borel-Moore homology, we did not find a cellular 
definition as above, except in Hattori \cite{Hattori} 
where he does
not call it the Borel-Moore homology.  He also gives a
definition using singular chains and shows that the two homology
groups are isomorphic.

There are several definitions of Borel-Moore homology
that may be associated to a (strict) simplicial complex.
The definition of the Borel-Moore homology for PL manifolds 
is found in Haefliger \cite{Haefliger}.  
There is also a sheaf theoretic definition in Iversen 
\cite{Iversen}.
More importantly, there is the general definition which concerns
the intersection homology.
However, we did not find a statement in the literature
and we did not check that the cellular
definition in Hattori (same as the one given in this article) is
isomorphic to the other Borel-Moore homology theories.
\newcommand{\be}{\mathbf{e}}
\chapter{The Bruhat-Tits building and apartments}
\label{sec:BT}
Let $d\ge 1$ be a positive integer.
In this chapter, we recall the definition of the 
Bruhat-Tits building of $\mathrm{PGL}_d$
over a nonarchimedean local field using 
lattices and subsimplicial complexes
called apartments.
(For the general theory of Bruhat-Tits building and
apartments, the reader is referred to \cite{BT} and 
the book \cite{Ab-Br}.)
Then we define the fundamental class of an apartment
in its $(d-2)$-nd Borel-Moore homology group.

Later, we define modular symbols to be the 
image of the fundamental classes of apartments 
associated with $F$-basis ($F$ is a global field)
in the Borel-Moore homology of quotients 
of the Bruhat-Tits building.

\section{The Bruhat-Tits building of $\mathrm{PGL}_d$}
In the following paragraphs, 
we recall the definition of
the Bruhat-Tits building of $\PGL_d$ over a 
nonarchimedean local field.
We recall that it is a strict simplicial complex.

\subsection{Notation}
Let $K$ be a nonarchimedean local field.
We let $\cO \subset K$ denote the ring of integers.
We fix a uniformizer $\varpi \in \cO$.
Let $d \ge 1$ be an integer.
Let $V=K^{\oplus d}$.  We regard it as the set of row vectors
so that $\GL_d(K)$ acts from the right by multiplication.

\subsection{the Bruhat-Tits building (\cite{BT}) using lattices}
\label{Bruhat-Tits}
We do not recall here the most general definition 
of the Bruhat-Tits buildings.   Let us give a definition
using lattices (see also \cite[\S 4]{Gra}) first, 
and then give a more explicit description 
for later use.
\subsubsection{}
An $\cO$-lattice in $V$ is 
a free $\cO$-submodule of $V$ of rank $d$.
We denote by $\Lat_{\cO}(V)$ the
set of $\cO$-lattices in $V$.
We regard the set $\Lat_{\cO}(V)$ as
a partially ordered set whose elements are ordered 
by the inclusions of $\cO$-lattices.
%
\subsubsection{}
Two $\cO$-lattices $L$, $L'$ of $V$
are called homothetic if $L = \varpi^j L'$
for some $j \in \Z$. Let $\Latbar_{\cO}(V)$
denote the set of homothety classes of 
$\cO$-lattices $V$. 

We denote by
$\cl$ the canonical surjection
$\cl: \Lat_{\cO}(V) \to
\Latbar_{\cO}(V)$.

\begin{definition}
We say that a subset $S$ of 
$\Latbar_{\cO}(V)$
is totally ordered if $\cl^{-1}(S)$ is 
a totally ordered subset of $\Lat_{\cO}(V)$.
\end{definition}

\subsubsection{}
The pair $(\Latbar_{\cO}(V), \Delta)$ 
of the set $\Latbar_{\cO}(V)$ 
and the set $\Delta$ of totally ordered finite nonempty subsets of
$\Latbar_{\cO}(V)$ forms
a strict simplicial complex.
\begin{definition}
\label{def:BT}
The Bruhat-Tits building of $\PGL_d$ over $K$
is the simplicial complex $\cBT_\bullet$ 
associated to this strict simplicial complex.
\end{definition}

\subsection{Explicit description of the building}
In the next paragraphs we explicitly describe
the simplicial complex $\cBT_\bullet$.
\subsubsection{}
For an integer $i \ge 0$, let
$\wt{\cBT}_i$ be the set of sequences $(L_j)_{j \in \Z}$
of $\cO$-lattices in $V$ indexed by $j \in \Z$
such that $L_j \supsetneqq L_{j+1}$ 
and $\varpi L_j=L_{j+i+1}$ hold for all $j\in \Z$.
In particular, if $(L_j)_{j \in \Z}$ is an element in 
$\wt{\cBT}_0$, then $L_j = \varpi^j L_0$ for
all $j \in \Z$. We identify the set
$\wt{\cBT}_0$ with the set $\Lat_{\cO}(V)$
by associating the $\cO$-lattice $L_0$
to an element $(L_j)_{j \in \Z}$ in $\wt{\cBT}_0$.
We say that two elements $(L_j)_{j \in \Z}$
and $(L'_j)_{j \in \Z}$ in $\wt{\cBT}_i$ are equivalent if
there exists an integer $\ell$ satisfying
$L'_j=L_{j+\ell}$ for all $j \in \Z$.
We will see below that the set of the equivalence classes in $\wt{\cBT}_i$
is identified with $\cBT_i$. 
For $i=0$, the identification $\wt{\cBT_0}
\cong \Lat_{\cO}(V)$ gives an identification
$\cBT_0 \cong \Latbar_{\cO}(V)$.

\subsubsection{}
Let $\sigma \in \cBT_i$ and take a representative
$(L_j)_{j \in \Z}$ of $\sigma$.
For $j \in \Z$, let us consider the class
$\cl(L_j)$ in $\Latbar_{\cO}(V)$.
Since $\varpi L_j = L_{j+i+1}$, we have
$\cl(L_j) = \cl(L_{j+i+1})$. Since
$L_j \supsetneqq L_k \supsetneqq \varpi L_j$
for $0 \le j < k \le i$, the elements
$\cl(L_0), \ldots, \cl(L_i) \in 
\Latbar_{\cO}(V)$
are distinct. Hence the subset 
$V(\sigma) = \{ \cl(L_j)\ |\ j \in \Z \} \subset \cBT_0$ 
has cardinality $i+1$ and
does not depend on the choice of $(L_j)_{j \in \Z}$.
It is easy to check that the map from
$\cBT_i$ to the set of finite subsets of 
$\Latbar_{\cO}(V)$ which sends
$\sigma \in \cBT_i$ to $V(\sigma)$ is injective
and that
the set $\{V(\sigma)\ |\ \sigma \in \cBT_i\}$
is equal to the set of totally ordered subsets of
$\Latbar_{\cO}(V)$ with cardinality
$i+1$. In particular, for any $j \in \{0,\ldots,i\}$ 
and for any subset $V' \subset V(\sigma)$
of cardinality $j+1$, there exists a unique element
in $\cBT_j$, which we denote by $\sigma \times_{V(\sigma)} V'$,
such that $V(\sigma \times_{V(\sigma)} V')$ is equal to $V'$.
Thus the collection $\cBT_\bullet = \coprod_{i \ge 0} \cBT_i$
together with the data $V(\sigma)$ and
$\sigma \times_{V(\sigma)} V'$
forms a simplicial complex which is canonically isomorphic
to the simplicial complex associated to the
strict simplicial complex 
$(\Latbar_{\cO}(V), \Delta)$ 
which we introduced in the
first paragraph of Section~\ref{Bruhat-Tits}.
We call the simplicial complex $\cBT_\bullet$
the Bruhat-Tits building of $\PGL_d$ over $K$.

\subsubsection{}
\label{sec:dimension}
The simplicial complex $\cBT_\bullet$ is of dimension
at most $d-1$, by which we mean that $\cBT_i$ is an empty
set for $i > d-1$. It follows from the fact that
$\wt{\cBT}_i$ is an empty set for $i > d-1$, which
we can check as follows. Let $i > d-1$ and assume 
that there exists an element 
$(L_j)_{j \in \Z}$ in $\wt{\cBT}_i$.
Then for $j=0,\ldots,i+1$, the quotient 
$L_j/L_{i+1}$ is a subspace of the $d$-dimensional 
$(\cO/\varpi \cO)$-vector 
space $L_0/L_{i+1}=L_0/\varpi L_0$. These subspaces
must satisfy $L_0/L_{i+1} \supsetneqq L_1/L_{i+1}
\supsetneqq \cdots \supsetneqq L_{i+1}/L_{i+1}$.
It is impossible since $i+1 > d$.

\section{Apartments}
\label{sec:apartments}

Here we recall the definition of apartment.
It is a simplicial subcomplex of the Bruhat-Tits 
building.   
We are interested only in the apartements of 
$\PGL_d$ of a nonarchimedean local field
and not of other algebraic groups.
To describe apartments, we do not use the general theory 
via root systems but give a simpler, ad hoc treatment,
particular to $\PGL_d$.
The readers are referred to \cite[p. 523, 10.1.7 Example]{Ab-Br}
for the general theory.

For example, when $d=2$, it is an infinite sequence of $1$-simplices.
When $d=3$, it is an $\mathbb{R}^{2}$ tiled by triangles ($2$-simplices).
The geometric realization is homeomorphic to $\mathbb{R}^{d-1}$.

\subsection{}
\label{sec:def apartment}
Set $A_0 = \Z^{\oplus d}/\Z(1,\ldots,1)$.
For two elements $x=(x_j), y=(y_j) \in \Z^{\oplus d}$,
we write $x \le y$ when $x_j \le y_j$ for all $1 \le j \le d$.
We regard $\Z^{\oplus d}$ as a partially ordered set with
respect to $\le$.
Let $\pi: \Z^{\oplus d} \to A_0$ denote the quotient map.

We say that a subset $S$ of $A_0$
is totally ordered if $\pi^{-1}(S)$ is 
a totally ordered subset of $\Z^{\oplus d}$.
The pair $(A_0, D)$ 
of the set $A_0$ 
and the set $D$ of totally ordered finite nonempty subsets of 
$A_0$ forms a strict simplicial complex.
We denote by $A_\bullet =(A_i)_{i \ge 0}$
the simplicial complex associated to the
strict simplicial complex $(A_0, \coprod_{i\ge 0} A_i)$.
We note that $A_i$ is an empty set for $i \ge d$, 
since by definition there is no totally ordered subset of $A_0$
with cardinality larger than $d$.

\subsection{} \label{sec:521}
Let $v_1, \dots, v_d$ be a basis of $V =K^{\oplus d}$.
We define a map $\iota_{v_1,\ldots,v_d} : 
A_\bullet \to \cBT_\bullet$ of simplicial complexes
as follows.

Let $\wt{\iota}_{v_1,\ldots,v_d}: \Z^{\oplus d} \to \wt{\cBT}_0$ 
denote the map that sends the element $(n_1,\ldots,n_d) \in \Z^d$
to the $\cO$-lattice 
$\cO\varpi^{n_1} v_1  
\oplus 
\cO \varpi^{n_2} v_2  
\oplus \dots \oplus
\cO \varpi^{n_d} v_d$.
The map $\wt{\iota}_{v_1,\ldots,v_d}$ is an order-embedding of partially
ordered sets, and induces a map
$\iota_{v_1,\ldots,v_d,0} : A_0 \to \cBT_0$
that makes the diagram
$$
\begin{CD}
\Z^{\oplus d}
@>{\wt{\iota}_{v_1,\ldots,v_d}}>>
\wt{\cBT}_0 \\
@VVV @VVV \\
A_0
@>{\iota_{v_1,\ldots,v_d,0}}>>
\cBT_0
\end{CD}
$$
of sets, where the vertical arrows are the quotient maps,
commutative and cartesian.
This implies that the map $\iota_{v_1,\ldots,v_d,0} : A_0 \to \cBT_0$
sends a totally ordered subset of $A_0$ to a totally
ordered subset of $\cBT_0$.
Hence the map $\iota_{v_1,\ldots,v_d,0} : A_0 \to \cBT_0$
induces a map $\iota_{v_1,\ldots,v_d} : A_\bullet \to \cBT_\bullet$ of
simplicial complexes.

It is easy to check that the map
$\iota_{v_1,\ldots,v_d,i} : A_i \to \cBT_i$
is injective for every $i \ge 0$.
We define a simplicial subcomplex 
$A_{v_1, \dots, v_d , \bullet}$ 
of $\cBT_\bullet$ to be the image of the
map $\iota_{v_1,\ldots,v_d}$ so that
$A_{v_1, \dots, v_d , i}$ is the image of
the map $\iota_{v_1,\ldots,v_d,i}$ for each $i \ge 0$.
We call the subcomplex 
$A_{v_1, \dots, v_d , \bullet}$ of $\cBT_\bullet$
the apartment in $\cBT_\bullet$ corresponding to the basis
$v_1,\ldots,v_d$. Since the map
$\iota_{v_1,\ldots,v_d,i}$ is injective for every $i \ge 0$,
the map $\iota_{v_1,\ldots,v_d}$ induces an isomorphism
$A_\bullet \xto{\cong} A_{v_1, \dots, v_d , \bullet}$
of simplicial complexes.

\section{the fundamental class} \label{sec:fundamental class}
We introduce a special element $\beta$ in
the group $H^\BM_{d-1}(A_\bullet,\Z)$, which is
an analogue of the fundamental class.

\subsection{ }
Let $\sigma \in A_i$ and let $V(\sigma) \subset A_0$ denote the
set of vertices of $\sigma$.
By definition, $V(\sigma)$ consists of exactly $i+1$ elements 
and the inverse image $\wt{V}(\sigma) = \pi^{-1}(V(\sigma))$ 
is a totally ordered subset of $\Z^d$.

\begin{lem} \label{lem:totally}
As a totally ordered set, $\wt{V}(\sigma)$ is 
isomorphic to $\Z$.
\end{lem}

\begin{proof}
Let $\Sigma: \Z^{\oplus d} \to \Z$ denote the map
that sends $(n_1,\ldots,n_d)$ to $n_1 + \cdots + n_d$.
Then the composite of the inclusion $\wt{V}(\sigma) \inj \Z^{\oplus d}$
with $\Sigma$ is order-preserving and injective since
$\wt{V}(\sigma)$ is totally ordered and for any
$x,y \in \Z^{\oplus d}$, the relations $x \le y$ and
$\Sigma(x) < \Sigma(y)$ implies $x < y$.
Since $\wt{V}(\sigma)$ is closed under the addition of
$\pm(1,\ldots,1)$, it follows that
$\wt{V}(\sigma)$ is isomorphic, as a totally ordered set, 
to a subset of $\Z$ which is unbounded both from below and from above.
This shows that $\wt{V}(\sigma)$ is isomorphic to $\Z$.
\end{proof}

Let $x \in V(\sigma)$ and
let us choose a lift $\wt{x} \in \wt{V}(\sigma)$ of $x$.
Lemma \ref{lem:totally} implies that
there exists a maximum element $\wt{x}'$ of $\wt{V}(\sigma)$
satisfying $\wt{x}' < \wt{x}$.
We set $e(\sigma,x) = \wt{x} - \wt{x}' \in \Z^d$.
Then $e(\sigma,x)$ does not depend on the choice of the lift $\wt{x}$
and we have $e(\sigma,x) > 0$ and 
$\sum_{x \in V(\sigma)} e(\sigma,x) = (1,\ldots,1)$.

\subsection{ }
Now let us assume that $\sigma \in A_{d-1}$.
Then the set $\{e(\sigma,x)\ |\ x \in V(\sigma)\}$ is equal to
the set $\{\be_1,\ldots,\be_d\}$ where
$\be_1, \ldots, \be_d$ denotes the standard $\Z$-basis of $\Z^{\oplus d}$.
This implies that, for any $i \in \{1,\ldots,d\}$,
there exists a unique element $x_i \in V(\sigma)$
satisfying $e(\sigma,x_i)=\be_i$. The map
$\wt{[\sigma]} :\{1,\ldots,d\} \to V(\sigma)$ that sends
$i$ to $x_i$ is bijective. Hence the map $\wt{[\sigma]}$
defines an element of $T(\sigma)$ 
(see Section~\ref{sec:orientation} for the definition).
We denote by $[\sigma]$ the class of $\wt{[\sigma]}$
in $O(\sigma)$.
We let $\wh{\beta}$ denote the element
$\wh{\beta} = (\beta_{\nu})_{\nu \in A_{d-1}'}$
in $\prod_{\nu \in A_{d-1}'} \Z$ where $\beta_\nu =1$ if
$\nu = [\sigma]$ for some $\sigma \in A_{d-1}$ and
$\beta_\nu=0$ otherwise.
%
%
We denote by $\beta$ the class of $\wh{\beta}$
in $(\prod_{\nu \in A_{d-1}'} \Z)_\pmone$.

\subsection{ }
Recall that we defined in Section~\ref{sec:def homology} a chain complex 
which computes the Borel-Moore homology of $A_\bullet$.

\begin{prop}
The element $\beta \in (\prod_{\nu \in A_{d-1}'} \Z)_\pmone$ 
is a $(d-1)$-cycle in the chain complex which computes
the Borel-Moore homology of $A_\bullet$.
\end{prop}
\begin{proof}
The assertion is clear for $d=1$ since
the $(d-2)$-nd component of the complex is zero.
Suppose that $d \ge 2$.
Let $\tau$ be an element in $A_{d-2}$.

Since $\sum_{x \in V(\tau)} e(\tau,x) = (1,\ldots,1)$,
there exists a unique vertex $y \in V(\tau)$ such that
$e(\tau,x)$ belongs to $\{\be_1,\ldots,\be_d\}$ for $x \neq y$
and $e(\tau,y)$ is equal to the sum of two distinct element
of $\{\be_1,\ldots,\be_d\}$.
Let us write $e(\tau,y) = \be_{j} + \be_{j'}$.
Let us choose a lift $\wt{y} \in \wt{V}(\tau)$ of $y$
and set $\wt{y}' = \wt{y} - e(\tau,y)$.
Then $\wt{y}' \in \wt{V}(\tau)$ and 
there are exactly two elements
in $\Z^{\oplus d}$ which is larger than $\wt{y}'$
and which is smaller than $\wt{y}$,
namely $\wt{y}-\be_j$ and $\wt{y}-\be_{j'}$.
We set $y_1 = \pi(\wt{y}-\be_j)$ and
$y_2 = \pi(\wt{y}-\be_{j'})$.
For $i=1,2$, let $\sigma_i \in A_{d-1}$ denote
the unique element 
satisfying $V(\sigma_i) = V(\tau) \cup \{y_i\}$.
It is easily checked that 
the set of the elements in $A_{d-1}$
which has $\tau$ as a face is equal to
$\{\sigma_1,\sigma_2\}$. 
Let $\iota: V(\sigma_1) \cong V(\sigma_2)$
denote the bijection such that
$\iota(x)=x$ for any $x \in V(\tau)$ and
$\iota(y_1)=y_2$.
Then the composite $\iota \circ \wt{[\sigma_1]}$
is equal to the composite $\wt{[\sigma_2]} \circ (jj')$,
where $(jj')$ denotes the transposition of $j$ and $j'$.
Since the signature of $(jj')$ is equal to $-1$,
it follows that the component 
in $(\prod_{\nu \in O(\tau)} \Z)_\pmone$ of the
image of $\beta$ under the boundary map 
$(\prod_{\nu \in A_{d-1}'} \Z)_\pmone 
\to (\prod_{\nu' \in A_{d-2}'} \Z)_\pmone$
is equal
to zero. This proves the claim.
\end{proof}

\begin{definition}
We refer to the class in $H^{\BM}_{d-1}(A_\bullet, \Z)$ 
defined by the $(d-1)$-cycle $\beta$ as 
the fundamental class of the apartment.
\end{definition}

\chapter{Arithmetic subgroups and modular symbols}
Let $d \ge 1$
and $F$ be a global field of positive characteristic.
We give the definition of 
our main object of study, an arithmetic 
subgroup $\Gamma \subset \GL_d(F)$.

Fixing a place $\infty$ of $F$ 
and denoting by $F_\infty$
 the completion at $\infty$
of $F$, arithmetic subgroups 
act on the Bruhat-Tits building
$\cBT_\bullet$
of $\PGL_d(F_\infty)$.

We are interested 
in the $(d-1)$-st Borel-Moore homology of the quotient
$\Gamma \backslash \cBT_\bullet$.
This is an analogue of $\Gamma\backslash \cH$ where 
$\cH$ is the upper half space, $\Gamma$ 
is a congruence subgroup of 
$\SL_d(\Z)$.   The first cohomology
group is known to be related to 
the space of elliptic modular forms.

In Section~\ref{sec:def modular symbol}, 
we define modular symbols in the Borel-Moore
homology group.   
Recall that a building is a union of subsimplicial complexes 
called apartements.  
They are indexed by bases of 
$F_\infty^{\oplus d}$, 
but we restrict to those coming 
from bases of $F^{\oplus d}$
(we say $F$-basis).    
A modular symbol is the 
image of the fundamental class of an apartment
associated with an $F$-basis.

We state our main result in Section~\ref{sec:state main theorem}.
It computes a bound on the index
of the subgroup generated
by modular symbols in 
the Borel-Moore homology of the quotient.

Recall that one property of an arithmetic subgroup $\Gamma$ 
of $\GL_d(\Q)$ is that 
there exists a subgroup of finite index $\Gamma' \subset 
\Gamma$ such that $\Gamma$ is torsion free.
In positive characteristic
of characteristic $p$, 
this is not expected to hold true.
Instead, there always exist a $p'$-torsion subgroup (i.e.,
the torsion subgroup is a $p$-group) of finite index.
We compute a bound on the index in Section~\ref{sec:p'torsion}.

\section{Arithmetic subgroups}
\label{sec:4.1.1}
We define arithmetic subgroups.
The main examples of arithmetic subgroups are
the congruence subgroups.
We verify the basic properties (see (1)-(5) below)
in Section~\ref{sec:def Gamma}.   
We then see 
that an arithmetic subgroup acts on the Bruhat-Tits
building with finite stabilizer groups
and 
the map from an apartment to a quotient is a 
locally finite map (so the pushforward map
of Borel-Moore homology groups is defined).

\subsection{}
\label{sec:def arithmetic}
Let us give the setup.   (This is common when considering 
Drinfeld modules.)
We let $F$ denote a global field of positive characteristic $p>0$.
Let $C$ be a proper smooth curve over a finite field
whose function field is $F$.
Let $\infty$ be a place of $F$ and let $K=F_\infty$ denote
the local field at $\infty$.
We let $A=H^0(C \setminus \{\infty\}, \cO_C)$.
Here we identified a closed point of $C$ and
a place of $F$.
We write $\wh{A}=\varprojlim_I A/I$, where
the limit is taken over the nonzero ideals of $A$.
We let  $\A^\infty=\wh{A} \otimes_A F$ denote the ring of 
finite adeles.

\begin{definition}
A subgroup $\Gamma \subset \GL_d(K)$ 
is called arithmetic subgroup
if there exists a compact open 
subgroup $\bK^\infty \subset \GL_d(\A^\infty)$
such that  $\Gamma=\GL_d(F) \cap \bK^\infty \subset \GL_d(K)$.
\end{definition}

\subsection{Congruence subgroups}
\label{sec:congruence}
Let $I \subset A$ be a nonzero ideal.   Set 
$\Gamma_I= 
\Ker
[
\GL_d(A) \to \GL_d(A/I)
]
$
where the map is the canonical map induced by
the projection $A \to A/I$.
Then $\Gamma_I$ is an arithmetic subgroup
because we can take 
$\bK^\infty$ to be
$\Ker
[\GL_d(\wh{A}) 
\to 
\GL_d(\wh{A}/I \wh{A})
]$.
These are also known as congruence subgroups.

\subsection{}
Let $\Gamma$ be an arithmetic subgroup.
Then $\Gamma \cap \mathrm{SL}_d(F)=\Gamma \cap \mathrm{SL}_d(K)$
is a subgroup of $\Gamma$ of finite index, and is an $S$-arithmetic subgroup
of $\mathrm{SL}_d$ over $F$ for $S=\{\infty\}$ 
in the paper of Harder \cite{Harder}.
\subsection{}
\label{sec:def Gamma}
Let $\Gamma \subset \GL_d(K)$ be a subgroup.
We consider the following Conditions (1) to (5) 
on $\Gamma$.
\begin{enumerate}
\item $\Gamma \subset \GL_d(K)$ is a discrete subgroup,
\item $\{\det(\gamma) \,|\, \gamma \in \Gamma \} \subset O_\infty^\times$
where $O_\infty$ is the ring of integers of $K$,
\item $\Gamma \cap Z(\GL_d(K))$ is finite.
\end{enumerate}
Let $A_\bullet=A_{v_1,\dots,v_d,\bullet}$ 
denote the apartment corresponding to a basis 
$v_1,\dots, v_d\in K^{\oplus d}$
(defined in Section~\ref{sec:521}).
\begin{enumerate}
\setcounter{enumi}{3} 
\item For any apartment 
$A_\bullet=A_{v_1,\dots, v_d,\bullet}$
with $v_1, \dots, v_d \in F^{\oplus d}$, 
the composition 
$A_\bullet \hookrightarrow \cB\cT_\bullet \to \Gamma \bsl \cB\cT_\bullet$ 
is quasi-finite,
that is, the inverse image of any simplex by this map is a finite set.
\item (Harder) The cohomology group $H^{d-1}(\Gamma, \Q)$
is a finite dimensional $\Q$-vector space.
\end{enumerate}
We verify below that the conditions above are 
satisfied for any arithmetic subgroup~$\Gamma$.

\subsection{}
\label{sec:finite stab}
These properties are used in the following way.
The condition (2) implies that each element in the stabilizer group of a simplex fixes the vertices of the simplex.
Under the condition (1), 
the condition (3) implies that the stabilizer of a simplex is finite.
This implies that the $\Q$-coefficient group homology of $\Gamma$ 
and the homology of $\Gamma \bsl |\cB\cT_\bullet|$ are isomorphic.
The condition (4) will be used to define a class in Borel-Moore homology 
of $\Gamma\bsl\cB\cT_\bullet$ starting from an apartment 
(Section \ref{sec:def modular symbol}).  
We do not use Condition (5) in this form but we record it here
because it is related to the finite dimensionality
of the space of cusp forms of fixed level.

\subsection{}
For the rest of this subsection, we give a proof
that these Conditions hold true.
\begin{prop}
For an arithmetic subgroup $\Gamma$,
Conditions (1)-(5) 
of Section~\ref{sec:def Gamma} hold.
\end{prop}
\begin{proof}
This follows from Lemmas below.
\end{proof}

\begin{lem}
For an arithmetic subgroup $\Gamma$,
Conditions (1), (2), and (3) hold.
\end{lem}
\begin{proof}
Condition (1) holds trivially.
We note that there exists an element $g \in \GL_d(\A^\infty)$
such that $g \bK^\infty g^{-1} \subset \GL_d(\wh{A})$.
Since $\det(\gamma) \in F^\times \cap \wh{A}^\times \subset O_\infty^\times$ for $\gamma \in \Gamma$, (2) holds.
Because $F^\times \cap \GL_d(\wh{A})$ is finite, (3) holds.
\end{proof}

\begin{lem}
\label{lem:arithmetic}
Let $\Gamma$ be an arithmetic subgroup. 
Then (4) holds.
\end{lem}

\begin{proof}
We show that the inverse image of each simplex 
of $\Gamma\bsl \cB\cT_\bullet$ 
under the map in (4) is finite.
For $0 \le i \le d-1$, 
the set of $i$-dimensional simplices $\cB\cT_i$ is 
is identified 
(see Section~\ref{sec:6.1.2} for the 
identification) with 
the coset $\GL_d(K)/\wt{\bK}_\infty$
for an open subgroup
$\wt{\bK}_\infty \subset
\GL_d(K)$
which contains
$K^\times \bK_\infty$
as a subgroup of finite index
for some compact open subgroup
$\bK_\infty \subset \GL_d(K)$.

Let $T \subset \GL_d$ denote the diagonal maximal torus.
The set of simplices of 
$A_\bullet$ of fixed dimension is identified
with the image of the map
\[
\coprod_{w \in S_d}
gwT(K)
\to
\GL_d(K)/\wt{\bK}_\infty
\]
for some $g \in \GL_d(F)$.

Since $S_d$ is a finite group, it then suffices to show that for any
$w \in S_d$, the map
$$
\mathrm{Image}[gwT(K) \to \GL_d(K)/K^\times \bK_\infty ]
\to \Gamma \bsl \GL_d(K) /K^\times \bK_\infty
$$
is quasi-finite.
The inverse image under the last map of the image of 
$gwt \in gwT(K)$ 
is isomorphic to
the set
\[
\begin{array}{ll}
\{\gamma \in\Gamma\,|\, 
\gamma gwt \in gwT(K) 
K^\times \bK_\infty\}
&=\Gamma \cap 
gwT(K) K^\times
\bK_\infty (gwt)^{-1}\\
&=\Gamma \cap 
(gw)T(K) t \bK_\infty
t^{-1}
(gw)^{-1}.
\end{array}
\]
Hence, if we let $g'=gw$
and $\bK'_\infty
=t\bK_\infty t^{-1}$,
this set equals
\[
\begin{array}{ll}
\Gamma\cap 
g'T(K) \bK'_\infty g'^{-1}
&
=\GL_d(F) \cap 
(\bK^\infty \times 
g'T(K) \bK'_\infty g'^{-1})
\\
&=
g'(\GL_d(F) \cap (g'^{-1}
\bK^\infty
g' \cap T(K)
\bK'_\infty)
g'^{-1}.
\end{array}
\]
The finiteness of this set is proved in the following lemma.
\end{proof}

\begin{lem}
For any compact open subgroup $\bK \subset \GL_d(\A)$, 
the set $\GL_d(F) \cap T(K) \bK$ is finite.
\end{lem}

\begin{proof}
Let $U=T(K) \cap \bK$.
Then $T(O_\infty) \supset U$ and is of finite index.
Note that there exist a non-zero ideal 
$I\subset A$ and an integer $N$ such that 
$\bK \subset I^{-1}\varpi_\infty^{-N} 
\mathrm{Mat}_d(\wh{A})\times \mathrm{Mat}_d(O_\infty)$
where $\varpi_\infty$ is a uniformizer in $O_\infty$.

Let $\alpha:T(K)/U \to T(K)/T(O_\infty) \cong \Z^{\oplus d}$ 
be the (quasi-finite) map induced by the 
inclusion $U \subset T(O_\infty)$.
For $h \in T(K)$, we write
$(h_1,\dots, h_d)=\alpha(h)$.
Then for $i=1,\ldots,d$, the 
$i$-th row of $h \bK$ is contained in
$(I^{-1}\wh{A} \times \varpi_\infty^{-N}\varpi_\infty^{h_i}O_\infty)^{\oplus d}$.
Hence, for sufficiently large $h_i$, 
the intersection  $h \bK \cap \GL_d(F)$ is empty.  
We then have, for sufficiently large $N'$,

\begin{equation}\label{lem q-finite}
\GL_d(F)\cap T(K) \bK 
=\coprod_{h\in T(K)/U, \atop 
h_1, \ldots, h_d \le N'} 
\GL_d(F) \cap h\bK.
\end{equation}

The adelic norm of the determinant of an element in $\GL_d(F)$ is 1, 
while that of an element in $h \bK$ is 
$|\det h\,|_\infty = \sum_{i=1}^d h_i$.
So (\ref{lem q-finite}) equals
\[\displaystyle 
\coprod_{h\in T(K)/U, h_i \le N', \sum h_i=0} 
\GL_d(F) \cap h \bK.
\]
The index set of the disjoint union above is finite 
since $\alpha$ is quasi-finite,
and $\GL_d(F) \cap h \bK$ is finite
since $\GL_d(F)$ is discrete and 
$h \bK$ is compact.
The claim follows.
 \end{proof}

\begin{lem}
\label{lem:Harder}
Let $\Gamma$ be an arithmetic subgroup.  
Then (5) holds.
\end{lem}
\begin{proof}
This follows from \cite[p.136, Satz 2]{Harder}.
 \end{proof}

\section{$p'$-torsion free subgroups}
\label{sec:p'torsion}
For an arithmetic subgroup (or congruence subgroup) 
$\Gamma$
of $\SL_d(\Z)$, there always exists
a torsion free subgroup $\Gamma' \subset \Gamma$
of finite index.
In our positive characteristic (of characteristic $p$)
setup, an arithmetic subgroup always contains
some $p$-torsion.   We therefore define 
an arithmetic subgroup to be $p'$-torsion
free if any torsion element has a $p$-power order.
A similar fact that any arithmetic subgroup
contains a $p'$-torsion free subgroup of 
finite index holds true.

\subsection{}
Let $p>0$ denote the characteristic of $F$.
\begin{definition}
We say that a subgroup $\Gamma \subset \GL_d(K)$ is $p'$-torsion free
if any element of $\Gamma$ of finite order is of order a power of $p$.
\end{definition}
\begin{lem} \label{lem:quasi-neat criterion1}
Let $\Gamma \subset \GL_d(F)$ be an arithmetic subgroup and
let $v \neq \infty$ be a place of $F$.
Let $F_v$ denote the completion of $F$ at $v$, and 
$\cO_v \subset F_v$ its ring of integers.
Suppose that there exists an element $g_v \in \GL_d(F_v)$ 
such that the image of $\Gamma$ in $\GL_d(F_v)$ 
is contained in $g_v (I_d + \wp_v M_d(\cO_v)) g_v^{-1}$,
where $I_d$ is the identity $d$-by-$d$ matrix,
$\wp_v$ is the maximal ideal of $\cO_v$, and
$M_d(\cO_v)$ denotes the ring of $d$-by-$d$ matrices
with coefficients in $\cO_v$.
Then $\Gamma$ is $p'$-torsion free.
\end{lem}

\begin{proof}
It suffices to show that any matrix $h \in I_d + \wp_v M_d(\cO_v)$
of finite order is of order a power of $p$.
Let us fix an algebraic closure $\overline{F}_v$ of $F_v$ and
let $\overline{\cO}_v$ denote its ring of integers.
Then $\overline{\cO}_v$ is a valuation ring.
Let $h \in I_d + \wp_v M_d(\cO_v)$ be an element of finite order. 
Then any eigenvalue $\alpha$ of $h$
in $\overline{F}_v$ belongs to $\overline{\cO}_v$ and is congruent to $1$
modulo the maximal ideal of $\overline{\cO}_v$.
Since $h$ is of finite order, $\alpha$ is a root of unity.
This implies $\alpha=1$. Hence $(h -I_d )^N = 0$ for sufficiently
large $N$. Choose a power $q$ of $p$ satisfying $q \ge N$.
Then $h^q - I_d = (h-I_d)^q =0$. This shows that the order of $h$
is a power of $p$.
\end{proof}

\begin{cor}\label{cor:quasi-neat criterion2}
Let $\Gamma \subset \GL_d(F)$ be an arithmetic subgroup and
let $v \neq \infty$ be a place of $F$.
Suppose that, as a subgroup of $\GL_d(F_v)$, the group
$\Gamma$ is contained in a pro-$p$ open compact subgroup 
$\bK_v$ of $\GL_d(F_v)$.
Then $\Gamma$ is $p'$-torsion free.
\end{cor}

\begin{proof}
Let us choose an open subgroup $\bK'_v$ of $\bK_v \cap 
(I_d + \wp_v M_d(\cO_v))$
such that $\bK'_v$ is a normal subgroup of $\bK_v$.
Then by Lemma \ref{lem:quasi-neat criterion1},
$\Gamma' = \Gamma \cap \bK'_v$ is $p'$-torsion free.
Since $\Gamma'$ is equal to the kernel of the composite
$\Gamma \inj \bK_v \surj \bK_v/\bK'_v$ and
$\bK_v/\bK'_v$ is a finite $p$-group,
$\Gamma'$ a normal subgroup of $\Gamma$ and the quotient 
$\Gamma/\Gamma'$ is a finite $p$-group.
It follows that $\Gamma$ is $p'$-torsion free.
\end{proof}

\begin{cor}
Let $I \subset A$ be a nonzero ideal.   Then $\Gamma_I$
(see Section~\ref{sec:congruence}) is $p'$-torsion free.     
\end{cor}
\begin{proof}
Let $v$ be a prime which divides $I$.   Then $\Gamma$ is contained 
in the compact open subgroup $I_d + \wp_v M_d(\cO_v)$ where 
$I_d$ is the identity matrix.   Hence 
the claim follows from Corollary~\ref{cor:quasi-neat criterion2}.
\end{proof}

\subsection{}As the compact open subgroups 
\[
\Ker[\GL_d(\widehat{A}) 
\to \GL_d(\widehat{A}/J\widehat{A})]
\]
as $J$ runs over the ideals form a 
fundamental system of neihborhoods of 
the identity matrix
$I_d$, an arithmetic subgroup $\Gamma$
 contains some congruence subgroup
$\Gamma_I$ and it is of finite index in 
$\Gamma$.

\section{Arithmetic quotients of the Bruhat-Tits building}
\label{sec:explicit beta2}
Let us define simplicial complex $\Gamma \bsl \cBT_\bullet$
for an arithmetic subgroup $\Gamma$
and check that the canonical quotient map 
is well defined.   

\subsection{}
We need a lemma.
\begin{lem} \label{lem:stabilizer}
Let $i \ge 0$ be an integer, let
$\sigma \in \cBT_i$ and let
$v,v' \in V(\sigma)$ be two vertices with
$v \neq v'$. Suppose that
an element $g \in \GL_d(K)$ satisfies
$|\det\, g|_\infty =1$. Then we have
$gv \neq v'$.
\end{lem}
\begin{proof}
Let $\wt{\sigma}$ be an element
$(L_j)_{j \in \Z}$ in $\wt{\cBT}_i$ 
such that the class of $\wt{\sigma}$ in $\cBT_i$ is equal
to $\sigma$. 
There exist two integers $j,j' \in \Z$
such that $v$, $v'$ is the class of $L_j$, $L_{j'}$,
respectively.
Assume that $g v =v'$. Then there exists an integer
$k \in \Z$ such that $L_j g^{-1}
= \varpi_\infty^{k} L_{j'} = L_{j' + (i+1) k}$.
Let us fix a Haar measure $d\mu$ of the $K$-vector space
$V_\infty=K^{\oplus d}$. 
As is well-known, the push-forward of $d\mu$
with respect to the automorphism $V_\infty \to V_\infty$
given by the right multiplication by $\gamma$ is equal
to $|\det\, \gamma|_\infty^{-1} d\mu$ 
for every $\gamma \in \GL_d(K)$.
Since $|\det\, g|_\infty =1$, it follows from the
equality $L_j g^{-1} = L_{j'+(i+1)k}$ that
the two $\cO_\infty$-lattices $L_j$ and $L_{j'+(i+1)k}$ have
a same volume with respect to $d\mu$.
Hence we have $j=j'+(i+1)k$, which implies 
$L_j = \varpi_\infty^k L_{j'}$. It follows that the class of
$L_j$ in $\cBT_0$ is equal to the class of $L_{j'}$,
which contradicts the assumption $v \neq v'$.
 \end{proof}

\subsection{quotients and the canonical maps}
Let 
$\Gamma \subset \GL_d(K)$
be an arithmetic subgroup.
It follows from Lemma~\ref{lem:stabilizer}
(using Condition (2) of Section~\ref{sec:def Gamma})
that for each $i \ge 0$
and for each $\sigma \in \cBT_i$,
the image of $V(\sigma)$ under the 
surjection $\cBT_0 \surj \Gamma \bsl \cBT_0$
is a subset of $\Gamma \bsl \cBT_0$
with cardinality $i+1$.
We denote this subset by $V(\cl(\sigma))$, since
it is easily checked that it depends only on the class 
$\cl(\sigma)$ of $\sigma$ in $\Gamma \bsl \cBT_i$.
Thus the collection
$\Gamma \bsl \cBT_\bullet =(\Gamma \bsl \cBT_i)_{i \ge 0}$ 
has a canonical structure of a simplicial complex such that
the collection of the canonical surjection
$\cBT_i \surj \Gamma \bsl \cBT_i$ 
is a map of simplicial complexes 
$\cBT_\bullet \surj \Gamma \bsl \cBT_\bullet$.

\section{Definition of modular symbols} 
\label{sec:def modular symbol}
We define modular symbols here.

\subsection{}
Notice that apartments are defined for any basis of 
$F_\infty^{\oplus d}$.   However, the ones of interest in 
number theory are those associated with the $F$-bases
(i.e., a basis of $F^{\oplus d}$ regarded as a basis 
of $F_\infty^{\oplus d}$).

\subsection{}
Let $v_1,\dots, v_d$ be an $F$-basis 
(that is, a basis of $F^{\oplus d}$
regarded as a basis of $F_\infty^{\oplus d}$).
We consider the composite
\begin{equation} \label{quasifinite}
A_\bullet \xto{\iota_{v_1,\ldots,v_d}}
\cBT_\bullet \to \Gamma \bsl \cBT_\bullet.
\end{equation}
Condition (4) implies that 
the map (\ref{quasifinite})  
is a finite map of simplicial complexes
in the sense of Section~\ref{sec:quasifinite}.
It follows 
that the map (\ref{quasifinite}) induces a homomorphism
$$
H^\BM_{d-1}(A_\bullet, \Z) \to
H^\BM_{d-1}(\Gamma \bsl \cBT_\bullet, \Z).
$$
We let $\beta_{v_1,\ldots,v_d}
\in H^\BM_{d-1}(\Gamma \bsl \cBT_\bullet, \Z)$ denote the
image under this homomorphism 
of the element $\beta \in H^\BM_{d-1}(A_\bullet, \Z)$ introduced
in Section~\ref{sec:fundamental class}.
We call this the class of the apartment
$A_{v_1,\dots, v_d,\bullet}$.

\begin{definition}
We let $\MS(\Gamma)_\Z \subset H^\BM_{d-1}(\Gamma \bsl \cBT_\bullet, \Z)$
denote the submodule generated by the classes $\beta_{v_1,\ldots,v_d}$
as $v_1, \dots, v_d$ runs over the set of ordered $F$-bases.
\end{definition}

\section{Statement of Main Theorem}
\label{sec:state main theorem}
We are ready to state our theorem.
The proof begins in Chapter~\ref{ch:pf for ums}
and ends in Chapter~\ref{ch:compare ms}.
\begin{thm}\label{lem:apartment}
Let $\Gamma \subset \GL_d(K)$
be an arithmetic subgroup.
We write $\MS(\Gamma)_\Z \subset 
H^\BM_{d-1}(\Gamma \bsl \cBT_\bullet, \Z)
$
for the submodule generated by the classes of apartments
associated to $F$-bases.
\begin{enumerate}
\item 
We have 
\[
H^\BM_{d-1}(\Gamma \bsl \cBT_\bullet, \Q)=MS(\Gamma)_\Z \otimes \Q
\]
\item Suppose that $\Gamma$ is $p'$-torsion free.   Set
$$
e(d) = 
(d-2)\left(1+ \frac{(d-1)(d-2)}{2}\right)
$$
Then 
\[
p^{e(d)}H^\BM_{d-1}(\Gamma \bsl \cBT_\bullet, \Z) \subset \MS(\Gamma)_\Z.
\]
\item Let $v \neq \infty$ be a prime of $F$, and
let $F_v$ denote the completion of $F$ at $v$.
Let $\bK_v$ be a pro-$p$ open compact subgroup of $\GL_d(F_v)$.
Let us consider the intersection $\Gamma' = \Gamma \cap \bK_v$
in $\GL_d(F_v)$.
Then 
\[
p^{e(d)} [\Gamma:\Gamma'] H^\BM_{d-1}(\Gamma \bsl \cBT_\bullet, \Z)
\subset \MS(\Gamma)_\Z.
\]
\item Let $v_0 \neq \infty$ be a prime of $F$ such that the cardinality 
$q_0$ of the residue field $\kappa(v_0)$ at $v_0$ is smallest 
among those at the primes $v \neq \infty$.
Set $N(d) = \prod_{i=1}^{d} (q_0^i-1)$.
Then 
\[
p^{e(d)} N(d) 
H^\BM_{d-1}(\Gamma \bsl \cBT_\bullet, \Z)
\subset 
\MS(\Gamma)_\Z
\].
\item Suppose that $d=2$. Then
\[
H^\BM_{1}(\Gamma \bsl \cBT_\bullet, \Z)
=\MS(\Gamma)_\Z.
\]
\end{enumerate}
\end{thm}

\chapter{Automorphic Forms with Steinberg at infinity}
Let $d \le 1$ and $F$ be a global field.
An automoprhic form for $\GL_d$ over $F$
is a $\C$-valued function on 
$\GL_d(F) \backslash \GL_d(\A_F)$,
where $\A_F$ is the ring of adeles,
satisfying some more conditions.
Instead of studying all automorphic forms, 
we study certain subset
consisting of those
satisfying a certain condition at a 
fixed place $\infty$ of $F$. 
We call them automorphic forms 
with Steinberg at infinity

In Section~\ref{sec:Laumon review},
we recall the basics and some 
results on automorphic forms from
Laumon's book \cite{Laumon2}.   The reader is 
also referred to \cite{BJ} and Cogdell's lectures 
\cite{CKM}.

In Section~\ref{sec:Steinberg autom},
we give the definition of the automorphic forms of 
interest to us.   These are the automorphic forms 
that appear when studying the cohomology
of Drinfeld modular varieties.

We can apply known results to
determine (Proposition~\ref{prop:66_3})
the direct sum decomposition,
as automorphic representation,
of the space of 
cusp forms that are Steinberg at infinity.
For the space of all automorphic forms
that are Steinberg at infinity,
we have Theorem~\ref{7_PROP1}.
This does not seem to follow from previously 
known results.   We remark that this space is 
not contained in the space of square integrable 
automorphic forms.   The proof of this theorem
is given in Chapter~\ref{ch:pf of Thm17}.

\section{Automorphic forms with values in $\C$}
\label{sec:Laumon review}
Let $F$ be the global field in positive characteristic.  
We fix a place $\infty$ of $F$.
We write $\A$ for the ring of adeles of $F$.

\subsection{}
The contents of this section is summarized in
the following two diagrams.  We define 
each object and explain the inclusions:
\\
\begin{tikzcd}
C_{\C,\chi}^\infty   \arrow[r, phantom, "\subset"]
&
C_{\C}^\infty   \arrow[r, phantom, "\subset"]
&
C_{\C},   
\\
C_{\C, \mathrm{cusp}}^\infty   \arrow[r, phantom, "\subset"]
&
C_{\C}^\infty   \arrow[u, ,sloped, phantom, "\subset"]
\end{tikzcd}
\begin{tikzcd}
\cA_{\C,\chi}^\infty   \arrow[r, phantom, "\subset"]
&
\cA_{\C}^\infty   
\\
\cA_{\C, \mathrm{cusp}}^\infty   \arrow[r, phantom, "\subset"]
&
\cA_{\C}^\infty   \arrow[u, sloped, phantom, "\subset"]
\end{tikzcd}

\subsection{}
We set 
\[
C_\C=\Hom(\GL_d(F)\backslash\GL_d(\A), \C).
\]
This is a $\GL_d(\A)$-module, where an element $g$ acts as 
$(gf)(x)=f(xg)$.
We set 
\[
C_\C^\infty=\displaystyle\bigcup_\bK C_\C^\bK
\]
where $C_\C^\bK$ denotes the invariants,
and $\bK \subset \GL_d(\A)$ runs
over compact open subgroups.   
This is the space of smooth vectors.
This is a $\GL_d(\A)$-submodule of $C_\C$.

Let $\chi_{\C, \infty}$ 
denote the abelian group of smooth complex characters
of $Z(F)\backslash Z(\A)$ \cite[p.4]{Laumon2}.
For $\chi \in \chi_{\C, \infty}$,
we denote by 
\[
C^\infty_{\C,\chi}=C^\infty_\chi(\GL_d(F)\backslash \GL_d(\A), \C)
\subset \C_\C^\infty
\]
the $\C$-vector subspace of the functions $\varphi \in C^\infty_\C$
such that
\[
\varphi(zm)=\chi(z)\varphi(m)
\]
for $z\in Z(\A), m \in \GL_d(\A)$ (\cite[p.9, 9.1.9]{Laumon2}.

\subsection{}
We let 
\[
\cA_\C \subset C_\C^\infty
\]
denote the space of automorphic forms.
We use the definitions from \cite[p.2, 9.1]{Laumon2}.
Then we set \cite[p.12, 9.1.14]{Laumon2}
\[
A_{\C, \chi}=\cA_\C \cap \C_{\C, \chi}^\infty
\]
Let $C^\infty_{\C, \mathrm{cusp}} \subset C^\infty_\C$ denote 
the $\C$ subvector space of cuspidal functions 
(see \cite[p.15, 9.2.3]{Laumon2} for definition).
We write 
\[
\cA_{\C, \mathrm{cusp}}=\cA_\C \cap C^\infty_{\C, cusp} \subset C^\infty_\C
\]
for the space of cusp forms.

For $\chi \in \chi_\C^\infty$, we 
set 
\[
\cA_{\C, \mathrm{cusp}, \chi}=\cA_{\C, \mathrm{cusp}} \cap \cA_{\C, \chi}
\]

\subsection{}
Let $\chi \in \chi_\C^\infty$ be a unitary character.
We let $\cA^2_{\C, \chi} \subset \cA_{\C, \chi}$
denote the $\C$-vector subspace 
of square-integrable automorphic forms 
(\cite[p.24, 9.3]{Laumon2}).
Lemma 9.3.3 of \cite{Laumon2} says
\[
\cA_{\C, \chi, \mathrm{cusp}} \subset \cA^2_{\C, \chi}.
\]

\section{Automorphic forms with Steinberg at infinity}
\label{sec:Steinberg autom}
\subsection{}
Let $\St_{d,\C}$ denote the Steinberg representation.
It is an admissible representation 
of $\GL_d(F_\infty)$.   
Later we will also use the $\Q$-vector space version $\St_{d,\Q}$.

\subsubsection{}
Let $\cI \subset \GL_d(\cO_\infty)$ denote the Iwahori subgroup.
It is known that the Iwahori fixed part
$\St_{d,\C}^\cI$ of the Steinberg representation
is a one dimensional vector space.   Take a nonzero element $v \in \St_{d,\C}$.
For a $\GL_d(F_\infty)$-module $M$, we set $M_\St$ to be the image by
the evaluation at $v$:
\[
M_\St=
\mathrm{Image}
[\Hom_{\GL_d(F_\infty)}
(\St_{d,\C}, M)
\to M] \subset M
\]
By the one-dimensionality, $M_\St$ does not depend on the 
choice of $v$.

We have the following lemma.
\begin{lem}
We have 
\\
(1) $C_{\C,\St}^\infty=\cA_{\C,\St}$,
\\
(2) $C_{\C,c, \St}^\infty=\cA_{\C, cusp, \St}$.
\end{lem}
\begin{proof}(Proof of (1))
The $\GL_d(F_\infty)$-representation 
generated by a vector of $C^\infty_{\C, \St}$ 
is the Steinberg representation.
Since the Steinberg representation is an admissible representation,
by definition of automorphic form (see, for example, \cite[Prop. 4.5, p.196]{BJ}), 
the vector is an automorphic form.  
\end{proof}

\subsubsection{}
(See \cite[p.35]{Laumon2})
Let $\chi \in \chi^\infty_\C$ 
and let us assume that $\chi$
is trivial on $F_\infty^\times$.
Then $\chi$ is automatically of finite order
(since $F_\infty^\times F^\times\backslash \A^\times$ 
is compact) therefore unitary.
Let $f \in \cA_{\C, \mathrm{cusp},\St}$ and let $\chi$ be 
its central quasi-character.   Then $\chi(F_\infty)=1$
because $F_\infty^\times$ acts trivially for 
the Steinberg representation.    Therefore, as seen above, $\chi$ 
is unitary.
We write $\cA_{\C, \mathrm{cusp}, \St, \chi} \subset \cA_{\C, \mathrm{cusp}, \St}$
for the subspace consisting of those cusp forms whose 
unitary central character is $\chi$.
\subsubsection{}
Let 
\[
\cA_{\chi, \mathrm{disc}}^2 \subset \cA^2_\chi
\]
denote the discrete spectrum (see \cite[9.3.4, p.25]{Laumon2}).
For a quasi-character $\chi \in \chi_\C^\infty$
such that $\chi(F_\infty^\times)=1$
(hence $\chi$ is unitary),
Corollary 9.5.6 of \cite[p.36]{Laumon2} implies
\[
\cA_{\C, \mathrm{cusp}, \chi, \St}
=\cA^2_{\chi, \mathrm{disc}, \St}.
\]

\subsubsection{}
Let 
$C_{\C,c}^\infty \subset C_{\C}^\infty$
denote the $\C$ subvector space of compactly supported functions.
Then Harder's theorem (see \cite[Theorem 9.2.6, p.16]{Laumon2}
implies that $\cA_{\C, \mathrm{cusp}} \subset C_{\C,c}^\infty$.

The elements of 
$C^\infty_{\C,c, \St}$ 
are automorphic and square-integrable
hence $C^\infty_{\C,c, \St} \subset \cA^2_{\mathrm{disc}}$.
Now, for a unitary character $\chi$, 
let us write 
$C^\infty_{\C,c,\chi, \St}=C^\infty_{\C,c,\St} \cap \cA^2_\chi$. 

\begin{proof}(Proof of (2))
From the discussions above, we obtain
\[
\cA_{\C, \mathrm{cusp}, \St, \chi} \subset
C_{\C,c,\chi, \St}^\infty
\subset
\cA^2_{\chi, \mathrm{disc}, \St}
=
\cA_{\C, \mathrm{cusp}, \St, \chi}.
\]
This proves (2).
\end{proof}

\section{Admissiblity}
The admissiblity follows from a result of Harder.
\begin{prop}
\label{prop:admissible}
The vector spaces 
$\cA_{\C, \mathrm{cusp},\St}$ and $\cA_{\C,\St}$ 
are admissible representations of $\GL_d(\A)$.
\end{prop}
\begin{proof}
Because of the condition at $\infty$,
we can use the result of Harder
\cite[p.198, Proposition 5.2]{BJ}
to see that $\cA_{\C, \St}$ is 
admissible.   Then the claim follows
for the subspace $\cA_{\C, \mathrm{cusp}, \St}$.
\end{proof}

\section{Decompositions}
\subsection{}
There are well-known theorems 
on the structure of the space of cusp forms,
which in turn imply the following.
\begin{prop} \label{prop:66_3}
As a representation of $\GL_d(\A^\infty)$,
$$
\cA_{\C, \mathrm{cusp}, \St} \cong
\bigoplus_\pi \pi^\infty,
$$
where $\pi = \pi^\infty \otimes \pi_\infty$ runs over
(the isomorphism classes of) the irreducible cuspidal
automorphic representations 
of $\GL_d(\A)$ such that $\pi_\infty$ is isomorphic to 
the Steinberg representation of $\GL_d(F_\infty)$.
\end{prop}
\begin{proof}
The theorems to use are
\cite[Theorem 9.2.14, p.22]{Laumon2}
due to Gelfand and Piatetski-Shapiro, and 
the multiplicity one theorem 
\cite[Remark 9.2.15, p.22]{Laumon2}
due to Shalika. 
\end{proof}

\subsection{}
We prove the following structure theorem 
for $\cA_{\C,\St}$ in Chapter~\ref{ch:pf of Thm17}.

\begin{thm}\label{7_PROP1}
Let $\pi = \pi^\infty \otimes \pi_\infty$ be an irreducible
smooth representation of $\GL_d(\A)$
such that $\pi^\infty$ appears as a subquotient of
$\cA_{\C,\St}$.
%
Then there exist an integer $r \ge 1$, 
a partition $d=d_1 + \cdots + d_r$ of $d$,
and irreducible cuspidal automorphic representations $\pi_i$
of $\GL_{d_i}(\A)$ for $i=1,\ldots,d$ which satisfy
the following properties:
\begin{itemize}
\item For each $i$ with $0 \le i \le r$,
the component $\pi_{i,\infty}$ at $\infty$ of $\pi_i$ is 
isomorphic to the
Steinberg representation of $\GL_{d_i}(F_\infty)$.
\item Let us write $\pi_i = \pi_i^\infty \otimes \pi_{i,\infty}$.
Let $P \subset \GL_d$ denote the standard parabolic subgroup
corresponding to the partition $d=d_1 + \cdots + d_r$.
Then $\pi^\infty$ is isomorphic to a subquotient of the unnormalized 
parabolic induction $\Ind_{P(\A^\infty)}^{\GL_d(\A^\infty)}
\pi_1^\infty \otimes \cdots \otimes \pi_r^\infty$.
\end{itemize}
Moreover for any subquotient $H$ of 
$\cA_{\C,\St}$
which is of finite length as a representation of $\GL_d(\A^\infty)$,
the multiplicity of $\pi$ in $H$ is at most one.
\end{thm}

\chapter{Double cosets for automorphic forms}
We introduce here the main arithmetic object of study
$X_{\bK,\bullet}$ for a compact open subgroup
$\bK \subset \GL_d(\A^\infty)$.   It is a (generalized)
simplicial complex.   This is an analogue of the 
$\C$-valued points 
of a Shimura variety as a double coset.

In Section~\ref{subsec:X_K}, we define a 
(generalized) simplicial complex in the form of 
double coset.   We show (Proposition~\ref{prop:homology autom})
that (the limit of)
the Borel-Moore homology 
is isomorphic to the space of our automorphic forms
and the homology is isomorphic to the subspace of 
cusp forms.
This gives a precise relation between 
the Borel-Moore homology/homology of the geometry
and the space of automorphic forms.  
In particular, we see how our modular symbols
are related to automorphic forms.
Later, we study this geometry to prove Theorem~\ref{7_PROP1}.

\section{Simplicial complexes for automorphic forms}
\label{subsec:X_K}
We give the definition of some double coset $X_{\bK, \bullet}$
associated with a compact open subgroup $\bK$.

\subsection{}
For an open compact subgroup
$\bK \subset \GL_d(\A^\infty)$,
we let 
$\wt{X}_{\GL_d,\bK,\bullet}$ denote the disjoint union
$\wt{X}_{\GL_d,\bK,\bullet}=(\GL_d(\A^\infty)/\bK) 
\times \cBT_{\bullet}$
of copies of the Bruhat-Tits building $\cBT_{\bullet}$ 
indexed by $\GL_d(\A^\infty)/\bK$.
We often omit the subscript $\GL_d$ on
$\wt{X}_{\GL_d,\bK,\bullet}$ when there is no
fear of confusion. The group $\GL_d(\A)=\GL_d(\A^\infty) 
\times \GL_d(F_\infty)$ acts on
the simplicial complex $\wt{X}_{\bK,\bullet}$ from the left.
$X_{\bK, \bullet}=\GL_d(F) \bsl \wt{X}_{\bK,\bullet}
=\GL_d(F) \bsl (\GL_d(\A^\infty) \times \cBT_{\bullet}/\bK$ of 
$\wt{X}_{\bK,\bullet}$ by the subgroup
$\GL_d(F) \subset \GL_d(\A)$.

For $0\le i \le d-1$, we let $X_{\bK,i}= X_{\GL_d,\bK,i}$ 
denote the quotient $X_{\bK,i} = \GL_d(F) \bsl 
\wt{X}_{\GL_d,\bK,i}.$
We let $X_{\bK,\bullet}=
\GL_d(F) \bsl \wt{X}_{\bK,\bullet}
=\GL_d(F) \bsl (\GL_d(\A^\infty) \times \cBT_{\bullet})/\bK$.

\subsection{}
The non-adelic description is as follows. 
We set $J_\bK =\GL_d(F) \bsl \GL_d(\A^\infty)/\bK$.
For each $j \in J_\bK$, we choose an element 
$g_j \in \GL_d(\A^\infty)$ in the double coset $j$
and set $\Gamma_j = \GL_d(F) \cap g_j \bK g_j^{-1}$.
Then the set $X_{\bK,i}$ is isomorphic to the disjoint union
$\coprod_j \Gamma_j \bsl \cBT_i$. For each $j$,
the group $\Gamma_j \subset \GL_d(F)$ is an arithmetic 
subgroup as defined in Section~\ref{sec:def arithmetic}.  It follows that the tuple
$X_{\bK,\bullet}=(X_{\bK,i})_{0\le i \le d-1}$ forms a simplicial complex 
which is isomorphic to the disjoint union 
$\coprod_{j\in J_\bK} \Gamma_j \bsl \cBT_\bullet$.

\subsection{}
Since  
the simplicial complex $\wt{X}_{\GL_d,\bK,\bullet}$ is locally
finite, it follows that the simplicial complex 
$X_{\bK,\bullet}$ is locally finite.
Hence, as in Section~\ref{sec:def homology},  for an abelian group $M$, we may consider the
cohomology groups with compact support $H_c^{*}(X_{\bK,\bullet},M)$
and the Borel-Moore homology groups $H_*^\mathrm{BM}(X_{\bK,\bullet},M)$. 
of the simplicial complex $X_{\bK,\bullet}$.
%

\subsection{}
Since the simplicial complex $X_{\bK,\bullet}$
has no $i$-simplex for $i \ge d$ as was remarked in
Section~\ref{sec:dimension}, 
it follows that the $(d-1)$-st homology and Borel-Moore homology
are isomorphic to the group of chains.
This implies that the map
$$
H_{d-1}(X_{\bK,\bullet},M) \to 
H_{d-1}^\mathrm{BM}(X_{\bK,\bullet},M)
$$
is injective for any abelian group $M$.
We regard $H_{d-1}(X_{\bK,\bullet},M)$ as a subgroup of
$H_{d-1}^\mathrm{BM}(X_{\bK,\bullet},M)$ via this map.

\section{Pull-back maps for homology groups}
We study some functoriality
for two open compact subgroups; then we are able to define the 
inductive limit of the Borel-Moore homology/homology groups.

\subsection{}
Let $\bK,\bK' \subset
\GL_d(\A^\infty)$ be open compact subgroups 
with $\bK' \subset \bK$. We denote by 
$f_{\bK',\bK}$ the natural projection map
$X_{\bK',i} \to X_{\bK,i}$.
Since $\bK'$ is a subgroup of
$\bK$ of finite index, it follows that 
for any $i$ with $0 \le i \le d-1$ 
and for any $i$-simplex $\sigma \in X_{\bK,i}$, 
the inverse image of $\sigma$ under the
map $f_{\bK',\bK}$ is a finite set.
Let $i$ be an integer with $0 \le i \le d-1$
and let $\sigma' \in X_{\bK',i}$. 
Let $\sigma$ denote the image of $\sigma'$
under the map $f_{\bK',\bK}$.
Let us choose an $i$-simplex $\wt{\sigma}'$ of
$\wt{X}_{\bK',\bullet}$ which is sent to $\sigma'$
under the projection map $\wt{X}_{\bK',\bullet}
\to X_{\bK',\bullet}$.
Let $\wt{\sigma}$ denote the image of $\wt{\sigma}'$
under the map $\wt{X}_{\bK',i} \to \wt{X}_{\bK,i}$.
We let
$$
\Gamma_{\wt{\sigma}'} = \{ \gamma \in \GL_d(F)\ |\ 
\gamma \wt{\sigma}' =\wt{\sigma}' \}
$$
and
$$
\Gamma_{\wt{\sigma}} = \{ \gamma \in \GL_d(F)\ |\ 
\gamma \wt{\sigma} =\wt{\sigma} \}
$$
denote the stabilizer group of $\wt{\sigma}'$ and
$\wt{\sigma}$, respectively.

\begin{lem}
Let the notation be as above.
\begin{enumerate}
\item The group $\Gamma_{\wt{\sigma}}$ is a finite group and
the group $\Gamma_{\wt{\sigma}'}$ is a subgroup of $\Gamma_{\wt{\sigma}}$.
\item The isomorphism class of the group $\Gamma_{\wt{\sigma}'}$
(\resp $\Gamma_{\wt{\sigma}}$) depends only on $\sigma'$
(\resp $\sigma$) and does not depends on the choice of
$\wt{\sigma}'$. 
\end{enumerate}
 
\end{lem}
\begin{proof}
For the case of a 0-dimensonal simplex of (1), 
see \cite[Proof of Theorem 0.8]{Gra}.
In the stabilizer group for a higher dimensional simplex,
the subgroup of elements that fix each vertex is 
of finite index.   This proves (1). 
The claim (2) can be checked easily.
\end{proof}

\begin{definition}
The lemma above shows in particular that the index
$[\Gamma_{\wt{\sigma}} : \Gamma_{\wt{\sigma}'}]$ is finite and
depends only on $\sigma'$ and $f_{\bK',\bK}$. We denote
this index by $e_{\bK',\bK}(\sigma')$ and call it the
ramification index of $f_{\bK',\bK}$ at $\sigma'$.
\end{definition}

\subsection{}
Let $M$ be an abelian group.
Let $i$ be an integer with $0 \le i \le d$.
We set $X'_{\bK,i}= \coprod_{\sigma \in X_{\bK,i}} O(\sigma)$.
The map $f_{\bK',\bK} : X_{\bK',\bullet} \to X_{\bK,\bullet}$
induces a map $X'_{\bK',i} \to X'_{\bK,i}$ which we denote
also by $f_{\bK'.\bK}$.
\subsubsection{}
Let $m = (m_{\nu})_{\nu \in X'_{\bK,i}}$ be an element
of the $\{\pm 1\}$-module $\prod_{\nu \in X'_{\bK,i}} M$.
We define the element $f^*_{\bK',\bK}(m)$
in $\prod_{\nu \in X'_{\bK',i}} M$ to be
$$
f^*_{\bK',\bK}(m) = (m'_{\nu'})_{\nu' \in X'_{\bK',i}}
$$
where for $\nu' \in O(\sigma') \subset X'_{\bK',i}$,
the element $m'_{\nu'} \in M$ is given by
$m'_{\nu'} = e_{\bK',\bK}(\sigma') m_{f_{\bK',\bK}(\nu')}$.
The following lemma can be checked easily.
\begin{lem}
Let the notation be as above.
\begin{enumerate}
\item The map $f^*_{\bK',\bK} : \prod_{\nu \in X'_{\bK,i}} M
\to \prod_{\nu' \in X'_{\bK',i}} M$ is a homomorphism of
$\{\pm 1\}$-modules.
\item The map $f^*_{\bK',\bK} : \prod_{\nu \in X'_{\bK,i}} M
\to \prod_{\nu' \in X'_{\bK',i}} M$ sends an element in 
the subgroup $\bigoplus_{\nu \in X'_{\bK,i}} M 
\subset \prod_{\nu \in X'_{\bK,i}} M$ to an element
in  $\bigoplus_{\nu \in X'_{\bK',i}} M$.
\item For $1 \le i \le d-1$, the diagrams
$$
\begin{CD}
\prod_{\nu \in X'_{\bK,i}} M @>{\wt{\partial}_{i,\prod}}>> 
\prod_{\nu \in X'_{\bK,i-1}} M \\
@V{f_{\bK',\bK}^*}VV @V{f_{\bK',\bK}^*}VV \\
\prod_{\nu' \in X'_{\bK',i}} M @>{\wt{\partial}_{i,\prod}}>> 
\prod_{\nu' \in X'_{\bK',i-1}} M
\end{CD}
$$
and
$$
\begin{CD}
\bigoplus_{\nu \in X'_{\bK,i}} M @>{\wt{\partial}_{i,\oplus}}>> 
\bigoplus_{\nu \in X'_{\bK,i-1}} M \\
@V{f_{\bK',\bK}^*}VV @V{f_{\bK',\bK}^*}VV \\
\bigoplus_{\nu' \in X'_{\bK',i}} M @>{\wt{\partial}_{i,\oplus}}>> 
\bigoplus_{\nu' \in X'_{\bK',i-1}} M
\end{CD}
$$
are commutative.
\end{enumerate}
 
\end{lem}

\section{Limits}
\subsection{}
For an abelian group $M$, we set
\[
H_*(X_{\lim,\bullet},M)=H_*(X_{\GL_d,\lim,\bullet},M)=\varinjlim_{\bK}H_*(X_{\bK,\bullet},M)
\]
and
\[
H_*^\mathrm{BM}(X_{\lim,\bullet},M)=H_*^\mathrm{BM}(X_{\GL_d,\lim,\bullet},M)
=
\varinjlim_{\bK}H_*^\mathrm{BM}(X_{\bK,\bullet},M).
\]
Here the transition maps in the inductive limits are 
given by $f^*_{\bK',\bK}$.


\subsection{}
For $g \in \GL_d(\A^\infty)$, we let 
$\wt{\xi}_g : \wt{X}_{\bK,\bullet}
\xto{\cong} \wt{X}_{g^{-1} \bK g,\bullet}$
denote the isomorphism of simplicial complexes
induced by the isomorphism $\GL_d(\A^\infty)/\bK \xto{\cong}
\GL_d(\A^\infty)/g^{-1}\bK g$ that sends 
a coset $h \bK$ to the coset $hg\cdot g^{-1}
\bK g$ and by the identity on $\cBT_{\bullet}$.
The isomorphism $\wt{\xi}_g$ induces an isomorphism 
$\xi_{g} :X_{\bK,\bullet}
\xto{\cong} X_{g^{-1} \bK g,\bullet}$ of
simplicial complexes. 
%
For two elements
$g, g' \in \GL_d(\A^\infty)$, we have
$\xi_{gg'}=\xi_{g'}\circ \xi_g$.

\subsection{}
The isomorphisms
$\xi_g$ for $g\in \GL_d(\A^\infty)$ give
rise to a smooth action 
(i.e., the stabilizer of each vector is a compact open 
subgroup) of the group $\GL_d(\A^\infty)$
on these inductive limits.
If $M$ is a torsion free abelian group, then
for each compact open subgroup $\bK \subset \GL_d(\A^\infty)$,
the homomorphism
$H_*(X_{\bK,\bullet},M) \to 
H_*(X_{\lim,\bullet},M)$
is injective and its image is equal to the
$\bK$-invariant part $H_*(X_{\lim,\bullet},M)^{\bK}$
of $H_*(X_{\lim,\bullet},M)$.
Similar statement holds for $H_*^\mathrm{BM}$.

\section{the Steinberg representation and harmonic cochains}
In this section, we consider coefficients of $\Q$-vector spaces.
We show how the limit of the Borel-Moore homology/homology 
corresponds to a sub $\GL_d(\A^\infty)$ representation
of the space of automorphic forms.  

\subsection{}
Let $\St_{d, \C}$ denote the Steinberg representation
as defined, for example, in \cite[p.193]{Laumon1}.
It is defined with coefficients in $\C$,
but it can also be defined with coefficients in $\Q$
in a similar manner.  We let $\St_{d, \Q}$ denote the
corresponding representation.

\subsection{}
\begin{definition}
A harmonic cochain with values in a $\Q$-vector space 
$M$ is 
defined as an element of 
$\Hom(H^{d-1}_c(\cBT_\bullet,\Q),M)$.
\end{definition}

\begin{lem}\label{lem:Steinberg}
For a $\Q$-vector space $M$, there is a canonical,
$\GL_d(F_\infty)$-equivariant isomorphism 
between the module of $M$-valued 
harmonic $(d-1)$-cochains and the module
$\Hom_{\Q}(\St_{d,\Q},M)$.
\end{lem}

\begin{proof}
It is shown in~\cite[6.2,6.4]{Borel} that
$\St_{d,\C}$ 
is canonically isomorphic to
$H^{d-1}_c(\cBT_\bullet,\C)$ as a 
representation of $\GL_d(F_\infty)$.
One can check that 
this map is defined over $\Q$. 
This proves the claim.
\end{proof}

\subsection{}
\label{sec:6.1.2}
We let $\cBT_{j,*}$ 
denote the quotient
$\wt{\cBT}_j/F_{\infty}^{\times}$.
This set 
is identified with the set of 
pairs $(\sigma, v)$ with 
$\sigma \in \cBT_j$
and $v \in \cBT_0$ a vertex of $\sigma$, 
which we call 
a pointed $j$-simplex.
Here the element 
$(L_i)_{i\in \Z} \mod K^\times$
of $\wt{\cBT_j}/K^\times$
corresponds to the pair
$((L_i)_{i\in \Z}, L_0)$
via this identification.
 
\subsection{}
We identify the set
$\wt{\cBT}_{0}$ with the coset
$\GL_d(K)/\GL_d(\cO)$
by associating to an element $g \in
\GL_d(K)/\GL_d(\cO)$ the lattice
$\cO_{V}g^{-1}$.
Let 
$\cI=\{(a_{ij})\in \GL_d(\cO) \,|\,
a_{ij}\,\mathrm{mod}\, \varpi =0 \ \text{if}\ i>j\}$
be the Iwahori subgroup.
Similarly, we identify the set $\wt{\cBT}_{d-1}$
with the coset $\GL_d(K)/\cI$
by associating
to an element $g\in \GL_d(K)/\cI$
the chain of lattices $(L_i)_{i \in \Z}$ characterized by
$L_i=\cO_{V}\Pi_ig^{-1}$ for $i=0,\dots,d$.
Here, for $i=0,\dots,d$,
we let $\Pi_i$ denote the diagonal $d\times d$ matrix
$\Pi_i=\diag(\varpi,\ldots,\varpi,1,\ldots,1)$
with $\varpi$ appearing $i$
times and $1$ appearing $d-i$ times.

\subsection{} 
Let $\bK \subset \GL_d(\A)$
be an open compact subgroup.
Let $M$ be a $\Q$-vector space.
Let $\cC^\bK(M)$ denote the 
($\Q$-vector)  space of locally constant 
$M$-valued functions on
$\GL_d(F)\bsl \GL_d(\A)/(\bK \times F_\infty^\times)$.
Let $\cC_c^\bK(M) \subset C^\bK(M)$
denote the subspace of compactly supported functions.
\begin{lem}\label{lem:Steinberg8}
\begin{enumerate}
\item There is a canonical isomorphism
$$
H_{d-1}^\mathrm{BM}(X_{\bK,\bullet},M) 
\cong \Hom_{\GL_d(F_\infty)} (\St_{d,\Q}, \cC^\bK(M)),
$$
where $\cC^\bK(M)$ 
denotes the space of locally constant $M$-valued
functions on $\GL_d(F)\bsl \GL_d(\A)/(\bK\times F_\infty^\times)$.
\item 
Let $v \in \St_{d, \Q}^\cI$ be a non-zero 
Iwahori-spherical vector. 
Then the image of the evaluation map
$$
\begin{array}{l}
\Hom_{\GL_d(F_\infty)}(\St_{d,\Q}, 
\cC^\bK(M)) \\
\to \Map(\GL_d(F)\bsl \GL_d(\A)/(\bK \times 
F_\infty^\times \cI),M)
\end{array}
$$
at $v$ is identified with the image of
the map
$$
\begin{array}{rl}
H_{d-1}^\mathrm{BM}(X_{\bK,\bullet},M) 
 & \to \Map(\GL_d(F)\bsl 
(\GL_d(\A^\infty)/\bK \times \cBT_{d-1,*}),M) \\
& \cong \Map(\GL_d(F)\bsl \GL_d(\A)/(\bK \times 
F_\infty^\times \cI),M).
\end{array}
$$
\end{enumerate}
\end{lem}
\begin{proof}
For a $\C$-vector space $M$,
(1) is proved in \cite[Section 5.2.3]{KY:Zeta elements},
and (2) is \cite[Corollary 5.7]{KY:Zeta elements}.
The proofs and the argument in loc. cit. 
work for a $\Q$-vector space $M$ as well.
\end{proof}

\begin{cor}
\label{cor:Steinberg8}
Under the isomorphism in (1), the subspace
$$
H_{d-1}(X_{\bK,\bullet},M) \subset 
H_{d-1}^\mathrm{BM}(X_{\bK,\bullet},M)
$$
corresponds to the subspace 
$$
\Hom_{\GL_d(F_\infty)}(\St_{d,\Q}, \cC_c^\bK(M))
\subset 
\Hom_{\GL_d(F_\infty)}(\St_{d,\Q}, \cC^\bK(M)).
$$
\end{cor}
\begin{proof}
This follows from Lemma~\ref{lem:Steinberg8} (2)
and the definition of the homology group
$H_{d-1}(X_{\bK, \bullet}, M)$.
\end{proof}

\subsection{}
The proof of the following lemma is straightforward
and is left to the reader.
\begin{lem}
Let the notation be as above.
\begin{enumerate}
\item 
Suppose that $\bK'$ is a normal subgroup of $\bK$.
Then the homomorphism $f^*_{\bK',\bK}$
induces an isomorphism $H^\mathrm{BM}_*(X_{\bK,\bullet},M) 
\cong H^\mathrm{BM}_*(X_{\bK',\bullet},M)^{\bK/\bK'}$
and a similar statement holds for 
$H_*$.
\item
Let $M$ be a $\Q$-vector space. Then the diagrams
$$
\begin{CD}
H^\mathrm{BM}_{d-1} (X_{\bK,\bullet}, M) @>{\cong}>>
\Hom_{\GL_d(F_\infty)}(\St_{d,\Q}, \cC^\bK(M)) \\
@V{f^*_{\bK',\bK}}VV @VVV \\
H^\mathrm{BM}_{d-1} (X_{\bK',\bullet}, M) @>{\cong}>>
\Hom_{\GL_d(F_\infty)}(\St_{d,\Q}, \cC^{\bK'}(M))
\end{CD}
$$
and
$$
\begin{CD}
H_{d-1} (X_{\bK,\bullet}, M) @>{\cong}>>
\Hom_{\GL_d(F_\infty)}(\St_{d,\Q}, \cC_c^\bK(M)) \\
@V{f^*_{\bK',\bK}}VV @VVV \\
H_{d-1} (X_{\bK',\bullet}, M) @>{\cong}>>
\Hom_{\GL_d(F_\infty)}(\St_{d,\Q}, \cC_c^{\bK'}(M))
\end{CD}
$$
are commutative.
Here the horizontal arrows are the isomorphisms
given in Lemma~\ref{lem:Steinberg8} 
and Corollary~\ref{cor:Steinberg8}, 
and the right vertical
arrows are the map induced by the quotient map
$\GL_d(F) \bsl \GL_d(\A) /(\bK' \times F_\infty^\times)
\to \GL_d(F) \bsl \GL_d(\A) /(\bK\times F_\infty^\times)$.
\end{enumerate}
 
\end{lem}

\begin{cor}
\label{prop:homology autom}
We have isomorphisms
\[
H_{d-1}(X_{\lim},\C) \cong \cA_{\C,cusp,\St}
\]
\[
H^\BM_{d-1}(X_{\lim},\C) \cong \cA_{\C,\St}
\]
\end{cor}
\begin{proof}
This follows from the previous lemma 
and the definitions.
\end{proof}

\chapter{Proof of Theorem \ref{7_PROP1}}
\label{ch:pf of Thm17}
We give a proof of Theorem~\ref{7_PROP1}.

In Section \ref{sec:locfree}, we introduce some terminology
for locally free $\cO_C$-modules of rank $d$ 
and then describe the
sets of simplices of $X_{\bK,\bullet}$
in terms of chains of locally free $\cO_C$-modules 
of rank $d$.  
In Section \ref{sec:Xsigma},
we follow Section 4 of \cite{Gra}
and define certain subsimplicial complexes of $X_{\bK, \bullet}$.

\section{Chains of locally free $\cO_C$-modules}
\label{sec:locfree}

\subsection{}
Let $\eta : \Spec F \to C$ denote the
generic point of $C$.
\begin{definition}
For each $g \in \GL_d(\A^\infty)$ and an 
$\cO_\infty$-lattice $L_\infty \subset 
\cO_\infty^{\oplus d}$,
we denote by $\cF[g,L_{\infty}]$ the
$\cO_C$-submodule of $\eta_* F^{\oplus d}$ 
characterized by the following properties:
\begin{itemize}
\item $\cF[g,L_{\infty}]$ is a
locally free $\cO_C$-module of rank $d$.
\item $\Gamma(\Spec A, \cF[g,L_{\infty}])$ 
is equal to the $A$-submodule 
$\wh{A}^{\oplus d} g^{-1} \cap F^{\oplus d}$ of 
$F^{\oplus d} = \Gamma(\Spec A,\eta_* F^{\oplus d})$.
\item Let $\iota_\infty$ denote the morphism $\Spec \cO_\infty \to C$.
Then $\Gamma(\Spec \cO_\infty, \iota_\infty^* \cF[g,L_{\infty}])$
is equal to the $\cO_\infty$-submodule $L_{\infty}$ of
$F_\infty^{\oplus d} = \Gamma(\Spec \cO_\infty, \iota_\infty^* \eta_* F^{\oplus d})$.
\end{itemize}
\end{definition}
\subsection{}
Let $\cF$ be a locally free $\cO_C$-modules
of rank $d$. Let $I \subset A$ be a
non-zero ideal. We regard the $A$-module
$A/I$ as a coherent $\cO_C$-module of
finite length. 

\begin{definition}
A level $I$-structure on
$\cF$ is a surjective homomorphism
$\cF \to (A/I)^{\oplus d}$ of $\cO_C$-modules.
\end{definition}
Let $\bK^\infty_I
\subset \GL_d(\wh{A})$ be 
the kernel of the
homomorphism $\GL_d(\wh{A}) \to \GL_d(\wh{A}/I \wh{A})$. 
The group $\GL_d(A/I) \cong
\GL_d(\wh{A})/\bK^\infty_I$ acts from
the left on the set of level $I$-structures on $\cF$,
via its left action on $(A/I)^{\oplus d}$.
(We regard $(A/I)^{\oplus d}$ as an $A$-module of
row vectors. The left action of $\GL_d(A/I)$ on
$(A/I)^{\oplus d}$ is described as $g\cdot b
= b g^{-1}$ for $g\in \GL_d(A/I)$, $b\in
(A/I)^{\oplus d}$.) 
\begin{definition}
For a subgroup
$\bK \subset \GL_d(\wh{A})$ containing $\bK_I^\infty$,
a level $\bK$-structure on $\cF$ is a
$\bK/\bK^\infty_I$-orbit of level
$I$-structures on $\cF$. 
\end{definition}
For an
open subgroup $\bK \subset \GL_d(\wh{A})$,
the set of level $\bK$-structures on $\cF$
does not depend, up to canonical isomorphisms,
on the choice of an ideal $I$ with
$\bK_I^\infty \subset \bK$.

\subsection{}
Let $\bK \subset \GL_d(\wh{A})$ be an open subgroup.
Let $(g,\sigma)$ be an $i$-simplex of
$\wt{X}_{\bK,\bullet}$.
Take a chain 
\[
\cdots \supsetneqq L_{-1} 
\supsetneqq L_0 \supsetneqq 
L_1 \supsetneqq \cdots
\]
 of $\cO_{\infty}$-lattices
of $F_{\infty}^{\oplus d}$ which represents $\sigma$.
To $(g,\sigma)$ we associate the
chain 
\[
\cdots \supsetneqq \cF[g,L_{-1}] \supsetneqq
\cF[g,L_0] \supsetneqq \cF[g,L_1] \supsetneqq \cdots
\]
of $\cO_C$-submodules of $\eta_* F^{\oplus d}$.
Then the set of $i$-simplices in $\wt{X}_{\bK,\bullet}$
is identified with the set of the equivalences classes
of chains 
\[
\cdots \supsetneqq \cF_{-1} \supsetneqq \cF_0 
\supsetneqq \cF_1 \supsetneqq \cdots
\]
of locally free $\cO_C$-submodules of rank $d$ of
$\eta_* \eta^* \cO_C^{\oplus d}$ with a level 
$\bK$-structure 
such that $\cF_{j-i-1}$ equals the
twist $\cF_{j}(\infty)$ as an $\cO_C$-submodule of
$\eta_* F^{\oplus d}$ with a level $\bK$-structure
for every $j\in \Z$.

Two chains $\cdots \supsetneqq \cF_{-1} \supsetneqq
 \cF_0 \supsetneqq \cF_1 \supsetneqq \cdots$ and 
$\cdots \supsetneqq \cF'_{-1} \supsetneqq \cF'_0 
\supsetneqq \cF'_1 \supsetneqq \cdots$ are equivalent 
if and only if there exists an integer $l$ such that
$\cF_{j} = \cF'_{j+l}$ as an $\cO_C$-submodule of
$\eta_* F^{\oplus d}$ with a level structure
for every $j\in \Z$.

\subsection{}
Let $g\in \GL_d(\A^\infty)$ and let 
$L_\infty$ be an $\cO_\infty$-lattice of $F_\infty^{\oplus d}$.
For $\gamma \in \GL_d(F)$, the two $\cO_C$-submodules
$\cF[g,L_\infty]$ and $\cF[\gamma g,\gamma L_\infty]$ 
are isomorphic
as $\cO_C$-modules. The set of $i$-simplices
in $X_{\bK,\bullet}$
is identified 
with the set of the equivalence classes
of chains
$\cdots \inj \cF_{1} \inj  \cF_0 \inj
\cF_{-1} \inj \cdots$ of injective non-isomorphisms
of locally free $\cO_C$-modules of rank $d$ with
a level $\bK$-structure such that the image of
$\cF_{j+i+1}\to \cF_j$ equals the image of
the canonical injection 
$\cF_{j}(-\infty)\inj \cF_j$ for every $j\in \Z$.

Two chains $\cdots \inj \cF_{1} \inj \cF_0 \inj
\cF_{-1} \inj \cdots$ and 
$\cdots \inj \cF'_{1} \inj \cF'_0 \inj
\cF'_{-1} \inj \cdots$ are equivalent if and only if
there exists an integer $l$ and an isomorphism
$\cF_{j} \cong \cF'_{j+l}$ of $\cO_C$-modules with
level structures for every $j\in \Z$ such that 
the diagram
$$
\begin{CD}
\cdots @>>> \cF_{1} @>>> \cF_0 @>>> \cF_{-1} @>>> \cdots \\
@. @V{\cong}VV @V{\cong}VV @V{\cong}VV @. \\
\cdots @>>> \cF'_{l+1} @>>> \cF'_l @>>> \cF'_{l-1} @>>> \cdots
\end{CD}
$$
is commutative.

\subsection{the Harder-Narasimhan polygons}
\label{sec:HN}
We use functions 
$\Delta_{p_\cF}$
from the theory of Harder-Narasimhan polygons.
Let us recall the properties we use below.

Let $\cF$ be a locally free $\cO_C$-module of 
rank $r$. For an $\cO_C$-submodule $\cF' \subset \cF$
(note that $\cF'$ is automatically locally free),
we set $z_{\cF}(\cF') = (\rank(\cF'),\deg(\cF'))
\in \Q^2$. It is known that there exists a unique convex, 
piecewise affine,
affine on $[i-1,i]$ for
$i=1,\ldots, r$, 
continuous function 
$p_{\cF} :[0,r] \to \R$ on the interval
$[0,r]$ such that the 
convex hull of the set $\{z_{\cF}(\cF')\ |\ 
\cF' \subset \cF\}$ in $\R^2$ equals
$\{(x,y)\ |\ 0\le x\le r,\, y\le p_{\cF}(x) \}$.
We define the function $\Delta p_{\cF}:
\{1,\ldots, r-1\} \to \R$ as
$\Delta p_{\cF}(i) = 2 p_{\cF}(i) - p_{\cF}(i-1)
- p_{\cF}(i+1)$. Then $\Delta p_{\cF}(i) \ge 0$
for all $i$. We note that for an invertible $\cO_C$-module
$\cL$, $\Delta p_{\cF \otimes \cL}$ equals $\Delta p_{\cF}$.
The theory of Harder-Narasimhan
filtration (\cite{HN}) implies that, 
if $i \in \Supp(\Delta p_{\cF})
= \{ i\ |\ \Delta p_{\cF}(i)>0 \}$, 
then there exists a
unique $\cO_C$-submodule $\cF' \subset \cF$
satisfying $z_{\cF}(\cF')=(i,p_{\cF}(i))$. 
We denote this $\cO_C$-submodule $\cF'$ by $\cF_{(i)}$.
The submodule $\cF_{(i)}$ has the following properties.
\begin{itemize}
\item If $i, j \in \Supp(\Delta p_{\cF})$ with $i\le j$,
then $\cF_{(i)}\subset \cF_{(j)}$ and $\cF_{(j)}/\cF_{(i)}$
is locally free.
\item If $i \in \Supp(\Delta p_{\cF})$, then
$p_{\cF_{(i)}}(x) = p_{\cF}(x)$
for $x \in [0,i]$ and $p_{\cF/\cF_{(i)}}(x-i)
=p_{\cF}(x) - \deg(\cF_{(i)})$ for 
$x \in [i,r]$.
\end{itemize}

\begin{lem}\label{7_diffdeg}
Let $\cF$ be a locally free
$\cO_C$-module of finite rank, and
let $\cF'\subset \cF$ be a $\cO_C$-submodule 
of the same rank. 
Then we have 
$0 \le p_{\cF}(i) -p_{\cF'}(i) \le 
\deg(\cF)-\deg(\cF')$
for $i =1,\ldots,\rank(\cF)-1$.
\end{lem}
\begin{proof}
Immediate from the definition of $p_{\cF}$.
\end{proof}

\begin{lem}\label{7_HNdiff}
Let $\cF$ be a locally free $\cO_C$-module
of rank $d$. Let $\cF' \subset \cF$ be
a n $\cO_C$-submodule of the same rank. 
Suppose that 
$\Delta p_{\cF}(i) > \deg(\cF) 
-\deg(\cF')$. Then we have 
$\cF'_{(i)} = \cF_{(i)} \cap \cF'$.
\end{lem}

\begin{proof}
It suffices to prove that
$\cF'_{(i)}\subset \cF_{(i)}$.
Assume otherwise. Let us consider
the short exact sequence
$$
0\to \cF'_{(i)}\cap
\cF_{(i)} \to \cF'_{(i)} \to
\cF'_{(i)}/(\cF'_{(i)}\cap
\cF_{(i)}) \to 0
$$
Let $r$ denote the rank of
$\cF'_{(i)}\cap \cF_{(i)}$. 
By assumption, $r$ is strictly
smaller than $i$. Hence
$$
\begin{array}{rl}
\deg(\cF'_{(i)}) & =
\deg(\cF'_{(i)}\cap
\cF_{(i)}) + 
\deg(\cF'_{(i)}/(\cF'_{(i)}\cap
\cF_{(i)})) \\
& \le p_{\cF}(r) + 
p_{\cF/\cF_{(i)}}(i-r) \\
& \le p_{\cF}(i) - (i-r)(p_{\cF}(i)-p_{\cF}(i-1))
+ (i-r)(p_{\cF}(i+1) -p_{\cF}(i)) \\
& = \deg(\cF_{(i)}) - (i-r) \Delta p_{\cF}(i) \\
& < \deg(\cF_{(i)}) - (\deg(\cF)-\deg(\cF')).
\end{array}
$$
On the other hand, Lemma~\ref{7_diffdeg}
shows that $\deg(\cF_{(i)}\cap \cF')\ge
\deg(\cF_{(i)}) - (\deg(\cF)-\deg(\cF'))$.
This is a contradiction.
\end{proof}

\section{Some subsimplicial complexes}
\label{sec:Xsigma}
Using the functions $\Delta_{p_\cF}$ defined 
above, we introduce some subsimplicial complexes of $X_{\bK,\bullet}$
and of $\wt{X}_{\bK,\bullet}$.
In \cite{Gra}, Grayson considered subsimplicial 
complexes of the building but here
we also consider subsimplicial complexes of the quotients.
We introduce three spaces.

\subsection{}
Given a subset $\cD \subset \{1,\ldots,d-1\}$
and a real number $\alpha > 0$,
we define the simplicial subcomplex
$X_{\bK,\bullet}^{(\alpha),\cD}$ of
$X_{\bK,\bullet}$ as follows:
A simplex of $X_{\bK,\bullet}$
belongs to $X_{\bK,\bullet}^{(\alpha),\cD}$
if and only if each of its vertices is
represented by a 
locally free $\cO_C$-module $\cF$
of rank $d$ with a level $\bK$-structure
such that $\Delta p_{\cF}(i) \ge \alpha$
holds for every $i \in \cD$.

Let $X_{\bK,\bullet}^{(\alpha)}$ denote
the union $X_{\bK,\bullet}^{(\alpha)}
= \bigcup_{\cD \neq \emptyset} X_{\bK,\bullet}^{(\alpha),\cD}$.

\subsection{}
Write $\cD=\{i_1, \dots, i_{r-1}\}$
with $i_1<\cdots <i_{r-1}$.
Let $\Flag_{\cD}$ denote the set
$$
\Flag_{\cD} = \{ f=
[0 \subset V_1 \subset \cdots \subset
V_{r-1} \subset F^{\oplus d}] \ |\ \dim(V_j)=i_j \}
$$
of flags in $F^{\oplus d}$.

Let $\wt{X}_{\bK,\bullet}^{(\alpha),\cD}$ denote
the inverse image of $X_{\bK,\bullet}^{(\alpha),\cD}$
by the morphism $\wt{X}_{\bK,\bullet}
\to X_{\bK,\bullet}$. 
For $f = [0 \subset V_1 \subset \cdots \subset
V_{r-1} \subset F^{\oplus d}]\in \Flag_{\cD}$, let 
$\wt{X}_{\bK,\bullet}^{(\alpha),\cD,f}$ 
denote the simplicial subcomplex of 
$\wt{X}_{\bK,\bullet}^{(\alpha),\cD}$
consisting of the simplices in $\wt{X}_{\bK,\bullet}$ 
whose representative
$\cdots \supsetneqq \cF_{-1} \supsetneqq 
\cF_0 \supsetneqq \cF_1
\supsetneqq \cdots$ satisfies
$\cF_{l,(i_j)} = \cF_l \cap \eta_* V_{i_j}$
for every $l\in\Z$, $j=1,\ldots,r-1$.
Lemma~\ref{7_HNdiff} implies that, for 
$\alpha > (d-1)\deg(\infty)$,
$\wt{X}_{\bK,\bullet}^{(\alpha),\cD}$ is decomposed
into a disjoint union $\wt{X}_{\bK,\bullet}^{(\alpha),\cD}
= \coprod_{f \in \Flag_{\cD}}
\wt{X}_{\bK,\bullet}^{(\alpha),\cD,f}$.

\subsection{}
For $g \in \GL_d(\A^\infty)$, we set 
$$
\wt{Y}_{\bK,\bullet}^{(\alpha),\cD,g} =
\wt{X}_{\bK,\bullet}^{(\alpha),\cD,f_0}
\cap (P_{\cD}(\A^\infty)g/(g^{-1}P_{\cD}(\A^\infty)g \cap \bK)
\times \cBT_{\bullet})
$$
and $Y_{\bK,\bullet}^{(\alpha),\cD,g} = 
P_{\cD}(F)\bsl \wt{Y}_{\bK,\bullet}^{(\alpha),\cD,g}$.
We omit the superscript $g$ on 
$\wt{Y}_{\bK,\bullet}^{(\alpha),\cD,g}$
and $Y_{\bK,\bullet}^{(\alpha),\cD,g}$ if $g=1$.
If we take a complete set $T \subset \GL_d(\A^\infty)$ 
of representatives of $P_{\cD}(\A^\infty)\bsl \GL_d(\A^\infty)$, 
then we have $\wt{X}_{\bK,\bullet}^{(\alpha),\cD,f_0} = 
\coprod_{g \in T} \wt{Y}_{\bK,\bullet}^{(\alpha),\cD,g}$.

\section{the finite adele actions}
\begin{lem}\label{7_betag}
For every $g \in \GL_d(\A^\infty)$
satisfying $g^{-1}\bK g \subset \GL_d(\wh{A})$,
there exists a real number $\beta_g \ge 0$
such that the isomorphism $\xi_g:X_{\bK,\bullet}
\xto{\cong} X_{g^{-1}\bK g,\bullet}$ sends
$X_{\bK,\bullet}^{(\alpha),\cD}$ to
$X_{g^{-1} \bK g,\bullet}^{(\alpha -\beta_g),\cD}
\subset X_{g^{-1}\bK g,\bullet}$
%
for all $\alpha > \beta_g$, 
and for 
all nonempty subset $\cD \subset \{1,\ldots,d-1 \}$.

\end{lem}
\begin{proof}
Take two elements $a,b \in \A^{\infty \times}
\cap \wh{A}$ such that both $ag$ and
$bg^{-1}$ lie in $\GL_d(\A^\infty)\cap
\Mat_d(\wh{A})$. Then for any
$h \in \GL_d(\A^\infty)$ we have 
$a
\wh{A}^{\oplus d}h^{-1} 
\subset \wh{A}^{\oplus d}g^{-1} h^{-1} \subset b^{-1} \wh{A}^{\oplus d}h^{-1}$.
This implies that,
for any vertex $x \in X_{\bK,0}$,
if we take suitable representatives 
$\cF_x$, $\cF_{\xi_{g}(x)}$ of the equivalence classes
of locally free $\cO_C$-modules 
corresponding to $x$, $\xi_{g}(x)$, then
there exists a sequence of injections 
$\cF_x(-\divi(a)) \inj \cF_{\xi_g(x)}
\inj \cF_x(\divi(b))$. Applying Lemma~\ref{7_diffdeg},
we see that there exists a positive real number $m_g>0$
not depending on $x$ such that 
$|p_{\cF_x}(i) - p_{\cF_{\xi_{g}(x)}}(i)| < m_g$ for all $i$.
Hence the claim follows.
\end{proof}

\begin{remark}
Any open compact subgroup of $\GL_d(\A^\infty)$
is conjugate to an open subgroup of $\GL_d(\wh{A})$.
The set of open subgroups
of $\GL_d(\wh{A})$ is cofinal in the inductive
system of all open compact subgroups of $\GL_d(\A^\infty)$.
Therefore, to prove Theorem~\ref{7_PROP1}, we may
without loss of generality assume that 
the group $\bK$ is contained in
$\GL_d(\wh{A})$, and we may
replace the inductive limit $\varinjlim_{\bK}$ in the definition of
$H_{d-1}^\mathrm{BM}(X_{\lim,\bullet},M)$ and $H_{d-1}(X_{\lim,\bullet},M)$
with the inductive limit $\varinjlim_{\bK\subset \GL_d(\wh{A})}$.
\end{remark}
From now on until the end of this section,
we exclusively deal with the subgroups 
$\bK \subset \GL_d(\A^\infty)$ contained in $\GL_d(\wh{A})$.
The notation $\varinjlim_{\bK}$ henceforth means
the inductive limit $\varinjlim_{\bK\subset \GL_d(\wh{A})}$.

Thus the group $\GL_d(\A^\infty)$ acts on 
$\varinjlim_{\bK} \varinjlim_{\alpha > 0} 
H^{*}(X^{(\alpha),\cD}_{\bK,\bullet},\Q)$
in such a way that the exact sequence (\ref{7_cohseq})
is $\GL_d(\A^\infty)$-equivariant.

\subsection{}
An argument similar to that in the proof of
Lemma~\ref{7_betag} shows that, 
for each $g \in \GL_d(\A^\infty)$ satisfying
$g^{-1}\bK g \subset \GL_d(\wh{A})$, there exists
a real number $\beta'_g > \beta_g$ such that
the isomorphism $\wt{\xi}_g$
sends $\wt{X}_{\bK,\bullet}^{(\alpha),\cD,f}$ to
$\wt{X}_{g \bK g^{-1},\bullet}^{(\alpha-\beta_g),\cD,f}
\subset \wt{X}_{g \bK g^{-1},\bullet}$
for $\alpha > \beta'_g$ and for any $f \in \Flag_{\cD}$.

\subsection{}
For $\gamma \in \GL_d(F)$, the action of
$\gamma$ on $\wt{X}_{\bK,\bullet}$
sends $\wt{X}_{\bK,\bullet}^{(\alpha),\cD,f}$ 
bijectively to 
$\wt{X}_{\bK,\bullet}^{(\alpha),\cD,\gamma f}$.
Let $f_0 = [0\ \subset F^{\oplus i_1}\oplus 
\{0\}^{\oplus d-i_1} 
\subset \cdots \subset F^{\oplus i_{r-1}}
\oplus \{0\}^{\oplus d -i_{r-1}}
\subset F^{\oplus d}]\in \Flag_{\cD}$ be
the standard flag. The group $\GL_d(F)$
acts transitively on $\Flag_{\cD}$ and its
stabilizer at $f_0$ equals $P_{\cD}(F)$.
Hence for $\alpha > (d-1)\deg(\infty)$,
$X_{\bK,\bullet}^{(\alpha),\cD}$ is
isomorphic to the quotient 
$P_{\cD}(F)\bsl 
\wt{X}_{\bK,\bullet}^{(\alpha),\cD,f_0}$.

For $g \in \GL_d(\A^\infty)$, we set 
$$
\wt{Y}_{\bK,\bullet}^{(\alpha),\cD,g} =
\wt{X}_{\bK,\bullet}^{(\alpha),\cD,f_0}
\cap (P_{\cD}(\A^\infty)g/(g^{-1}P_{\cD}(\A^\infty)g \cap \bK)
\times \cBT_{\bullet})
$$
and $Y_{\bK,\bullet}^{(\alpha),\cD,g} = 
P_{\cD}(F)\bsl \wt{Y}_{\bK,\bullet}^{(\alpha),\cD,g}$.
We omit the superscript $g$ on 
$\wt{Y}_{\bK,\bullet}^{(\alpha),\cD,g}$
and $Y_{\bK,\bullet}^{(\alpha),\cD,g}$ if $g=1$.

\newcommand{\gbar}{\overline{g}}
We note that, when $\bK$, $\alpha$ and $\cD$ are fixed, 
$\wt{Y}_{\bK,\bullet}^{(\alpha),\cD,g}$ and
$Y_{\bK,\bullet}^{(\alpha),\cD,g}$ 
depend only on the class $\gbar = P_{\cD}(\A^\infty)g$ 
of $g$ in $P_{\cD}(\A^\infty)\bsl \GL_d(\A^\infty)$.
By abuse of notation, we denote 
$\wt{Y}_{\bK,\bullet}^{(\alpha),\cD,g}$ and
$Y_{\bK,\bullet}^{(\alpha),\cD,g}$ 
by $\wt{Y}_{\bK,\bullet}^{(\alpha),\cD,\gbar}$ and
$Y_{\bK,\bullet}^{(\alpha),\cD,\gbar}$, respectively.
Then we have $\wt{X}_{\bK,\bullet}^{(\alpha),\cD,f_0} = 
\coprod_{\gbar \in P_{\cD}(\A^\infty)\bsl \GL_d(\A^\infty)} 
\wt{Y}_{\bK,\bullet}^{(\alpha),\cD,\gbar}$.
Hence we have 
\begin{equation} \label{eq:decomposition}
X_{\bK,\bullet}^{(\alpha),\cD} = 
\coprod_{\gbar \in P_{\cD}(\A^\infty)\bsl \GL_d(\A^\infty)} 
Y_{\bK,\bullet}^{(\alpha),\cD,\gbar}
\end{equation}
for $\alpha > (d-1) \deg(\infty)$.

\subsection{}
We use a covering spectral sequence 
\begin{equation}\label{7_specseq}
E_1^{p,q} = \bigoplus_{\sharp \cD = p+1}
H^q(X^{(\alpha),\cD}_{\bK,\bullet},\Q)
\Rightarrow H^{p+q}(X^{(\alpha)}_{\bK,\bullet},\Q)
\end{equation}
with respect to
the covering $X_{\bK,\bullet}^{(\alpha)} = \bigcup_{1\le i \le d-1}
X_{\bK,\bullet}^{(\alpha),\{i\}}$ of $X_{\bK,\bullet}^{(\alpha)}$. 
For $\alpha' \ge \alpha>0$, the inclusion 
$X_{\bK,\bullet}^{(\alpha),\cD} \to
X_{\bK,\bullet}^{(\alpha'),\cD}$
induces a morphism of spectral sequences.
Taking the inductive limit, we obtain the 
spectral sequence 
$$
E^{p,q}_1
= \bigoplus_{\sharp \cD = p+1}
\varinjlim_{\alpha} H^q(X^{(\alpha),\cD}_{\bK,\bullet},\Q)
\Rightarrow
\varinjlim_{\alpha} H^{p+q}(X^{(\alpha)}_{\bK,\bullet},\Q).
$$
For $g \in \GL_d(\A^\infty)$ satisfying
$g^{-1}\bK g \subset \GL_d(\wh{A})$, let $\beta_g$ be as
in Lemma~\ref{7_betag}. Then for $\alpha >\beta_g$ 
the isomorphism $\xi_g :X_{\bK,\bullet}\xto{\cong}
X_{g\bK g^{-1},\bullet}$ induces a homomorphism
from the spectral sequence (\ref{7_specseq}) for 
$X_{\bK,\bullet}^{(\alpha)}$ to that 
for $X_{\bK,\bullet}^{(\alpha -\beta_g)}$. 
Passing to the inductive limit with respect to $\alpha$
and then passing to the inductive limit with respect to 
$\bK$, we obtain the left action of the group 
$\GL_d(\A^\infty)$ on the spectral sequence
\begin{equation}\label{7_limspecseq}
E^{p,q}_1
= \bigoplus_{\sharp \cD = p+1}
\varinjlim_{\bK} 
\varinjlim_{\alpha} H^q(X^{(\alpha),\cD}_{\bK,\bullet},\Q)
\Rightarrow
\varinjlim_{\bK} 
\varinjlim_{\alpha} H^{p+q}(X^{(\alpha)}_{\bK,\bullet},\Q).
\end{equation}

\section{Finiteness and an application}
We prove that the complement of 
the boundary is finite.  
Then we express the Borel-Moore
homology and cohomology with 
compact support as a limit of 
relative homology and cohomology
respectively.
\begin{lem}\label{7_lem:finiteness}
For any $\alpha >0$, the set of the simplices 
in $X_{\bK,\bullet}$ not belonging to 
$X_{\bK,\bullet}^{(\alpha)}$ is finite.
\end{lem}
\begin{proof}
Let $\cP$ denote the set of continuous,
convex functions $p':[0,d]\to \R$ with
$p'(0)=0$ such that $p'(i)\in \Z$ and 
$p'$ is affine on $[i-1,i]$ for $i=1,\ldots,d$.
It is known that 
for any $r \ge 1$ and $f \in \Z$,
there are only a finite number of isomorphism classes
of semi-stable locally free $\cO_C$-modules of 
rank $r$ with degree $f$. Hence by the theory of
Harder-Narasimhan filtration, for any $p' \in \cP$,
the set of the isomorphism classes
of locally free $\cO_C$-modules $\cF$ with
$p_{\cF} = p'$ is finite.
Let us give an action of the group $\Z$ on the 
set $\cP$, by setting $(a\cdot p')(x)=p'(x)+ a\deg(\infty)x$
for $a \in \Z$ and for $p' \in \cP$.
Then $p_{\cF(a \infty)}= a\cdot p_{\cF}$ for
any $a \in \Z$ and for any locally free $\cO_C$-module
$\cF$ of rank $d$. For $\alpha >0$ let 
$\cP^{(\alpha)} \subset \cP$ denote the set 
of functions $p' \in \cP$ 
with $2p'(i)- p'(i-1)-p'(i+1) \le \alpha$
for each $i \in \{1,\ldots, d-1\}$.
An elementary argument shows
that the quotient $\cP^{(\alpha)}/\Z$ is a finite set,
whence the claim follows.
\end{proof}

\subsection{}
\label{sec:limit BM isom}
Lemma~\ref{7_lem:finiteness} implies that 
$H_{d-1}^\mathrm{BM}(X_{\bK,\bullet},\Q)$ is canonically
isomorphic to the projective limit
$\varprojlim_{\alpha >0}
H_{d-1}(X_{\bK,\bullet},X_{\bK,\bullet}^{(\alpha)}; \Q)$
and $H_c^{d-1}(X_{\bK,\bullet},\Q)$ is canonically
isomorphic to the inductive limit
$\varinjlim_{\alpha >0}
H^{d-1}(X_{\bK,\bullet},X_{\bK,\bullet}^{(\alpha)}; \Q)$.
Thus from the (usual) long exact sequence of relative homology,
we have an exact sequence 
\begin{equation}\label{7_cohseq}
\varinjlim_{\alpha >0} 
H^{d-2}(X^{(\alpha)}_{\bK,\bullet},\Q)
\to H_c^{d-1}(X_{\bK,\bullet},\Q)
\to H^{d-1}(X_{\bK,\bullet},\Q)
\to \varinjlim_{\alpha >0} 
H^{d-1}(X^{(\alpha)}_{\bK,\bullet},\Q).
\end{equation}

\section{Some isomorphisms}
\begin{prop}\label{7_PROP1b}
For $\alpha' \ge \alpha >(d-1)\deg(\infty)$, 
the homomorphism
$H^*(X^{(\alpha)}_{\bK,\bullet},\Q)
\to H^*(X^{(\alpha')}_{\bK,\bullet},\Q)$
is an isomorphism.
\end{prop}

\begin{lem}\label{7_contract}
For any $g \in \GL_d(\A^\infty)$, 
the simplicial complex
$\wt{X}_{\bK,\bullet}^{(\alpha),\cD,f_0}
\cap (\{ g\bK \}\times \cBT_{\bullet})$ is 
non-empty and contractible.
\end{lem}
\begin{proof}
Since $\wt{X}_{\bK,\bullet}^{(\alpha),\cD,f_0}
\cap (\{ g\bK \}\times \cBT_{\bullet})$
is isomorphic to $\wt{X}_{\GL_d(\wh{A}),\bullet}^{(\alpha),\cD,f_0}
\cap (\{ g\GL_d(\wh{A}) \}\times \cBT_{\bullet})$,
we may assume that $\bK = \GL_d(\wh{A})$.
We set $X=\wt{X}_{\GL_d,\GL_d(\wh{A}),\bullet}^{(\alpha),\cD,f_0}
\cap (\{ g\GL_d(\wh{A}) \}\times \cBT_{\GL_d,\bullet})$.

We proceed by induction on $d$, 
in a manner similar to that in the proof of Theorem~4.1 of 
\cite{Gra}.
Let $i\in \cD$ be the minimal element and set $d'=d-i$.
We define the subset $\cD' \subset \{1,\ldots,d'-1 \}$
as $\cD' = \{ i' -i\ |\ i' \in \cD, i'\neq i \}$.
We define $f'_0 \in \Flag_{\cD'}$ as the image of
the flag $f_0$ in $F^{\oplus d}$ with respect to the
the projection $F^{\oplus d} \surj 
F^{\oplus d}/(F^{\oplus i}\oplus \{0\}^{\oplus d'})
\cong F^{\oplus d'}$.
Take an element $g' \in \GL_{d'}(\A^\infty)$ such that
the quotient $\wh{A}^{\oplus d} g^{-1}/
(\wh{A}^{\oplus d}g^{-1}\cap 
(\A^{\infty \oplus i} \oplus \{0\}^{\oplus d'}))$
equals $\wh{A}^{\oplus d'}g^{\prime -1}$ as
an $\wh{A}$-lattice of $\A^{\infty \oplus d'}$.
We set $X'=\wt{X}_{\GL_{d'}, 
\GL_{d'}(\wh{A}),\bullet}^{(\alpha),\cD',f'_0}
\cap (\{ g' \GL_{d'}(\wh{A}) \}\times \cBT_{\GL_{d'},\bullet})$
if $\cD'$ is non-empty. Otherwise we set $X'
= \wt{X}_{\GL_{d'}, \GL_{d'}(\wh{A}),\bullet}
\cap (\{ g' \GL_{d'}(\wh{A}) \}\times \cBT_{\GL_{d'},\bullet})$.
By induction hypothesis, $|X'|$ is contractible.
There is a canonical morphism
$h:X \to X'$ which sends an $\cO_C$-submodule $\cF[g,L_\infty]$ of
$\eta_* F^{\oplus d}$ to the 
$\cO_C$-submodule $\cF[g,L_\infty]/\cF[g,L_\infty]_{(i)}$
of $\eta_* F^{\oplus d'}$.
Let $\epsilon : \Vertex(X) \to \Z$ and
$\epsilon' : \Vertex(X') \to \Z$ denote the maps that
send a locally free $\cO_C$-module $\cF$ to the integer
$[p_{\cF}(1)/\deg(\infty)]$.
We fix an
$\cO_C$-submodule $\cF_0$ of $\eta_* F^{\oplus d}$
whose equivalence class belongs to $X$.
By twisting $\cF_0$ by some power of $\cO_C(\infty)$
if necessary, we may assume that 
$p_{\cF_0}(i) - p_{\cF_0}(i-1) > \alpha$.
We fix a splitting $\cF_0 = \cF_{0,(i)}\oplus \cF'_0$.
This splitting induces an isomorphism
$\varphi : \eta_* \eta^* \cF'_0 \cong \eta_* F^{\oplus d}$.
Let $h' : X' \to X$ denote the morphism that sends
an $\cO_C$-submodule $\cF'$ of
$\eta_* \eta^* F^{\oplus d'}$ to the
$\cO_C$-submodule $\cF_{0,(i)}(\epsilon'(\cF')\infty) 
\oplus \varphi^{-1}(\cF')$ of $\eta_* F^{\oplus d}$.
For each $n \in \Z$, define a morphism
$G_n : X \to X$ by
sending an $\cO_C$-submodule $\cF$ of
$\eta_* \eta^* F^{\oplus d}$ to the
$\cO_C$-submodule $\cF_{0,(i)}((n+\epsilon(\cF))\infty) 
+ \cF$ of $\eta_* F^{\oplus d}$. Then the
argument in~\cite[p. 85--86]{Gra} shows that
$f$ and $|h'|\circ|h| \circ f$ are homotopic
for any map $f:Z \to |X|$ from a compact space $Z$
to $|X|$. Since the map $|h'|\circ|h| \circ f$ factors
through the contractible space $|X'|$, $f$ is
null-homotopic. Hence $|X|$ is contractible.
\end{proof}

\begin{proof}[Proof of Proposition~\ref{7_PROP1b}]
For any simplex $\sigma$ in $\wt{X}_{\bK,\bullet}$,
the stabilizer group $\Gamma_\sigma \subset \GL_d(F)$
is finite, as remarked in Section~\ref{sec:finite stab}.
Hence by Lemma~\ref{7_contract}, both 
$H^*(Y_{\bK,\bullet}^{(\alpha),\cD,g},\Q)$ and
$H^*(Y_{\bK,\bullet}^{(\alpha'),\cD,g},\Q)$ are canonically
isomorphic to the same group $H^*(P_{\cD}(F),
\Map(P_{\cD}(\A)g/(g^{-1}P_{\cD}(\A^\infty)g \cap \bK),\Q))$ 
for any non-empty subset $\cD \subset \{1,\ldots, d-1\}$
and for $g \in \GL_d(\A^\infty)$.
This shows that $H^*(Y_{\bK,\bullet}^{(\alpha),\cD,g},\Q)
\to H^*(Y_{\bK,\bullet}^{(\alpha'),\cD,g},\Q)$ is an isomorphism.
One can decompose, by using \eqref{eq:decomposition},
the $\Q$-vector space $H^*(X_{\bK,\bullet}^{(\alpha),\cD},\Q)$
as the product
\begin{equation} \label{eq:decomposition2}
H^*(X_{\bK,\bullet}^{(\alpha),\cD},\Q)
= \prod_{\gbar \in P_{\cD}(\A^\infty)\bsl \GL_d(\A^\infty)} 
H^*(Y_{\bK,\bullet}^{(\alpha),\cD,\gbar},\Q).
\end{equation}
Hence the homomorphism 
$H^*(X_{\bK,\bullet}^{(\alpha),\cD},\Q)
\to H^*(X_{\bK,\bullet}^{(\alpha'),\cD},\Q)$ is an isomorphism.
\end{proof}

\section{Proof of Theorem \ref{7_PROP1}}

\subsection{}
For a subset $\cD$ of $\{1,\ldots,d-1\}$,
we define the algebraic
groups $P_{\cD}$, $N_{\cD}$ and
$M_{\cD}$ as follows. We write 
$\cD= \{i_1,\ldots,i_{r-1} \}$, with
$i_0=0 < i_1 < \cdots < i_{r-1} <i_r=d$ 
and set $d_j = i_j -i_{j-1}$ for
$j=1,\ldots,r$.
We define $P_{\cD}$, $N_{\cD}$ and
$M_{\cD}$ as the standard parabolic
subgroup of $\GL_d$ of type
$(d_1,\ldots,d_r)$, the unipotent radical
of $P_{\cD}$, and the quotient group $P_{\cD}/N_{\cD}$
respectively. We identify the group $M_{\cD}$ 
with $\GL_{d_1}\times\cdots \times \GL_{d_r}$.

Let us consider the smooth $\GL_d(\A^\infty)$-module,
$\varinjlim_{\bK} \varinjlim_{\alpha}
H^*(X_{\bK,\bullet}^{(\alpha),\cD},\Q)$.
For a fixed $\bK$, we have
\begin{equation} \label{eq:decomposition3}
\varinjlim_{\alpha}
H^*(X_{\bK,\bullet}^{(\alpha),\cD},\Q)
= \prod_{\gbar \in \gbar \in P_{\cD}(\A^\infty)\bsl \GL_d(\A^\infty)} 
\varinjlim_{\alpha}
H^*(Y_{\bK,\bullet}^{(\alpha),\cD,\gbar},\Q).
\end{equation}
since $H^*(Y_{\bK,\bullet}^{(\alpha),\cD,\gbar},\Q)
\to H^*(Y_{\bK,\bullet}^{(\alpha'),\cD,\gbar},\Q)$ is an isomorphism
for $\alpha' \ge \alpha > (d-1)\deg(\infty)$.
We note that
$\varinjlim_{\bK} \varinjlim_{\alpha}
H^*(Y_{\bK,\bullet}^{(\alpha),\cD,g},\Q)$
is a smooth $g^{-1} P_\cD(\A^\infty) g$-module
for any $g \in \GL_d(\A^\infty)$.
Via \eqref{eq:decomposition3} we regard
$\varinjlim_{\bK} \varinjlim_{\alpha}
H^*(X_{\bK,\bullet}^{(\alpha),\cD},\Q)$
as a submodule of
$\prod_{\gbar \in P_{\cD}(\A^\infty)\bsl \GL_d(\A^\infty)} 
\varinjlim_{\bK} \varinjlim_{\alpha}
H^*(Y_{\bK,\bullet}^{(\alpha),\cD,\gbar},\Q)$.

Let $g \in \GL_d(\A^\infty)$.
For an open compact subgroup
$\bK \subset \GL_d(\wh{A})$
satisfying $g^{-1}\bK g \subset \GL_d(\wh{A})$,
there exists
a real number $\beta'_g > \beta_g$ such that
the isomorphism $\wt{\xi}_g$
sends $\wt{Y}_{\bK,\bullet}^{(\alpha),\cD,g'}$ 
to $\wt{Y}_{g \bK g^{-1},\bullet}^{(\alpha-\beta_g),\cD,g'g}
\subset \wt{X}_{g \bK g^{-1},\bullet}$
for $\alpha > \beta'_g$, for any $f \in \Flag_{\cD}$
and for any $g' \in \GL_d(\A^\infty)$.
This induces a morphism
$\xi_{g,g',\bK}: Y_{\bK,\bullet}^{(\alpha),\cD,g'} 
\to Y_{g \bK g^{-1},\bullet}^{(\alpha-\beta_g),\cD,g'g}$
of (generalized) simplicial complexes.
By varying $\alpha$, we have a homomorphism
$$
\xi_{g,g',\bK}^* :
\varinjlim_{\alpha}
H^*(Y_{g \bK g^{-1},\bullet}^{(\alpha),\cD,g'g},\Q)
\to 
\varinjlim_{\alpha}
H^*(Y_{\bK,\bullet}^{(\alpha),\cD,g'},\Q).
$$
The homomorphism $\xi_{g,g',\bK}^*$ is an isomorphism
since the homomorphism
$\xi_{g^{-1},g'g,g \bK g^{-1}}^*$ 
gives its inverse.
By varying $\bK$, we obtain an isomorphism
$$
\xi_{g,g'}^* :
\varinjlim_{\bK} \varinjlim_{\alpha}
H^*(Y_{\bK,\bullet}^{(\alpha),\cD,g'g},\Q)
\xto{\cong}
\varinjlim_{\bK} \varinjlim_{\alpha}
H^*(Y_{\bK,\bullet}^{(\alpha),\cD,g'},\Q).
$$

Let us choose a complete set $T \subset \GL_d(\A^\infty)$
of representatives such that $1 \in T$.
Then $\varinjlim_{\bK} \varinjlim_{\alpha}
H^*(X_{\bK,\bullet}^{(\alpha),\cD},\Q)$
is a submodule of
$\prod_{g \in T} 
\varinjlim_{\bK} \varinjlim_{\alpha}
H^*(Y_{\bK,\bullet}^{(\alpha),\cD,g},\Q)$.
We set 
$$
H^{*,\cD}_\Q = 
\varinjlim_{\bK} \varinjlim_{\alpha}
H^*(Y_{\bK,\bullet}^{(\alpha),\cD},\Q).
$$
The isomorphism $\xi_{g,1}^* :
\varinjlim_{\bK} \varinjlim_{\alpha}
H^*(Y_{\bK,\bullet}^{(\alpha),\cD,g},\Q)
\xto{\cong}H^{*,\cD}_\Q$
for each $g \in T$ gives an
isomorphism
$$
\prod_{g \in T} 
\varinjlim_{\bK} \varinjlim_{\alpha}
H^*(Y_{\bK,\bullet}^{(\alpha),\cD,g},\Q)
\cong \prod_{g \in T} 
H^{*,\cD}_\Q.
$$
Via this isomorphism
we regard  $\varinjlim_{\bK} \varinjlim_{\alpha}
H^*(X_{\bK,\bullet}^{(\alpha),\cD},\Q)$
as a submodule of
$\prod_{g \in T} H^{*,\cD}_\Q$.

Let $g' \in \GL_d(\A^\infty)$ be an arbitrary element.
For each $g \in T$, let us write
$gg' = h_g g_1$ with $g_1 \in T$ and
$h_g \in P_\cD(\A^\infty)$.
Then, with respect to the inclusion
$\varinjlim_{\bK} \varinjlim_{\alpha}
H^*(X_{\bK,\bullet}^{(\alpha),\cD},\Q)
\inj
\prod_{g \in T} H^{*,\cD}_\Q$,
the action of $g'$ on
$\varinjlim_{\bK} \varinjlim_{\alpha}
H^*(X_{\bK,\bullet}^{(\alpha),\cD},\Q)$
is compatible with
the automorphism $\theta(g')$ of 
$\prod_{g \in T} H^{*,\cD}_\Q$
that sends $(x_g)_{g \in T}$
to $(\xi_{h_g,1}^*(x_{g_1}))_{g \in T}$.
The automorphisms $\theta(g')$ for various $g'$
give an action of $\GL_d(\A^\infty)$
on $\prod_{g \in T} H^{*,\cD}_\Q$
from the left and an element $x$ of
$\prod_{g \in T} H^{*,\cD}_\Q$
belongs to the
submodule $\varinjlim_{\bK} \varinjlim_{\alpha}
H^*(X_{\bK,\bullet}^{(\alpha),\cD},\Q)$
if and only if $x$ is invariant under
some open compact subgroup $\bK \subset \GL_d(\A^\infty)$.
Observe that the smooth part of
$\prod_{g \in T} H^{*,\cD}_\Q$
with respect to the action of $\GL_d(\A^\infty)$
introduced above
is equal to the (unnormalized) parabolic induction
$\Ind_{P_{\cD}(\A^\infty)}^{\GL_d(\A^\infty)}
\varinjlim_{\bK} \varinjlim_{\alpha} 
H^* (Y_{\bK,\bullet}^{(\alpha),\cD},\Q)$.
Thus we obtain an isomorphism
$$
\varinjlim_{\bK} \varinjlim_{\alpha}
H^*(X_{\bK,\bullet}^{(\alpha),\cD},\Q)
\cong \Ind_{P_{\cD}(\A^\infty)}^{\GL_d(\A^\infty)}
\varinjlim_{\bK} \varinjlim_{\alpha} 
H^* (Y_{\bK,\bullet}^{(\alpha),\cD},\Q).
$$
It is straightforward to check that this isomorphism
is independent of the choice of $T$.

\begin{prop}\label{7_prop2}
Let the notations be above.
Then as a smooth $\GL_d(\A^\infty)$-module,
$\varinjlim_{\bK} \varinjlim_{\alpha >0}
H^*(X_{\bK,\bullet}^{(\alpha),\cD},\Q)$ is
isomorphic to 
$$
\Ind_{P_{\cD}(\A^\infty)}^{\GL_d(\A^\infty)}
\bigotimes_{j=1}^r 
\varinjlim_{\bK_j \subset \GL_{d_j}(\wh{A})} 
H^*(X_{\GL_{d_j},\bK_j,\bullet},\Q),
$$
where the group $P_{\cD}(\A^\infty)$ acts on
${\displaystyle \bigotimes_{j=1}^r 
\varinjlim_{\bK_j\subset \GL_{d_j}(\wh{A})}} 
H^*(X_{\GL_{d_j},\bK_j,\bullet},\Q)$
via the quotient 
$P_{\cD}(\A^\infty) \to M_{\cD}(A^\infty)
= \prod_j \GL_{d_j}(\A^\infty)$,
and $\Ind_{P_{\cD}(\A^\infty)}^{\GL_d(\A^\infty)}$
denotes the parabolic induction 
unnormalized by the modulus function.
\end{prop}
We give a proof in the following subsections.

\subsection{}
For $j=1,\ldots,r$, let $\bK_j 
\subset \GL_{d_i}(\A^\infty)$ denote the
image of $\bK \cap P_{\cD}(\A^\infty)$ by the composition
$P_{\cD}(\A^\infty) \to M_{\cD}(\A^\infty)
\to \GL_{d_i}(\A^\infty)$.

We define the continuous map
$\wt{\pi}_{\cD,j}: 
|\wt{Y}_{\bK,\bullet}^{(\alpha),\cD}| \to 
|\wt{X}_{\GL_{d_j},\bK_j,\bullet}|$
of topological spaces in the following way.
Let $\sigma$ be an $i$-simplex in 
$\wt{Y}_{\bK,\bullet}^{(\alpha),\cD}$.
Take a chain 
$\cdots \supsetneqq \cF_{-1} \supsetneqq \cF_0
\supsetneqq \cF_1 \supsetneqq \cdots$
of $\cO_C$-modules representing $\sigma$.
For $l\in \Z$ we set $\cF_{l,j} = \cF_{l,(i_j)}/\cF_{l,(i_{j-1})}$,
which is an $\cO_C$-submodule of $\eta_* F^{\oplus d_j}$.
We set $S_j = \{ l \in \Z \ |\ \cF_{l,j}\neq \cF_{l+1,j} \}$.
Define the map $\psi_j : \Z \to S_j$ as
$\psi_j(l) = \min \{ l'\ge l\ |\ l' \in S_j \}$.
Take an order-preserving bijection 
$\varphi_j :S_j \xto{\cong} \Z$.
For $l \in\Z$ set $\cF'_{l} =\cF_{\varphi_j^{-1}(l), j}$.
Then the chain $\cdots \supsetneqq \cF'_{-1} \supsetneqq \cF'_0
\supsetneqq \cF'_1 \supsetneqq \cdots$ defines a
simplex $\sigma'$ in $\wt{X}_{\GL_{d_j},\bK_j,\bullet}$.
We define a continuous map $|\sigma| \to |\sigma'|$
as the affine map sending the vertex of $\sigma$ 
corresponding to $\cF_l$ to the vertex of $\sigma'$
corresponding to $\cF'_{\varphi_j \circ \psi_j(l)}$.
Gluing these maps, we obtain a continuous map
$\wt{\pi}_{\cD,j}: 
|\wt{Y}_{\bK,\bullet}^{(\alpha),\cD}| \to 
|\wt{X}_{\GL_{d_j},\bK_j,\bullet}|$.
We set $\wt{\pi}_{\cD} =
(\wt{\pi}_{\cD,1},\ldots,\wt{\pi}_{\cD,r})
:|\wt{Y}_{\bK,\bullet}^{(\alpha),\cD}| \to 
\prod_{j=1}^{r} 
|\wt{X}_{\GL_{d_j},\bK_j,\bullet}|$.
This continuous map descends to the 
continuous map 
$\pi_{\cD}: |Y_{\bK,\bullet}^{(\alpha),\cD}| \to 
\prod_{j=1}^{r} |X_{\GL_{d_j},\bK_j,\bullet}|$.

\subsection{}
If $g \in P_{\cD}(\A^\infty)$ and
$g^{-1}\bK g \subset \GL_d(\wh{A})$, then the isomorphism
$\xi_g :X_{\bK,\bullet} \xto{\cong}
X_{g^{-1}\bK g,\bullet}$ sends 
$Y_{\bK,\bullet}^{(\alpha),\cD}$ inside
$Y_{g^{-1}\bK g,\bullet}^{(\alpha -\beta_g),\cD}$.
If we denote by $(g_1,\ldots,g_r)$ 
the image of $g$ in $M_{\cD}(\A^\infty)
= \prod_{j=1}^r \GL_{d_j}(\A^\infty)$,
then the diagram
$$
\begin{CD}
|Y_{\bK,\bullet}^{(\alpha),\cD}| @>{\xi_g}>>
|Y_{g^{-1}\bK g,\bullet}^{(\alpha-\beta_g),\cD}| \\
@V{\pi_{\cD}}VV @V{\pi_{\cD}}VV \\
\prod_{j=1}^r |X_{\GL_{d_j},\bK_j,\bullet}|
@>{(\xi_{g_1},\ldots,\xi_{g_r})}>> \prod_{j=1}^r
|X_{\GL_{d_j},g_j^{-1} \bK_j g_j,\bullet}|
\end{CD}
$$
is commutative.

\subsection{}
With the notations as above, 
suppose that the open compact
subgroup $\bK \subset \GL_d(\A^\infty)$
has the following property.
\begin{equation} \label{7_property}
\text{the homomorphism }
P_{\cD}(\A^\infty)\cap \bK \to
\bK_1 \times \cdots \times \bK_r
\text{ is surjective.}
\end{equation}
For a simplicial complex $X$, we set
$I_X = \Map(\pi_0(X),\Q)$, where $\pi_0(X)$ is 
the set of the connected components
of $X$. Let us consider the following commutative
diagram.
\begin{equation}\label{7_CD}
\begin{CD}
H^*(M_{\cD}(F), 
\Map(\prod_{j=1}^r \pi_0(X_{\GL_{d_j},\bK_j,\bullet}),\Q))
@>>>
H^*(P_{\cD}(F), I_{Y_{\bK,\bullet}^{(\alpha),\cD}})
\\
@VVV @VVV \\
H^*_{M_{\cD}(F)}
(\prod_{j=1}^r |\wt{X}_{\GL_{d_j},\bK_j,\bullet}|,\Q)
@>>> H^*_{P_{\cD}(F)}(|Y_{\bK,\bullet}^{(\alpha),\cD}|,\Q)\\
@AAA @AAA \\
H^*(\prod_{j=1}^r |X_{\GL_{d_j},\bK_j,\bullet}|,\Q)
@>>> H^*(|Y_{\bK,\bullet}^{(\alpha),\cD}|,\Q).
\end{CD}
\end{equation}
Here $H^*_{M_{\cD}(F)}$ and $H^*_{P_{\cD}(F)}$ denote
the equivariant cohomology groups.

\subsection{}

\begin{prop}\label{7_prop3}
All homomorphisms in the above
diagram (\ref{7_CD}) are isomorphisms.
\end{prop}

\begin{proof}
We prove that the upper horizontal arrow
and the four vertical arrows are isomorphisms.

First we consider the upper horizontal arrow.
\begin{lem}\label{thislemma}
For $q\ge 1$, the group
$H^q (N_{\cD}(F), I_{Y_{\bK,\bullet}^{(\alpha),\cD}})$
is zero.
\end{lem}
\begin{proof}[Proof of Lemma \ref{thislemma}]
For each $x \in N_{\cD}(F) \bsl 
\pi_0(Y_{\bK,\bullet}^{(\alpha),\cD})$,
take a lift $\wt{x} \in \pi_0(Y_{\bK,\bullet}^{(\alpha),\cD})$
of $x$ and let $N_x \subset N_{\cD}(F)$ denote the stabilizer
of $\wt{x}$. Then the group
$H^*(N_{\cD}(F), I_{Y_{\bK,\bullet}^{(\alpha),\cD}})$
is isomorphic to the direct product 
$$\prod_{x \in N_{\cD}(F) \bsl 
\pi_0(Y_{\bK,\bullet}^{(\alpha),\cD})}
H^*(N_x ,\Q).$$ We note that the group $N_{\cD}(F)$ is
a union $N_{\cD}(F) = \bigcup_{i} U_i$ of 
finite subgroups of $p$-power order where $p$ is 
the characteristic of $F$.  (This follows easily from
\cite[p.2, 1.A.2 Lemma]{KeWe} or from
\cite[p.60, 1.L.1 Theorem]{KeWe}.)
Hence $N_x = \bigcup_{i} (U_i \cap N_x)$.

We claim $H^{j}(N_x,\Q) =0$ for $j\ge 1$.
Because the projective system 
of cochain complexes that 
compute $H^*(U_i \cap N_x, M)$,
where $M$ is a $\Q$-vector space,
satisfies the Mittag-Leffler condition,
the cohomology 
$H^{j}(N_x,\Q)$ equals the 
projective limit $\varprojlim_i H^{j}(U_i \cap N_x, \Q)$
(see \cite{WeHA} p.83, Ex 3.5.2, 3.5.8).
Since $U_i \cap N_x$ is a finite group,
their higher cohomology with $\Q$
coefficient vanishes.   This proves 
the claim, hence the lemma.
\end{proof}
We note that $\pi_0(X_{\GL_{d_j},\bK_j,\bullet})$ is
canonically isomorphic to
$\GL_{d_i}(\A^\infty)/\bK_j$
for $j=1,\ldots,r$, and Lemma~\ref{7_contract} implies
that $I_{Y_{\bK,\bullet}^{(\alpha),\cD}}$ is
canonically isomorphic to
$\Map(P_{\cD}(\A^\infty)/
(P_{\cD}(\A^\infty) \cap \bK),\Q)$.
Since $N_{\cD}(F)$ is dense in 
$N_{\cD}(\A^\infty)$, the group
$H^0(N_{\cD}(F), I_{Y_{\bK,\bullet}^{(\alpha),\cD}})$
is canonically isomorphic to the group
$\Map(M_{\cD}(\A^\infty)/\prod_{j=1}^r \bK_j,\Q)$.
Hence the upper horizontal arrow
of the diagram (\ref{7_CD}) is an isomorphism.

Next we consider the vertical arrows.
Each connected component of
$\wt{X}_{\GL_{d_j},g_j^{-1} \bK_j g_j,\bullet}$
is contractible since it is isomorphic to 
the Bruhat-Tits building for $\GL_{d_j}$.
Recall that the simplicial complex
$X_{\GL_{d_j},g_j^{-1} \bK_j g_j,\bullet}$
is the quotient of 
$\wt{X}_{\GL_{d_j},g_j^{-1} \bK_j g_j,\bullet}$
by the action of $\GL_{d_j}(F)$.
For any simplex $\sigma$ in 
$\wt{X}_{\GL_{d_j},g_j^{-1} \bK_j g_j,\bullet}$,
the stabilizer group $\Gamma_\sigma \subset \GL_{d_j}(F)$
of $\sigma$ is finite, as remarked in Section~\ref{sec:finite stab}.
Hence the left two vertical arrows in the diagram (\ref{7_CD})
are isomorphisms.
Similarly, bijectivity of the two right
vertical arrows in the diagram (\ref{7_CD})
follows from Lemma~\ref{7_contract}. Thus we have
a proof of Proposition~\ref{7_prop3}.
\end{proof}

\begin{proof}[Proof of Proposition~\ref{7_prop2}]
Let us consider the lower horizontal arrow
in the diagram (\ref{7_CD}). By Proposition~\ref{7_prop3} 
it is an isomorphism. We note that the compact
open subgroups $\bK \subset \GL_d(\wh{A})$
with property (\ref{7_property}) form a 
cofinal subsystem of the
inductive system of all open compact subgroups
of $\GL_d(\A^\infty)$. Therefore,
passing to the inductive limits with respect to
$\alpha$ and $\bK$ with property (\ref{7_property}), we have
$\varinjlim_\bK \varinjlim_\alpha
H^*(Y_{\bK,\bullet}^{(\alpha),\cD},\Q)
\cong \bigotimes_{j=1}^r \varinjlim_{\bK_j} 
H^*(X_{\GL_{d_j},\bK_j,\bullet},\Q)$ as desired.
\end{proof}

\subsection{}
Let $\bK, \bK' \subset \GL_d(\A^\infty)$ be two
compact open subgroups with $\bK' \subset \bK$.
The pull-back morphism from the cochain complex
of $X_{\bK,\bullet}$ to that of $X_{\bK',\bullet}$
preserves the cochains with finite supports.
Thus we have pull-back homomorphisms 
$H^*_c(X_{\bK,\bullet},\Q) \to H^*_c(X_{\bK',\bullet},\Q)$
which is compatible with
the usual pull-back homomorphism 
$H^*(X_{\bK,\bullet},\Q) \to H^*(X_{\bK',\bullet},\Q)$.
For an abelian group $M$, we let
$H^*(X_{\lim,\bullet},M)=H^*(X_{\GL_d,\lim,\bullet},M)$ and
$H_c^*(X_{\lim,\bullet},M)=H_c^*(X_{\GL_d,\lim,\bullet},M)$
denote the
inductive limits
$\varinjlim_{\bK}H^*(X_{\bK,\bullet},M)$ and
$\varinjlim_{\bK}H_c^*(X_{\bK,\bullet},M)$,
respectively.
If $M$ is a $\Q$-vector space, then
for each compact open subgroup 
$\bK \subset \GL_d(\A^\infty)$,
the homomorphism
$H^*(X_{\bK,\bullet},M) \to H^*(X_{\lim,\bullet},M)$
is injective and its image is equal to the
$\bK$-invariant part $H^*(X_{\lim,\bullet},M)^{\bK}$
of $H^*(X_{\lim,\bullet},M)$.
Similar statement holds for $H_c^*$.
It follows from Proposition~\ref{prop:admissible} that
the inductive limits $H^{d-1}(X_{\lim,\bullet},\Q)$
and $H_c^{d-1}(X_{\lim,\bullet},\Q)$ are admissible
$\GL_d(\A^\infty)$-modules, and are isomorphic to
the contragradient of $H_{d-1}(X_{\lim,\bullet},\Q)$
and $H_{d-1}^\mathrm{BM}(X_{\lim,\bullet},\Q)$,
respectively.

\subsection{Proof of Theorem~\ref{7_PROP1}}
Since $\St_{d,\C}$ is self-contragradient,
it follows from the compatibility of
the normalized parabolic inductions with taking contragradient
that it suffice to prove that any irreducible subquotient of
$H^{d-1}_c(X_{\lim,\bullet},\C)$
satisfies the properties in the statement of Theorem~\ref{7_PROP1}.
Let $\pi$ be an irreducible subquotient of
$H^{d-1}_c(X_{\lim,\bullet},\C)$.
Then Proposition~\ref{7_prop2} combined with 
the spectral sequence (\ref{7_limspecseq}) shows that
there exists a subset $\cD \subset \{1,\ldots,d-1\}$
such that $\pi^\infty$ is isomorphic to a subquotient of
$\Ind_{P_{\cD}(\A^\infty)}^{GL_d(\A^\infty)}
\bigotimes_{j=1}^r \varinjlim_{\bK_j}
H^{d_j -1}(X_{\GL_{d_j},\bK_j,\bullet},\C)$.
Here $r=\sharp \cD +1$, and $d_1,\ldots,d_{r} \ge 1$ 
are the integers satisfying $\cD=\{d_1,d_1+d_2,\ldots,
d_1+\cdots+d_{r-1} \}$ and $d_1 +\cdots + d_{r} =d$.
By Proposition~\ref{prop:66_3}, $\pi^\infty$ is isomorphic 
to a subquotient of the non-$\infty$-component
of the induced representation from $P_{\cD}(\A)$ to
$\GL_d(\A)$ of an irreducible cuspidal automorphic representation
$\pi_1 \otimes \cdots \otimes \pi_r$ of $M_{\cD}(\A)$ whose component 
at $\infty$ is isomorphic to
the tensor product of the Steinberg representations. 

It remains to prove the claim of the multiplicity.
The Ramanujan-Petersson conjecture proved by Lafforgue
\cite[p.6, Th\'{e}or\`{e}me]{Laf}
shows that each place $v$ of $F$, the representation $\pi_{i,v}$
is tempered. Hence for almost all places $v$ of $F$,
the representation $\pi_v$ of $\GL_d(F_v)$ is unramified
and its associated Satake parameters $\alpha_{v,1},\cdots ,\alpha_{v,d}$
have the following property: for each $i$ with $1 \le i \le r$,
exactly $d_i$ parameters of $\alpha_{v,1},\cdots ,\alpha_{v,d}$
have the complex absolute value $q_v^{a_i/2}$ where
$q_v$ denotes the cardinality of the residue field at $v$
and $a_i = \sum_{i<j\le r} d_j - \sum_{1 \le j <i} d_j$.
This shows that the subset $\cD$ is uniquely determined by $\pi$.
It follows from the multiplicity one theorem and the strong multiplicity
one theorem that the
cuspidal automorphic representation 
$\pi_1 \otimes \cdots \otimes \pi_r$ of 
$M_{\cD}(\A)$
is also uniquely determined by $\pi$. 

Hence it suffices to show the following lemma.
\begin{lem}
The representation $\Ind_{P_{\cD}(F_v)}^{\GL_d(F_v)} 
\pi_{1,v} \otimes \cdots \otimes \pi_{r,v}$ of $\GL_d(F)$ is of
multiplicity free for every place $v$ of $F$.
\end{lem}

\begin{proof}
For $1 \le i \le r$, let $\Delta_i$ denote the multiset of
segments corresponding to the representation 
$\pi_{i,v}\otimes |\det(\ )|_v^{a_i/2}$
in the sense of \cite{Zelevinsky}. We denote by $\Delta_i^t$
the Zelevinski dual of $\Delta_i$.
Let $i_1, i_2$ be integers with $1 \le i_1 < i_2 \le r$
and suppose that there exist a segment in $\Delta_{i_1}^t$ 
and a segment in $\Delta_{i_2}^t$ which are linked.
Since $\pi_{i_1,v}$ and $\pi_{i_2,v}$ are tempered,
it follows that $i_2 = i_1 +1$ and that there exists a character 
$\chi$ of $F_v^\times$ such that both $\pi_{i_1,v} \otimes \chi$
and $\pi_{i_2,v} \otimes \chi$ are the Steinberg representations.
In this case the multiset $\Delta_{i_j}^t$ consists of
a single segment for $j=1,2$ and 
the unique segment in $\Delta_{i_1}^t$ and the
unique segment in $\Delta_{i_2}^t$ are juxtaposed.
Thus the claim is obtained by applying the formula in 
\cite[9.13, Proposition, p.201]{Zelevinsky}.
\end{proof}
This finishes the proof of Theorem~\ref{7_PROP1}. 
\qed

\chapter{Universal Modular Symbols}
We recall here the definition of universal modular symbols following 
Ash and Rudolph \cite[\S 2]{AR}.    
Let $F$ be a field and $d \ge 1$.
Let $q_1, \dots, q_d$ be an ordered basis 
(a basis with the order fixed) 
of $F^{\oplus d}$.
A universal modular symbol 
$[q_1, \dots, q_d]$ associated with the ordered 
basis is then 
a $(d-2)$-nd 
homology class 
of the Tits building $T_{F^{\oplus d}}$ 
of the algebraic group $\mathrm{SL}_d$ over $F$.

The treatment here may look slightly 
different from the paper \cite{AR}
because we put our emphasis on posets.
By definition, the Tits building $T_{F^{\oplus d}}$ 
is (the classifying space of) the poset
of $F$-subspaces of $F^{\oplus d}$.
Its barycentric subdivision is then the classifying
space of the poset of flags in $F^{\oplus d}$.
On the other hand, the boundary of the barycentric subdivision 
of the standard $(d-1)$-simplex $\Delta_{d-1}$
is the classifying space of the poset of subsets of 
$\{1,\dots, d\}$, excluding $\{1,\dots, d\}$ and $\emptyset$.
A choice of an ordered basis will give a morphism from 
this poset to the poset of flags.    This in turn induces 
a map of homology and the image of the fundamental
class of the homology of the boundary of the standard simplex
is defined as the universal symbol corresponding to the 
ordered basis.

We mention some results of \cite{AR} in Section~\ref{sec:AR results}.
Of course, our main aim is to consider the analogue of their Proposition 3.2;
this will be covered in Chapters~\ref{ch:pf for ums}.   
We use that the universal modular symbols generate  $H_{d-2}(T_{F^{\oplus d}})$.
Their main result (see Theorem~\ref{thm:AR main}) is applicable when our base ring $A$ is a Euclidean domain.

\section{the Tits building for the special linear groups}
\label{sec:Tits building}
Let us recall the definition of the Tits building.  

\subsection{}
Let $F$ be a field.   (We will take $F$ to be our ground global field of 
positive characteristic for our application.)
Let $d \ge 1$ be a positive integer, and $W$ be a $d$-dimensional
vector space over $F$.
We denote by $Q(W)$  the poset of nonzero proper vector subspaces $0 \subsetneqq W_0 \subsetneqq W$.
Let $\cP_\tot(Q(W))$ denote the poset of totally ordered finite subsets of $Q(W)$.
By definition, the Tits building $T_W$ is the simplicial complex 
$(Q(W), \cP_\tot(Q(W)))$ associated with
the poset $Q(W)$.
In the terminology of \cite[\S 1]{Quillen},
$T_W$ is the classifying space of $Q(W)$.

\subsection{}
We use the poset of flags in $W$, which gives rise to 
the barycentric subdivision of $T_W$.  
A flag in $W$ is a sequence,
for some $1 \le i \le d$,
\[
0=W_0 \subsetneq W_1 \subsetneq W_2 \subsetneq \dots \subsetneq
W_i \subsetneq W
\]
of $F$-subvector spaces.   
Let $F(W)$ denote the set of flags in $W$.   
For two flags $F_1, F_2$, we set $F_1 \le F_2$ if 
$F_2$ is a refinement of $F_1$.   
With this order, $F(W)$ is a poset.
We have a canonical identification 
$F(W)=\cP_\tot(Q(W))$ of two posets.

We denote by $T'_W$ the simplicial complex $(F(W), \cP_\tot(F(W)))$
where $\cP_\tot(F(W))$ is the set of 
nonempty totally ordered finite subsets of $F(W)$.
Then $T'_W$ is the barycentric subdivision of $T_W$.
(See proof of Lemma 1.9 in \cite{Gra}.)

\subsection{}
Let us take an isomorphism $W\cong F^{\oplus d}$ and regard 
$W$ as the set of column vectors. 
Then $\GL_d(F)$ acts on $W$ by multiplication from the left,
on the poset $Q(W)$, on $T_W$ and on $T'_W$.

\subsection{}
As mentioned in \cite[p.165, 2.2.3]{FKS}, $F(W)$ is isomorphic as a poset
to the set of parabolic subgroups (ordered by inclusion) 
of any of $GL_d$, $PGL_d$, and $SL_d$ over $F$.
This gives the description of $T_W$ as the Tits building 
of the semisimple algebraic group $SL_d$ over $F$.   (See \cite[p.242, Section 2]{AR}.)

\section{the boundary of the first barycentric subdivision of 
a standard simplex}
\label{sec:barycentric subdivision}
We now give the simplicial complex which describes 
the boundary of the first barycentric subdivision 
of the standard $(d-1)$-simplex.
This description is essentially the same as 
that in \cite[p.243]{AR}.

\subsection{}
Let $B$  be a nonempty finite set.
Let $\cP'(B)$ denote the set of subsets
$J \subset B$ satisfying $J \neq \emptyset, B$.
We regard $\cP'(B)$ as a partially ordered set
with respect to the inclusions.

Let $\cP_\tot(\cP'(B))$ denote the set of
non-empty totally ordered subsets of $\cP'(B)$.
Then the pair $(\cP'(B), \cP_\tot(\cP'(B)))$ forms
a finite simplicial complex.  It is the classifying space
of the poset $\cP'(B)$.

Let $d=|B|$ denote the cardinality of $B$.
Then this simplicial complex is isomorphic to the
boundary of the first 
barycentric subdivision of the standard $(d-1)$-simplex.
The set of vertices of this standard $(d-1)$-simplex $\Delta$
is the set $B$.   

\subsection{}
The following is the picture of the boundary of the first barycentric subdivision of 
the standard $2$-simplex:

\noindent
  
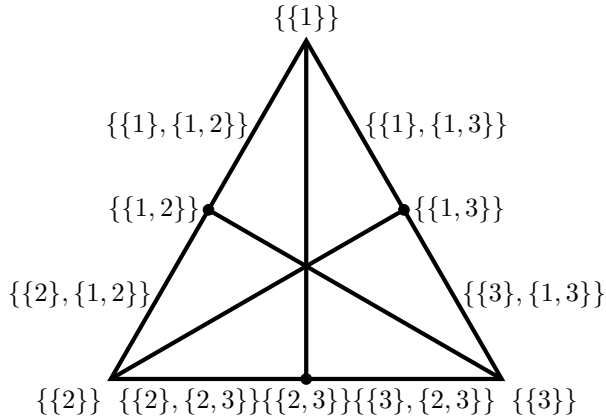
\begin{tikzpicture}
    \path
      (90:3cm) coordinate (a) node[above] {$\{\{1\}\}$}
      (210:3cm) coordinate (b) node[below left] {$\{\{2\}\}$}
      (-30:3cm) coordinate (c) node[below right] {$\{\{3\}\}$};

    \draw [ultra thick] (a) -- (b) -- (c) --cycle;
    
    \draw [ultra thick, x=(a),y=(b),z=(c)] (1,0,0) -- (0,.5,.5);
    \draw [ultra thick, x=(a),y=(b),z=(c)] (0,1,0) -- (.5,0,.5);
    \draw [ultra thick, x=(a),y=(b),z=(c)] (0,0,1) -- (.5,.5,0);
    
\filldraw[black, x=(a),y=(b),z=(c)] (0,.5,.5) circle (2pt) node[anchor=north] {$\{\{2,3\}\}$};
\filldraw[black, x=(a),y=(b),z=(c)] (.5,0,.5) circle (2pt) node[anchor=west] {$\{\{1,3\}\}$};
    \filldraw[black, x=(a),y=(b),z=(c)] (.5,.5,0) circle (2pt) node[anchor=east] {$\{\{1,2\}\}$};
    
    \filldraw[black, x=(a),y=(b),z=(c)] (.75,0,.25) node[anchor=west] {$\{\{1\},\{1,3\}\}$};
    \filldraw[black, x=(a),y=(b),z=(c)] (.25,0,.75) node[anchor=west] {$\{\{3\},\{1,3\}\}$};
    
    \filldraw[black, x=(a),y=(b),z=(c)] (0,.8,.2) node[anchor=north] {$\{\{2\},\{2,3\}\}$};
    \filldraw[black, x=(a),y=(b),z=(c)] (0,.2,.8) node[anchor=north] {$\{\{3\},\{2,3\}\}$};
    
        \filldraw[black, x=(a),y=(b),z=(c)] (.75,.25,0) node[anchor=east] {$\{\{1\},\{1,2\}\}$};
        \filldraw[black, x=(a),y=(b),z=(c)] (.25,.75,0) node[anchor=east] {$\{\{2\},\{1,2\}\}$};
  \end{tikzpicture}

\subsection{an orientation}
If $B$ is totally ordered, then the order gives an orientation of  $\Delta$.
This orientation in turn gives an orientation of the subdivision, and of the boundary of the barycentric 
subdivision.   Let us describe this explicitly. 

The number of $(d-1)$-dimensional simplices in 
$(\cP'(B), \cP_\tot(\cP'(B)))$ is equal to $d!$. They are the simplices $\sigma_g$ 
whose vertices are 
\[
\{
g(1)\},
\{g(1), g(2)\},
\dots,
\{g(1), g(2), \dots, g(d-1)\}
\]
for each $g \in S_d$ in the $d$-th symmetric group.

The orientation of $\sigma_g$ that we use below (which is the one determined from the orientation of $\Delta$) is as follows.   If $g$ is an even permutation, then the orientation is to be given by the increasing order
$\{
g(1)\},
\{g(1), g(2)\},
\dots,
\{g(1), g(2), \dots, g(d-1)
\}$.
If $g$ is odd, then the orientation is the opposite of the order above.

This choice of orientation (or the total order) gives an element (the fundamental class)
in 
\[H_{d-2}((\cP'(B), \cP_\tot(\cP'(B))), \Z).\]

\section{the universal modular symbols}
\label{sec:univ ms}
We define the universal modular symbols of Ash-Rudolph 
(\cite[Definition 2.1]{AR}).
\subsection{}
Let $d \ge 1$ be an integer.
Let $B=\{1,\dots, d\}$ be a totally ordered finite set of cardinality $d$.
Let $F$ be a field and $W$ be a $d$-dimensional $F$-vector space.

Let $q_1, \dots, q_d \in W$ be a basis of $W$.

We define a morphism of posets
\[
\phi=\phi_{q_1, \dots, q_d}:
\cP'(B) \to Q(W)
\]
as follows.
Let $1 \le i \le d$ and $J=\{j_1,\dots, j_i\} \subset B$ belonging 
to $\cP'(B)$.   Then we set
\[
\phi(J)=\langle q_{j_1}, \dots, q_{j_i} \rangle \subsetneqq W
\]
to be the vector subspace spanned by 
$\{q_{j_1}, \dots, q_{j_i}\}$.

The morphism of posets induces a morphism of classifying spaces
which we also denote by $\phi=\phi_{q_1,\dots, q_d}$:
\[
\phi=\phi_{q_1, \dots, q_d}:
((\cP'(B), \cP_\tot(\cP'(B))))
\to 
T_W=
(Q(W), \cP_\tot(Q(W))).
\]

\subsection{}
We write 
$[q_1,\dots, q_d]$
for the image of the fundamental class in 
\[
H_{d-2}((F(W), \cP_\tot(F(W))), \Z)
\]
by the pushforward by $\phi_{q_1, \dots, q_d}$.

For convenience, we set 
$[q_1, \dots, q_d]=0$ if 
$q_1, \dots, q_d$ do not form a basis.

These elements are called  
universal modular symbols.

When $W=F^{\oplus d}$, for $Q \in \GL_d(F^{\oplus d})$,
we write $[Q]$ for the universal modular symbol
$[q_1,\dots, q_d]$ where $q_i$ is the $i$-th row for 
$1 \le i \le d$.

\section{Some results of Ash and Rudolph}
\label{sec:AR results}
We record here some results of Ash and Rudolph
on universal modular symbols 
which may or may not hold in our case.
The first statements in the case of $\GL_2$
were given by Manin \cite{Manin1}.

\subsection{}
We have a description of a $\Z$-basis of the homology. 
\begin{prop}
[Prop 2.3, p.244, \cite{AR}]
\label{prop:MS Z-basis}
Let $U$ denote the 
subgroup of $\GL(W)$ 
consisting of unipotent upper triangular matrices.
Then the symbols $[Q]$
as $Q$ runs through $U$,
make up a $\Z$-basis of $H_{d-2}(T_W, \Z)$.
\end{prop}
\begin{proof}
See loc.\ cit.
\end{proof}

\subsection{}
The following proposition gives some 
basic properties of the universal modular symbols.
\begin{prop}[{\cite[Prop 2.2, p.243]{AR}}]
The universal modular symbols enjoy the following properties:

\noindent
1. It is anti-symmetric.
\\
2. $[aq_1, q_2 \dots, q_d]=[q_1, \dots, q_d]$ for any nonzero $a \in F$.
\\
3. If $q_1, \dots, q_{d+1}$ are all nonzero, then
\[
\sum_{i=1}^{d+1}
(-1)^{i+1}
[q_1, \dots, \widehat{q}_i, \dots, q_{d+1}]=0.
\]
\\
4. If $A \in \GL_d(F)$, then 
$[AQ]=A \cdot [Q]$,
where the dot denotes the action of $\GL_F(W)$ 
on the homology of $T_W$.
\end{prop}
\begin{proof}
See {\cite[Prop 2.2, p.243]{AR}}.
\end{proof}

\subsection{}
We mention their main theorem. 
Here $T_d$ is the Tits building 
for subspaces in $L^n$
with $L$ the field of fractions of 
a Euclidean domain $\Lambda$.
\begin{thm}[{\cite[Thm 4.1, p.247]{AR}}]
\label{thm:AR main}
As $Q$ runs over $\mathrm{SL}(d, \Lambda)$,
the universal modular 
symbols $[Q]$ generate $H_{d-2}(T_d; \Z)$.
\end{thm}
\begin{proof}
See \cite[Thm 4.1]{AR}
\end{proof}

Ash and Rudolph use this theorem for $\Lambda=\Z$.  
For our setup, the ring $A$ is the analogue of $\Z$
but this is in general not a Euclidean domain.
We do not know to what extent 
this type of statement holds true.

\chapter{On finite $p$-subgroups of arithmetic subgroups}
\section{Introduction}
The result of this chapter will be used in Section~\ref{sec:d argument}.
\subsection{}
We are interested in the group homology of 
stabilizer groups of arithmetic subgroup action
on the Bruhat-Tits building $\cBT_\bullet$.
These homology groups appear as the $E_1$-terms
of the spectral sequence in Section~\ref{sec:2nd ss}.
In order to obtain the estimate of the exponent, 
we do not study the differentials of the spectral
sequence directly, but we give an 
estimate of the exponent of 
the $E_1$-terms, using the 
fact that group homology of a finite group 
is killed by the order of the group.   
In this chapter, we treat the $p$-part.

\subsection{}
Let $\Gamma$ be an arithmetic 
subgroup acting on $\cBT_\bullet$.   
Suppose that $\Gamma$ is a pro-$p$ group,
where $p$ is the characteristic of $F$.
Let 
$\sigma$ be a simplex of $\cBT_\bullet$
and let $\Gamma_\sigma$ be the stabilizer
of $\sigma$.  Then $\Gamma_\sigma$
is a finite $p$-subgroup of $\GL_d(F)$.
We will bound the exponent of the 
homology groups of $\Gamma_\sigma$.
We do not use the fact that it is a stabilizer group,
but merely the fact that it is a $p$-subgroup of $\GL_d(F)$.

In Lemma~\ref{lem:existence of flag} 
below, it is shown that there exists a flag 
in $F^{\oplus d}$ which is stabilized by $\Gamma_\sigma$.
This in turn gives a filtration of the stabilizer group by normal 
subgroups such that successive quotients are either trivial
or an elementary abelian $p$-group (Corollary~\ref{cor:length d-1}).
Since the group homology of an elementary abelian 
$p$-group is killed by $p$, 
we obtain a certain bound of the order of the homology groups 
by the length of the filtration.

\section{Lemmas}

\begin{lem} \label{lem:p-group}
Let $H$ be a finite $p$-group.
Suppose that there exists a decreasing filtration
$$
H=H^0 \supset H^1 \supset \cdots \supset H^{\ell-1} \supset
H^{\ell} = \{1\}
$$
of $H$ by normal subgroups of $H$ such that
$H^{i-1}/H^i$ is an elementary abelian $p$-group
for $i=1,\ldots, \ell$.
Let $R$ be a principal ideal domain and let
$\chi: H \to R^\times$ be a character of $H$.
Then for any integer $s \ge 1$, the $s$-th homology group
$H_s(H,\chi)$ is an abelian group killed by
$p^{1+s(\ell-1)}$.
\end{lem}

\begin{proof}
First we prove the claim for $\ell=1$.
In this case $H$ is an elementary abelian $p$-group.
Since the claim follows from the 
K\"unneth theorem if $\chi$ is trivial, we may assume
that $\chi$ is non-trivial.
Let $H'$ denote the kernel of $\chi$.
Since $H$ is an elementary abelian $p$-group,
there exists a cyclic subgroup $H'' \subset H$
of order $p$ such that the map $H' \times H'' \to H$ 
given by the multiplication in $H$ is
an isomorphism of groups.
Via this isomorphism, the character $\chi$ is regarded as
the external tensor product over $R$ 
of the trivial character of $H'$ and the restriction 
$\chi|_{H''}$ of $\chi$ to $H''$.
Hence it follows from the splitting of the short exact
sequence in the K\"unneth theorem that
$H_s(H,\chi)$ is killed by $p$.

We prove the claim for $\ell>1$ by induction on $\ell$.
Let us consider the Hochschild-Serre spectral sequence
(cf. \cite[p.171, VII, (6.3)]{Brown})
$$
E^2_{s,t} = H_s(H/H^{\ell-1},H_t(H^{\ell-1},\chi))
\Longrightarrow  H_{s+t}(H,\chi).
$$
The claim for the elementary abelian $H^{\ell-1}$ shows 
that $E^2_{s,t}$ is killed by $p$ for $t \ge 1$.
Since the claim is known for
$H/H^{\ell-1}$ by the inductive hypothesis, 
$E^2_{s,0}$ is killed by $p^{1+s(\ell-2)}$.
Hence $E^r_{q,0}$ is killed by $p^{1+q(\ell-2)}$
and $E^r_{s-i,i}$ are killed by $p$ for $i=1,\ldots,s$ and for any $r \ge 2$.
Thus the spectral sequence above shows that
$H_s(H,\chi)$ is killed by 
$p^{1+s(\ell-2)} \cdot \prod_{i=1}^s p = p^{1+s(\ell-1)}$.
\end{proof}

\begin{lem} \label{lem:existence of flag}
Let $V$ be a non-zero, finite dimensional $F$-vector space, and
$H$ a finite $p$-group which acts $F$-linearly on $V$ from the right.
Then there exists a flag
$0 = V_0 \subsetneqq V_1 \subsetneqq \cdots \subsetneqq V_\ell =  V$
in $V$ by $F$-linear subspaces
such that $V_i \cdot H = V_i$ and $H$ acts trivially on 
$V_i/V_{i-1}$ for $i=1,\ldots,\ell$
\end{lem}

\begin{proof}
First we prove that the $H$-invariant part $V^H$ is non-zero
by induction on the order $p^m$ of $H$. If $m \le 1$, then 
$H$ is generated by an element $h \in H$. 
Let $\rho(h) \in \GL_F(V)$ denote the action of $h$ on $V$.
Since $h^p =1$, we have $(\rho(h)-\id_V)^p=0$.
This implies that $1$ is an eigenvalue of $\rho(h)$. Hence $V$ has 
a non-zero $H$-invariant vector. Suppose that $m \ge 2$.
Since any nontrivial finite $p$-group has a non-trivial center, 
there exists an element $h$ of order $p$ in the center of $H$.
Let $W \subset V$ denote the subspace of $h$-invariant vectors.
Then $W \neq 0$ and $H/\langle h \rangle$ acts on $W$.
Hence by inductive hypothesis, $V^H = W^{H/\langle h \rangle}$
is non-zero.

Next we prove the claim by induction on $d = \dim_F V$.
The claim is clear when $d=1$ since $V^H \neq \{0\}$ implies
that $H$ trivially acts on $V$. Suppose that $d \ge 2$.
Set $W = V/V^H$. If $W=\{0\}$, then
$H$ acts trivially on $V$ and the claim is clear. Suppose otherwise.
Then the action of $H$ on $V$ induces a right action
of $H$ on $W$. By inductive hypothesis, there exists a flag
$0 = W_0 \subsetneqq W_1 \subsetneqq \cdots \subsetneqq W_{\ell'} =  W$
in $W$ by $F$-linear subspaces
such that $W_i \cdot H = W_i$ and $H$ acts trivially on 
$W_i/W_{i-1}$ for $i=1,\ldots,\ell'$.
Then the preimage $V_i$ of $W_{i-1}$ under the quotient map
$V \surj W$ for $i=1,\ldots,\ell=\ell'+1$
gives a desired flag in $V$.
\end{proof}

\begin{cor} \label{cor:length d-1}
Let $H \subset \GL_d(F)$ be a finite subgroup that is a finite $p$-group.
Then there exists an integer $\ell \le d$ and a decreasing filtration
$$
H=H^0 \supset H^1 \supset \cdots \supset H^{\ell-2} \supset
H^{\ell-1} = \{1\}
$$
of $H$ by normal subgroups of $H$ such that
$H^{i-1}/H^i$ is either trivial or an elementary abelian $p$-group
for $i=1,\ldots, \ell-1$.
\end{cor}

\begin{proof}
It follows from Lemma \ref{lem:existence of flag} that there 
exists a flag
$0 = V_0 \subsetneqq V_1 \subsetneqq \cdots \subsetneqq V_\ell =  F^d$
in $F^d$ by $F$-linear subspaces
such that $V_i \cdot H = V_i$ and $H$ acts trivially on 
$V_i/V_{i-1}$ for $i=1,\ldots,\ell$.
Since $F^d$ is $d$-dimensional, we have $\ell \le d$.
For $i=0,\ldots,\ell-1$, set $H^i = \Ker(H \to \GL_F(V_{i+1})) \subset H$.
Then $H^i$ is a normal subgroup of $H$ and we have
$H=H^0 \supset H^1 \supset \cdots \supset H^{\ell-1}= \{1\}$.
Let $0 \le i \le \ell-2$.
Then for $h \in H^i$, the $F$-linear map $V_{i+2} \to V_{i+1}$ that sends
$v \in V$ to $vh -v$ induces an $F$-linear map 
$\varphi_h : V_{i+2}/V_{i+1} \to V_{i+1}$. Let
$\varphi : H^i \to \Hom_F(V_{i+2}/V_{i+1},V_{i+1})$
denote the map that sends $h \in H^i$ to $\varphi_h$.
Then one can check easily that the map $\varphi$
is a homomorphism of groups and that $H^{i+1}$ is equal to
the kernel of $\varphi$. This shows that $H^i/H^{i+1}$ is an 
elementary abelian $p$-group since it is isomorphic to
a subgroup of $\Hom_F(V_{i+2}/V_{i+1},V_{i+1})$.
Hence the groups $H^0, \ldots, H^{\ell-1}$ have
the desired property.
\end{proof}

\begin{cor}
\label{cor:finite p-group homology}
Let $H \subset \GL_d(F)$ 
be a finite subgroup that is a finite $p$-group.  
Let $\chi:H \to \Z^\times$ be a character.
Then the homology groups
$H_s(H, \chi)$ is killed by 
$p^{1+s(d-2)}$ 
for $s \geqslant 1$.
\end{cor}
\begin{proof}
This follows from Lemma~\ref{lem:p-group}
and Corollary~\ref{cor:length d-1}.
\end{proof}

\chapter{Some spectral sequences}
In this section, we mention two spectral sequences that
we use in Chapter~\ref{ch:pf for ums}.   They are ordinary ones 
that are found in textbooks such as \cite{Brown}
or \cite{WeHA}.   Because there are minor differences,
we record them here.   The reader may safely skip this 
chapter.

\section{Cellular structure on some products}
\label{seq:CW_prod}
Let $Y_\bullet$ be a simplicial complex and $Y'_\bullet$ a simplicial set.
The geometric realizations 
$|Y_\bullet|$ and $|Y'_\bullet|$
have canonical structures of CW complexes, induced by 
the triangulations given by the simplices in
$Y_\bullet$ and $Y'_\bullet$, respectively.
The complex computing the cellular (co)homology groups
of the CW complex $|Y_\bullet|$
is identical to
the complex in Section \ref{sec:def homology}
computing the (co)homology groups of $Y_\bullet$.

Suppose that $|Y_\bullet|$ is locally compact.
Then, by \cite[Theorem A.6, p.525]{Hatcher}, 
the structures of CW complexes on 
$|Y_\bullet|$ and $|Y'_\bullet|$ naturally
induce a structure of CW complex on
the direct product $|Y_\bullet| \times |Y'_\bullet|$ 
of topological spaces, where the cells of dimension $i$
are indexed by the disjoint union
$$
\coprod_{j=0}^i Y_j \times Y'^{\nd}_{i-j}
$$
for any integer $i \ge 0$, where
$Y'^{\nd}_{i-j} \subset Y'_{i-j}$ denotes the
set of non-degenerate $(i-j)$-simplices in $Y'_\bullet$,
and the closure of
any cell whose index is in $Y_j \times Y'^\nd_{i-j}$
is isomorphic to the direct product of
a $j$-simplex and an $(i-j)$-simplex equipped with
a natural structure of CW complex.
We call this structure of CW complex on 
$|Y_\bullet| \times |Y'_\bullet|$ the product cellular
structure on $|Y_\bullet| \times |Y'_\bullet|$.

For a CW complex $Z$, we denote by 
$C_\bullet^\cell(Z,\Z)$
the complex computing the cellular homology groups
of $Z$. Since the two projections
$|Y_\bullet| \times |Y'_\bullet| \to |Y_\bullet|$
and
$|Y_\bullet| \times |Y'_\bullet| \to |Y'_\bullet|$
are cellular maps, they induce maps
$C_\bullet^\cell(|Y_\bullet| \times |Y'_\bullet|,\Z) 
\to C_\bullet(|Y_\bullet|,\Z)$
and
$C_\bullet^\cell(|Y_\bullet| \times |Y'_\bullet|,\Z) 
\to C_\bullet(|Y'_\bullet|,\Z)$
of complexes.
We will later use the terminology in this
paragraph when
$Y_\bullet = \cB\cT_\bullet$ and $Y'_\bullet = E\Gamma_\bullet$.
In this case, we note that $Y_\bullet = \cB\cT_\bullet$
is locally compact.

\section{Equivariant homology}
\label{Equivariant homology}		
Let $\Gamma \subset \GL_d(K)$
be an arithmetic subgroup.
We define the simplicial set 
(not a simplicial complex)
$E\Gamma_\bullet$ as follows.
We define $E\Gamma_n=\Gamma^{n+1}$ to be the $(n+1)$-fold direct product of $\Gamma$ for $n\ge 0$.
The set $\Gamma^{n+1}$ is naturally regarded as the set of maps of sets
$\Map(\{0,\dots, n\}, \Gamma)$ and from this one obtains naturally the 
structure of a simplicial set.
We let $|E\Gamma_\bullet|$ denote the geometric realization of $E\Gamma_\bullet$.  
Then $|E\Gamma_\bullet|$ is contractible.
We let $\Gamma$ act diagonally on each $E\Gamma_n$ $(n\ge 0)$.  
The induced action on $|E\Gamma_\bullet|$ is free.

Let $M$ be a topological space on which $\Gamma$ acts.
The diagonal action of $\Gamma$ on $M\times |E\Gamma_\bullet|$ is free.
We let $H_*^\Gamma(M, B)= H_*(\Gamma\bsl(M \times |E\Gamma_\bullet|), B)$ where $B$ is a coefficient ring,
and call it the equivariant homology of $M$ with coefficients in $B$.  
We also use the relative version, 
and define equivariant cohomology in a similar manner. 

\section{the Lyndon-Hochschild-Serre spectral sequence}
Let $Y_\bullet \subset Z_\bullet$ be simplicial compexes
with compatible $\Gamma$-action.    
We have (see \cite[p.172, VII 7]{Brown}) 
the following spectral sequence.
\[
E^2_{p,q}=
H_p(\Gamma, H_q(Z_\bullet, Y_\bullet; \Z)) 
\Rightarrow 
H_{p+q}^\Gamma(Z_\bullet, Y_\bullet; \Z)
\]

\section{Another spectral sequence}
\label{sec:another ss}
We use another spectral sequence which also converges to 
equivariant homology groups.   We make a slight change from 
\cite[VII.7]{Brown} because we use particular
cell structures.   The reader may skip this section.

\subsection{}
Let $\Gamma$ be an arithmetic subgroup.
We let $Y_\bullet$ be a simplicial complex with $\Gamma$-action.
We consider the product 
$|Y_\bullet| \times |E\Gamma_\bullet|$
as in Section \ref{seq:CW_prod}.

Since $\Gamma$ acts freely on the set $E\Gamma_i$ for every $i \ge 0$,
it also acts freely on the set of $i$-dimensional cells 
with respect to the product
cellular structure on $|Y_\bullet| \times |E\Gamma_\bullet|$, where
$\Gamma$ acts diagonally on the product 
$|Y_\bullet|  \times |E\Gamma_\bullet|$.

Hence the product cellular structure on $Y_\bullet \times |E\Gamma_\bullet|$
induces a structure of CW complex on the quotient 
$\Gamma \bsl  |Y_\bullet| \times |E\Gamma_\bullet|$, 
which enable us to consider the complex
$C_\bullet^\cell(\Gamma \bsl  |Y_\bullet| \times |E\Gamma_\bullet|)$
computing the cellular homology groups of 
$\Gamma \bsl  |Y_\bullet| \times |E\Gamma_\bullet|$.

\subsection{}
Consider the quotient maps
\[
|Y_\bullet| \to \Gamma\bsl |Y_\bullet|
\]
and 
\[
|Y_\bullet|\times |E\Gamma_\bullet| 
\to \Gamma\bsl |Y_\bullet|\times |E\Gamma_\bullet|.
\]
These induce morphisms of complexes
\[
C_\bullet^\cell(|Y_\bullet|, \Z)
\to 
C_\bullet^\cell(\Gamma\bsl|Y_\bullet|, \Z)
\]
and
\[
C_\bullet^\cell(|Y_\bullet| \times |E\Gamma|, \Z) 
\to
C_\bullet^\cell(\Gamma \bsl |Y_\bullet| \times |E\Gamma|, \Z).
\]
It is easy to see that the induced maps
\[
C_\bullet^\cell(|Y_\bullet|, \Z)_\Gamma
\cong
C_\bullet^\cell(\Gamma\bsl|Y_\bullet|, \Z)
\]
and
\[
C_\bullet^\cell(|Y_\bullet| \times |E\Gamma|, \Z)_\Gamma 
\cong
C_\bullet^\cell(\Gamma \bsl |Y_\bullet| \times |E\Gamma|, \Z),
\]
where the subscript $\Gamma$ denotes $\Gamma$-coinvariants,
are isomorphisms.

\subsection{}
Since $|E\Gamma_\bullet|$ is contractible,
the complex 
$C_\bullet^\cell(|E\Gamma_\bullet|,\Z)$ of
$\Z[\Gamma]$-modules is quasi-isomorphic
to $\Z$ with trivial $\Gamma$-action, regarded
as a complex concentrated at degree $0$.

Since $\Gamma$ acts freely on the set of $i$-dimensional cells 
in $|E\Gamma_\bullet|$ for every $i \ge 0$, the complex
$C_\bullet^\cell(|E\Gamma_\bullet|,\Z)$ gives
a $\Z[\Gamma]$-free resolution of 
$\Z$ with trivial $\Gamma$-action.
(We note that $C_\bullet^\cell(|E\Gamma_\bullet|,\Z)$ is not
identical to the standard $\Z[\Gamma]$-free resolution of 
$\Z$, since the degenerate simplices in $E\Gamma_\bullet$
are removed in $C_\bullet^\cell(|E\Gamma_\bullet|,\Z)$.)
\subsection{}
By construction, we can identify the complex
$C_\bullet^\cell(|Y_\bullet| \times |E\Gamma_\bullet|,\Z)$
with the simple complex associated with the
double complex
$$
C_{\bullet,\bullet'} = C_\bullet^\cell(|Y_\bullet|, \Z) \otimes_\Z
C_{\bullet'}(|E\Gamma_\bullet|,\Z)
$$
of $\Z[\Gamma]$-modules (with respect to a suitable sign convention).

\subsection{}
Let $C_{\bullet, 0}$ 
denotes the double complex whose entry at the
bidegree $(i,j)$ is equal to $C_{i,0}$ if $j=0$, and is zero otherwise.

The projection map
$|Y_\bullet| \times |E\Gamma_\bullet| \to |Y_\bullet|$
induces the map of complexes of 
$\Z[\Gamma]$-modules:
$$
C_\bullet^\cell(|Y_\bullet| \times |E\Gamma_\bullet|, \Z) 
\to C_\bullet^\cell(|Y_\bullet|, \Z).
$$

This map is induced by the projection map 
\[
C_{\bullet,\bullet'} \surj C_{\bullet, 0}
\]
of double complexes.

\subsection{}
As we remarked in the previous paragraph, the complex
$C_\bullet^\cell(|E\Gamma_\bullet|,\Z)$ gives 
a $\Z[\Gamma]$-free resolution of $\Z$.
Hence for each $j \ge 0$,
the complex $C_j^\cell(|Y_\bullet|, \Z) \otimes_\Z 
C_\bullet^\cell(|E\Gamma_\bullet|,\Z)$
gives a $\Z[\Gamma]$-free resolution
of the $\Z[\Gamma]$-module $C_j^\cell(Y_\bullet, \Z)$.
Hence the double complex $(C_{\bullet,\bullet'})_\Gamma$
induces a spectral sequence
\begin{equation} 
\label{eq:ss1}
E^{1}_{s,t} = H_t(\Gamma,C_s^\cell(|Y_\bullet|,\Z))
\Longrightarrow 
H_{s+t}^\Gamma(Y_\bullet, \Z)
\end{equation}
where 
\[
\begin{array}{rl}
H_{s+t}^\Gamma(Y_\bullet, \Z)
=&
H_{s+t}(C_\bullet^\cell(|Y_\bullet| \times |E\Gamma_\bullet|,\Z)_\Gamma, \Z))
\\
=& H_{s+t}(C_\bullet^\cell(\Gamma \bsl |Y_\bullet| \times |E\Gamma_\bullet|, \Z)).
\end{array}
\]

By definition, we have 
\[
E^2_{s,0} 
= H_s((C_\bullet^\cell(|Y_\bullet|, \Z))_\Gamma) 
\cong H_s(C_\bullet^\cell(\Gamma \bsl |Y_\bullet|, \Z))= 
H_s(\Gamma \bsl Y_\bullet, \Z).
\]

\subsection{}
For each $i \ge 0$ and for each $\sigma \in \Gamma\bsl |Y_i|$,
let us choose a representative $\wt{\sigma} \in Y_i$ of $\sigma$.
Let $\Gamma_{\wt{\sigma}} \subset \Gamma$ denote 
the stabilizer of $\wt{\sigma}$.
As we mentioned in Section \ref{sec:def Gamma}, it follows from 
the conditions (1) and (3) that $\Gamma_{\wt{\sigma}}$ is a finite
subgroup of $\Gamma$. The group $\Gamma_{\wt{\sigma}}$
may act non-trivially on the set $O(\wt{\sigma})$ 
of orientations of $\wt{\sigma}$. Hence we obtain
a character $\chi_{\wt{\sigma}} : \Gamma_{\wt{\sigma}} \to \{\pm 1\}$ where
an element $\gamma \in \Gamma_{\wt{\sigma}}$ is sent to $-1$ under
$\chi$ if and only if $\gamma$ acts non-trivially on $O(\wt{\sigma})$.
Then it follows from the construction that
for each $i \ge 0$, the $\Z[\Gamma]$-module 
$C_i^\cell(|Y_\bullet|,\Z)$ 
is isomorphic to the direct sum
$$
C_i^\cell(|Y_\bullet|, \Z) \cong
\bigoplus_{\sigma \in \Gamma\bsl Y_i}
\Z[\Gamma] \otimes_{\Z[\Gamma_{\wt{\sigma}}]} \chi_{\wt{\sigma}}
$$
where $\chi_{\wt{\sigma}} = \Z$ with the action of $\Gamma_{\wt{\sigma}}$
given by $\chi_{\wt{\sigma}}$. By Shapiro's lemma, this 
induces an isomorphism
$$
E^1_{s,t} \cong \bigoplus_{\sigma \in \Gamma\bsl Y_s}
H_t(\Gamma_{\wt{\sigma}}, \chi_{\wt{\sigma}})
$$
of abelian groups.

\subsection{}
Let $Y_\bullet$ be as above and let $Z_\bullet$
be a subsimplicial complex with $\Gamma$-action.
Then we obtain the following spectral sequence
in a similar manner:
\[
E^{1}_{s,t} = H_t(\Gamma,C_s^\cell(|Y_\bullet|,\Z)/C_s^\cell(|Z_\bullet|,\Z))
\Longrightarrow 
H^\Gamma_{s+t}(Y_\bullet, Z_\bullet; \Z).
\]
The $E^1$ terms are
\[
E_{s,t}^1 \cong \bigoplus_{\sigma \in \Gamma \bsl Y_s \setminus
\Gamma \bsl Z_s}
H_t(\Gamma_{\wt{\sigma}},\chi_{\wt{\sigma}}).
\]
\chapter{Proof for universal modular symbols}
\label{ch:pf for ums}
In this chapter, we prove Propositions~\ref{prop:thm13(2) AR}, \ref{prop:thm13(3) AR}, \ref{prop:thm13(4) AR}, \ref{prop:thm13(5) AR} and Corollary \ref{cor:thm13(1) AR}.
These are the statements of our main theorem (Theorem~\ref{lem:apartment})
with our modular symbols replaced by (the image of) the universal modular symbols.
We show that our modular symbols coincide with the universal modular symbols 
in Chapter~\ref{ch:compare ms}.

The content of this chapter is outlined in Section~\ref{sec:ums outline}.

\section{Modular symbols for automorphic forms}
Modular symbols are elements of 
$H_{d-1}(\Gamma \backslash \cBT_\bullet, \Z)$
for some arithmetic subgroup $\Gamma$.
For a possible application to automorphic forms,
we may look at all the connected components.
Recall the the space $\cA_{\Z,\St}^\bK$ 
of $\bK$-invariant $\Z$-valued automorphic 
forms is isomorphic to
$H_{d-1}^\BM(X_{\bK, \bullet}, \Z)$.
We have a non-adelic description of 
$X_{\bK, \bullet}$.
We set 
$J_\bK=\GL_d(F)\backslash \GL_d(\A^\infty)/\bK$,
take a set $\{g_j\}$ of representatives 
of $J_\bK$,
and set 
$\Gamma_j=\GL_d(F) \cap g_j \bK g_j^{-1}$.
Then 
\[
X_{\bK,\bullet}
 \cong \coprod_{j \in J_\bK}
\Gamma_j\backslash \cBT_\bullet.
\]

For a compact open subgroup $\bK \subset \GL_d(\A_F)$,
we set 
\[
\MS(\bK)=
\bigoplus_j \MS(\Gamma_j) \subset
\bigoplus_j H_{d-1}(\Gamma_j \backslash \cBT_\bullet,\Z)
\cong 
H_{d-1}(X_{\bK,\bullet}, \Z)
=
\cA^\bK_{\St, \Z}.
\]
This is the space of modular symbols for 
$\bK$-invariant automorphic forms 
$\cA^\bK_{\St, \Z}$.

We are interested in the size of the cokernel
of the map
$\MS(\bK) \subset \cA^\bK_{\St, \Z}$.
For this, it suffices to study the cokernel of 
$\MS(\Gamma) \subset H_{d-1}^\BM(\Gamma \backslash \cBT_\bullet, \Z)$
for an arithmetic group $\Gamma$.   (Because each $\Gamma_j$
is an arithmetic subgroup.)

\section{Connection with the Tits building}
\label{sec:Tits}
Let $\Gamma$ be an arithmetic subgroup.
By definition, there exists a compact open 
subgroup $\bK \subset \GL_d(\A)$ such that 
$\Gamma=\bK \bigcap \GL_d(F)$.

Recall the non-adelic description
from Section~\ref{subsec:X_K}.  
We set 
$J_\bK=\GL_d(F)\backslash \GL_d(\A^\infty)/\bK$,
take a set $\{g_j\}$ of representatives 
of $J_\bK$,
and set 
$\Gamma_j=\GL_d(F) \cap g_j \bK g_j^{-1}$.
Then 
\[
X_{\bK,i} \cong \coprod_{j \in J_\bK}
\Gamma_j\backslash \cBT_i.
\]
Let $e \in \Gamma \subset \GL_d(\A)$ 
be the identity element.
Then $\Gamma=\Gamma_e$.  
We will identify 
$\Gamma \backslash \cBT_\bullet$
with a connected component 
$\Gamma_e \backslash \cBT_\bullet 
\subset X_{\bK, \bullet}$.

We also consider $\cBT_\bullet$ as a subset
$\cBT_\bullet \subset \widetilde{X}_{\bK, \bullet}$
by identifying it with $\{e\} \times \cBT_\bullet$.
We have a map
\[
\cBT_\bullet \subset \widetilde{X}_{\bK,\bullet}
\to X_{\bK, \bullet}
\]
which factors through
$\Gamma \backslash \cBT_\bullet$.

Let $\alpha$ be a positive real number.
Let $\cD \subset \{1, \dots, d-1\}$ be a subset.
Let $f=[0\subset V_1 \subset \dots \subset V_{r-1} \subset F^{\oplus d}]
\in \Flag_\cD$ be a flag.
We defined a subsimplex
$\tilde{X}_{\bK, \bullet}^{(\alpha), \cD, f}$
in Section~\ref{sec:Xsigma}.
We define 
\[
\cBT_\bullet^{(\alpha), \cD, f}
=\cBT_\bullet \cap 
\widetilde{X}_{\bK, \bullet}^{(\alpha), \cD, f}
\subset 
\widetilde{X}_{\bK,\bullet}^{(\alpha)}.
\]
We set
\[
\cBT_\bullet^{(\alpha)}
=
\bigcup_{\cD, f} 
\cBT_\bullet^{(\alpha), \cD, f}.
\]

\begin{lem}
\label{lem:contractible single}
The space $\cBT_\bullet^{(\alpha), \cD, f_0}$
is non-empty and contractible.
\end{lem}
\begin{proof}
This is Lemma \ref{7_contract} stated for 
$\cBT_\bullet$.
\end{proof}

Recall we defined $T_{F^{\oplus d}}$ to be 
the Tits building of $\mathrm{SL}_d$
over $F$.
\begin{cor}[see Corollary 4.2 \cite{Gra}]
\label{cor:Grayson htpy eq}
There is a $\Gamma$-equivariant homotopy equivalence
\[
|\cBT_{\bullet}^{(\alpha)}|
\cong 
|T_{F^{\oplus d}}|
\]
\end{cor}
\begin{proof}
The proof is the same as that in loc.\ cit.
using Lemma~\ref{lem:contractible single} above.
\end{proof}

\section{the Solomon-Tits theorem}
\begin{thm}[Solomon-Tits]
\label{thm:ST}
If $d \ge 2$, 
the Tits building
$T_{F^{\oplus d}}$ has the homotopy type of 
bouquet of $d-2$ spheres.
\end{thm}

\begin{proof}
See \cite[Thm 5.1]{Gra}.
\end{proof}
\begin{lem}
\label{lem:Solomon}
We have 
\[
H_{i}(\cBT_\bullet, \cBT_\bullet^{(\alpha)}; \Z)
\cong
\left\{
\begin{array}{ll}
H_{d-2}(T_{F^{\oplus d}}, \Z)   & i=d-1, \\
0                                      & i\neq d-1.
\end{array}
\right.
\]
\end{lem}

\begin{proof}
Since $\cBT_\bullet$ is contractible 
(see \cite[Thm 2.1]{Gra}),
the claim follows 
from Corollary~\ref{cor:Grayson htpy eq}
and the Solomon-Tits theorem
(Theorem \ref{thm:ST}).
\end{proof}

\section{the Lyndon-Hochschild-Serre spectral sequence}
\label{sec:LHS ss}
We have 
the following spectral sequence.
\[
E^2_{p,q}=
H_p(\Gamma, H_q(\cBT_\bullet, \cBT_\bullet^{(\alpha)}; \Z)) 
\Rightarrow 
H_{p+q}^\Gamma(\cBT_\bullet, \cBT_\bullet^{(\alpha)}; \Z)
\]

By Lemma \ref{lem:Solomon}, 
we have $E^2_{p,q}=0$ unless $q=d-1$.
Hence the spectral sequence collapses at $E^2$ and 
we have moreover
\[
E^2_{0,d-1} \cong E_{d-1}.
\]
Composing with the canonical surjection to the coinvariants,
we obtain a surjection
\[
H_{d-1}(\cBT_\bullet, \cBT_\bullet^{(\alpha)}; \Z)
\rightarrow
H_0(\Gamma, H_{d-1}(\cBT_\bullet, \cBT_\bullet^{(\alpha)}; \Z) )
\cong
H_{d-1}^\Gamma(\cBT_\bullet, \cBT_\bullet^{(\alpha)}; \Z).
\]

\section{the second spectral sequence}
\label{sec:2nd ss}
Let us write $X=\cBT_\bullet$ and 
$X'=\cBT_\bullet^{(\alpha)}$ for short.
\subsection{}
We use the following spectral sequence
(\cite[\S VII 7]{Brown}, see also Section~\ref{sec:another ss})
\[
E_{p,q}^1=
\bigoplus_{\sigma \in \Sigma_p}
H_q(\Gamma_\sigma, \chi_\sigma)
\Rightarrow
H_{p+q}^\Gamma(X, X'; \Z)
\]
where $\Sigma_p$
is the set of $p$-simplices of 
$X \setminus X'$.
Note that by definition
\[
E_{i,0}^2=H_{i}(X,X'; \Z)
\]
for $0 \le i \le d$.

Because $E_{p,q}^1=0$ for $p \le -1$,  
we have 
\[
E_{d-1,0}^\infty=\cdots=E_{d-1,0}^d.
\]
Because $E_{p,q}^1=0$ for $q \le -1$,
we have
\[
E_{d-1,0}^k=\Ker\, d^{k-1}
\]
where $d^{k-1}: E_{d-1,0}^k \to E^k_{d-1+k,k-1}$
is the differential. 

Composing with the canonical map
$E_{d-1} \to E_{d-1,0}^2$, we obtain a homomorphism
\[
H^\Gamma_{d-1}(X, X'; \Z) \to
H_{d-1}(\Gamma \backslash X,\Gamma \backslash X'; \Z)
\]
whose image is $E_{d-1,0}^d$.

\subsection{}
\label{sec:d argument}
Let us give some estimate on the size of 
the image $E_{d-1,0}^d$.
Let $e_{p,q}^r$ denote the exponent of $E_{p,q}^r$.
For each $k$, we have an exact sequence
\[
E_{d-1,0}^{k+1} \to E_{d-1,0}^k \xto{d^{k-1}} E_{d-1+k,k-1}^k
\]
hence 
\[
e_{d-1+k,k-1}^k E_{d-1,0}^k \subset \Ker\, d^{k-1}=E_{d-1,0}^{k+1} 
\subset E_{d-1,0}^k
\]
Therefore 
\[
\prod_k e_{d-1+k,k-1}^k E_{d-1,0}^2 \subset E_{d-1}^d \subset E_{d-1,0}^2.
\]
It is clear that $e_{p,q}^r$ divides $e_{p,q}^1$ for all $p,q,r$,
hence
\[
\prod_k e_{d-1+k,k-1}^1 E_{d-1,0}^2 \subset E_{d-1}^d \subset E_{d-1,0}^2.
\]
It follows from Corollary \ref{cor:finite p-group homology} that 
$E_{s,t}^1$ is killed by 
$p^{1+t(d-2)}$ for $t\ge 1$.   That is,
$e_{s,t}^1$ divides $p^{1+t(d-2)}$. 
We therefore have 
\[
p^{e(d)}E_{d-1}^2 \subset E_{d-1}^d
\subset E_{d-1}^2
\]
where $e(d)=(d-2)\left(1+\frac{(d-1)(d-2)}{2}\right)$.

We arrive at the following lemma.
\begin{lem}
\label{lem:exponent prop}
Suppose $\bK \subset \GL_d(\A^\infty)$ be 
a pro-$p$ compact open subgroup.
Let $\Gamma=\GL_d(F) \cap \bK$.
Then 
\[
H_{d-1}^\Gamma(X, X';\Z)
\to
H_{d-1}(\Gamma \backslash X, \Gamma\backslash X';\Z)\]
is injective 
and the cokernel is annhilated by $p^{e(d)}$.
\end{lem}

\section{Corestriction and transfer}
\subsection{}
Let $\Gamma \supset \Gamma'$ be arithmetic subgroups. 
For a $\Z[\Gamma]$-module $M$, we denote by 
$\cores_M : M_{\Gamma'} \to M_{\Gamma}$ the map induced by the identity map
of $M$.

We define the transfer map 
$\tr_M : M_{\Gamma} \to M_{\Gamma'}$ 
(cf.\ \cite[III, 9.\ (B), p.\ 81]{Brown})
by 
\[
\tr(\overline{m}^\Gamma)
=
\sum_{g\in \Gamma' \bsl \Gamma}
\overline{gm}^{\Gamma'}
\]
where $\overline{m}^\Gamma$ 
and $\overline{m}^{\Gamma'}$
denote the class of $m \in M$
in the coinvariants
$M_\Gamma$ and 
$M_{\Gamma'}$
respectively.
By definition, the composite $\cores_M \circ \tr_M$ is equal the map
given by the multiplication by $[\Gamma:\Gamma']$.

\subsection{}
Let $M_{2, \bullet} \to M_{1, \bullet}$ 
be a morphism of 
complexes of $\Z[\Gamma]$-modules.
Since $\cores_M$ and $\tr_M$ are functorial in $M$, we have
a commutative diagram
$$
\begin{CD}
(M_{2,\bullet})_{\Gamma}
@>{\tr_{M_{2,\bullet}}}>>
(M_{2,\bullet})_{\Gamma'}
@>{\cores_{M_{2,\bullet}}}>>
(M_{2,\bullet})_{\Gamma} \\
@VVV @VVV @VVV \\
(M_{1,\bullet})_{\Gamma}
@>{\tr_{M_{1,\bullet}}}>>
(M_{1,\bullet})_{\Gamma'}
@>{\cores_{M_{1,\bullet}}}>>
(M_{1,\bullet})_{\Gamma}
\end{CD}
$$
of complexes of abelian groups.

\subsection{}
Let $Y_\bullet \supset Z_\bullet$
be simplicial complexes with (compatible)
$\Gamma$-actions.
Set $M_{1,\bullet} = 
C_\bullet^\cell(Y_\bullet,\Z)/
C_\bullet^\cell(Z_\bullet,\Z)$ 
and
$M_{2,\bullet}=C_\bullet^\cell(|Y_\bullet| \times |E\Gamma_\bullet|,\Z)/
C_\bullet^\cell(|Z_\bullet| \times |E\Gamma_\bullet|,\Z)$.

The quotient maps of spaces give a morphism
$M_{2, \bullet} \to M_{1, \bullet}$
of complexes.

By taking the $(d-1)$-st homology, we obtain the following
commutative diagram of abelian groups:
$$
\begin{CD}
H_{d-1}^\Gamma(Y_\bullet, Z_\bullet; \Z)
@>>>
H_{d-1}^{\Gamma'}(Y_\bullet, Z_\bullet; \Z)
@>>>
H_{d-1}^\Gamma(Y_\bullet, Z_\bullet; \Z) \\
@V{\beta_\Gamma}VV @V{\beta_{\Gamma'}}VV 
@VV{\beta_\Gamma}V \\
H_{d-1}(\Gamma \bsl Y_\bullet, \Gamma \bsl Z_\bullet; \Z)
@>>>
H_{d-1}(\Gamma \bsl Y_\bullet, \Gamma \bsl Z_\bullet;\Z)
@>>>
H_{d-1}(\Gamma \bsl Y_\bullet, \Gamma \bsl Z_\bullet;\Z).
\end{CD}
$$
Here $\beta_{\Gamma}$ and $\beta_{\Gamma'}$ 
are the maps induced by the quotient maps
for the groups $\Gamma$ and $\Gamma'$ respectively.

By taking the cokernels of the vertical arrows, we obtain
the maps
$$
\Coker\, \beta_{\Gamma} \xto{\alpha_1}
\Coker\, \beta_{\Gamma'} \xto{\alpha_2}
\Coker\, \beta_{\Gamma}.
$$
By construction, the composite $\alpha_2 \circ \alpha_1$
is equal to the map given by the multiplication by $[\Gamma:\Gamma']$.

\subsection{}
\begin{cor}
\label{cor:exponent alpha}
Let $v \neq \infty$ be a prime of $F$, and
let $F_v$ denote the completion of $F$ at $v$.
Let $\bK_v$ be a pro-$p$ open compact subgroup of $\GL_d(F_v)$.
Let us consider the intersection $\Gamma' = \Gamma \cap \bK_v$
in $\GL_d(F_v)$.
Then 
\[
H_{d-1}^\Gamma(X, X';\Z)
\to
H_{d-1}(\Gamma \backslash X, \Gamma\backslash X';\Z)\]
is injective 
and the cokernel is annhilated by $p^{e(d)} [\Gamma:\Gamma']$.
\end{cor}
\begin{proof}
It follows from Lemma~\ref{lem:exponent prop} 
that 
$\Coker \beta_\Gamma'$
is killed by $p^{e(d)}$.
Hence by the discussion above,
$\Coker \beta_\Gamma$
is killed by $p^{e(d)} [\Gamma:\Gamma']$.
This implies the claim.
\end{proof}

\subsection{}
We prove the analogue of Theorem~\ref{lem:apartment}(4).
\begin{cor}
\label{cor:pedNd}
 Let $v_0 \neq \infty$ be a prime of $F$ such that the cardinality 
$q_0$ of the residue field $\kappa(v_0)$ at $v_0$ is smallest 
among those at the primes $v \neq \infty$.
Set $N(d) = \prod_{i=1}^{d} (q_0^i-1)$.
Then 
the cokernel of the injective map
\[
H_{d-1}^\Gamma(X, X';\Z)
\to
H_{d-1}(\Gamma \backslash X, \Gamma\backslash X';\Z)
\]
is annihilated by
$p^{e(d)} N(d)$.
\end{cor}
\begin{proof}
Let $\alpha$ be a positive real number.
We use the discussion above with
$Y_\bullet=X_\bullet$
and $Z_\bullet=X^{(\alpha)}_\bullet$.

Since $\Gamma$ is an arithmetic subgroup,
$\Gamma$ is contained, as a subgroup of $\GL_d(F_{v_0})$,
in a compact open subgroup of $\GL_d(F_{v_0})$.
Let $\cO_{v_0}$ denote the ring of integers in $F_{v_0}$.
Since any maximal compact subgroup of $\GL_d(F_{v_0})$
is a conjugate of $\GL_d(\cO_{v_0})$, there exists
$g \in \GL_d(F_{v_0})$ such that $g^{-1} \Gamma g$
is contained in $\GL_d(\cO_{v_0})$. Let $\kappa(v_0)$
denote the residue field at $v_0$ and choose a
$p$-Sylow subgroup $P$ of $\GL_d(\kappa(v_0))$.
Let $\Gamma'$ denote the inverse image of $P$ under the
composite
$$
\Gamma \xto{f_1} \GL_d(\cO_{v_0}) \xto{f_2} \GL_d(\kappa(v_0)),
$$
where $f_1$ is the map that sends $\gamma \in \Gamma$ to
$g^{-1} \gamma g$ and $f_2$ is the map induced by
the ring homomorphism $\cO_{v_0} \surj \kappa(v_0)$.
Since $g f_2^{-1}(P) g^{-1}$ is a pro-$p$ compact open
subgroup of $\GL_d(F_{v_0})$, it follows from Corollary~\ref{cor:exponent alpha}
that the cokernel of the map $\beta_\Gamma$ above is killed by
$p^{e(d)} [\Gamma:\Gamma']$. Since $[\Gamma:\Gamma']$ divides
$[\GL_d(\kappa(v_0)):P] = N(d)$, the claim follows.
\end{proof}

\section{the Mittag-Leffler condition}
\label{sec:Mittag-Leffler}
Let $\alpha \in \Z$ be a positive integer.
We write $X=\cBT_\bullet$ and 
$X^{(\alpha)}=\cBT_\bullet^{(\alpha)}$ for 
short.
Let
\[
A^{(\alpha)}
=
\mathrm{Image}
(H_{d-1}(X,X^{(\alpha)})
\to
H_0(\Gamma, H_{d-1}(X, X^{(\alpha)}))
\to
H_{d-1}(\Gamma\backslash X, \Gamma\backslash X^{(\alpha)}
))
\]
For $\alpha < \alpha'$, the diagram
\[
\begin{CD}
H_{d-1}(X, X^{(\alpha)})   @>>>   
H_{d-1}(\Gamma\backslash X, \Gamma\backslash X^{(\alpha)}
))
\\
@VVV  @VVV 
\\
H_{d-1}(X, X^{(\alpha')})   @>>>   
H_{d-1}(\Gamma\backslash X, \Gamma\backslash X^{(\alpha')}))
\end{CD}
\]
is commutative (the maps are those induced by 
the inclusion map $X^{(\alpha')} \subset X^{(\alpha)}$).
Since the left vertical map is an isomorphism
(both are isomorphic to $H_{d-2}(T_{F^{\oplus d}}, \Z)$),
we obtain that $A^{(\alpha)}$ surjects to $A^{(\alpha')}$.
Hence the projective system
\[
\{A^{(\alpha)}\}_{\alpha}
\]
satisfies the Mittag-Leffler condition.

We set 
\[
C^{(\alpha)}
=
\mathrm{Coker}
(H_{d-1}(X,X^{(\alpha)})
\to
H_0(\Gamma, H_{d-1}(X, X^{(\alpha)})
\to
H_{d-1}(\Gamma\backslash X, \Gamma\backslash X^{(\alpha)})
))
\]
so that we have an exact sequence 
\[
0 
\to 
A^{(\alpha)}
\to 
H_{d-1}(\Gamma\backslash X, \Gamma\backslash X^{(\alpha)})
\to
C^{(\alpha)}
\to 
0
\]
for each $\alpha$.
Note that each $C^{(\alpha)}$ is a finite abelian group.

Because the Mittag-Leffler condition is satisfied, the following sequence
\[
0 
\to 
\varprojlim_\alpha A^{(\alpha)}
\to 
\varprojlim_\alpha H_{d-1}(\Gamma\backslash X, \Gamma\backslash X^{(\alpha)})
\to
\varprojlim_\alpha  C^{(\alpha)}
\to 
0
\]
is also exact
(see \cite[p.83]{WeHA}).

It follows that the cokernel 
of the map
\begin{eqnarray}
\label{eqn:almost all}
H_{d-2}(T_{F^{\oplus d}},\Z)
\to
\varprojlim_\alpha 
H_{d-1}(\Gamma \backslash \cBT_\bullet, \Gamma\backslash \cBT_\bullet^{(\alpha)};\Z)
\end{eqnarray}
is annihilated by 
the least common multiple of 
the exponents of $C^{(\alpha)}$.
We turn the corollaries above into
the following propositions.

\section{Proof}
\subsection{}
\begin{lem}
\label{lem:BM limit isom}
We have 
\[
H_{d-1}^\BM
(\Gamma \backslash \cBT_\bullet,
\Z)
\cong
\varprojlim_\alpha 
H_{d-1}(\Gamma \backslash \cBT_\bullet, \Gamma\backslash \cBT_\bullet^{(\alpha)};\Z).
\]
\end{lem}
\begin{proof}
There is a canonical map 
$H_{d-1}^\BM(\Gamma \backslash \cBT_\bullet, \Z)
\to H_{d-1}(\Gamma \backslash \cBT_\bullet, \Gamma\backslash \cBT_\bullet^{(\alpha)}; \Z)$ which 
assemble to give a map to the limit.
Then the claim follows from Lemma~\ref{7_lem:finiteness}.
\end{proof}

\subsection{}
\label{sec:AR ms}
We consider the composite of the map \eqref{eqn:almost all}
and the map in Lemma \ref{lem:BM limit isom}:
\[
H_{d-2}(T_{F^{\oplus d}},\Z)
\to 
H_{d-1}^\BM(\Gamma\backslash \cBT_\bullet, \Z).
\]
Let us write $\MS(\Gamma)_{AR}$
for the image of the map.
By Proposition \ref{prop:MS Z-basis}, this is the submodule generated
by the images of the universal 
modular symbols.   We write $AR$ 
to mean Ash-Rudolph modular symbols.

\subsection{}
Let us prove the following proposition, 
which will imply Theorem~\ref{lem:apartment}(2) later.
\begin{prop}
\label{prop:thm13(2) AR}
Suppose $\bK \subset \GL_d(\A^\infty)$ be 
a pro-$p$ compact open subgroup.
Let $\Gamma=\GL_d(F) \cap \bK$.
Then 
\[
p^{e(d)} H^\BM_{d-1}(\Gamma \bsl \cBT_\bullet, \Z)
\subset \MS(\Gamma)_{AR}.
\]
\end{prop}
\begin{proof}
This follows from Lemmas~\ref{lem:exponent prop}, \ref{lem:Solomon},
the surjection at the end of Section~\ref{sec:LHS ss},
and the discussion in Section~\ref{sec:Mittag-Leffler}
on the inverse limit.
\end{proof}

\subsection{}
Let us prove the following proposition, 
which will imply Theorem~\ref{lem:apartment}(3) later.
\begin{prop}
\label{prop:thm13(3) AR}
Let $v \neq \infty$ be a prime of $F$, and
let $F_v$ denote the completion of $F$ at $v$.
Let $\bK_v$ be a pro-$p$ open compact subgroup of $\GL_d(F_v)$.
Let us consider the intersection $\Gamma' = \Gamma \cap \bK_v$
in $\GL_d(F_v)$.
Then 
\[
p^{e(d)} [\Gamma:\Gamma'] H^\BM_{d-1}(\Gamma \bsl \cBT_\bullet, \Z)
\subset \MS(\Gamma)_{AR}.
\]
\end{prop}
\begin{proof}
This follows from Corollary~\ref{cor:exponent alpha}, 
Lemma~\ref{lem:BM limit isom}
and the discussion in Section~\ref{sec:Mittag-Leffler}
on the inverse limit.
\end{proof}

\subsection{}
We prove the analogue of Theorem~\ref{lem:apartment}(4).
\begin{prop}
\label{prop:thm13(4) AR}
 Let $v_0 \neq \infty$ be a prime of $F$ such that the cardinality 
$q_0$ of the residue field $\kappa(v_0)$ at $v_0$ is smallest 
among those at the primes $v \neq \infty$.
Set $N(d) = \prod_{i=1}^{d} (q_0^i-1)$.
Then 
\[
p^{e(d)} N(d) 
H^\BM_{d-1}(\Gamma \bsl \cBT_\bullet, \Z)
\subset 
\MS(\Gamma)_{AR}
\].
\end{prop}
\begin{proof}
This follows from Corollary~\ref{cor:pedNd}, 
Lemma~\ref{lem:BM limit isom}
and the discussion in Section~\ref{sec:Mittag-Leffler}
on the inverse limit.
\end{proof}

\begin{cor}
\label{cor:thm13(1) AR}
We have 
\[
H^\BM_{d-1}(\Gamma \bsl \cBT_\bullet, \Q)
=
\MS(\Gamma)_{AR}\otimes \Q
\]
\end{cor}
\begin{proof}
This follows immediately from the proposition.
\end{proof}

\section{the case $d=2$}
When $d=2$, we do not need the computation of the 
homology groups of the stabilizer groups.
Hence we obtain the best result that the exponent is 1.
\begin{prop}
\label{prop:thm13(5) AR}
When $d=2$, we have
\[
H^\BM_{1}(\Gamma \bsl \cBT_\bullet, \Z)
=\MS(\Gamma)_{AR}.
\]
\end{prop}
\begin{proof}
When $d=2$, in the argument given in Section~\ref{sec:d argument},
we only need to look at one edge map 
in the spectral sequence.   It is easy to see that the 
edge map is surjective.   Hence the claim follows.
\end{proof}

\chapter{Comparison of modular symbols}
\label{ch:compare ms}
We defined two kinds of modular symbols.
Our modular symbols come from fundamental classes of 
apartments, whereas the modular symbols of Ash-Rudolph
(or those coming from the universal modular symbols)
originate in the homology of the Tits building.
In this chapter, we show that these two definitions
coincide.   We refer to Section~\ref{sec:outline comparison}
for the outline.

We have a proof of our main theorem (Theorem~\ref{lem:apartment}).

\section{Statement of comparison proposition}
The goal of this chapter is to prove the 
following proposition.

Let $q_1, \dots, q_d$ be an ordered basis
of $F^{\oplus d}$.   In Section~\ref{sec:def modular symbol},
we defined (our) modular symbol in the Borel-Moore
homology of an arithmetic quotient.   Let us write $[q_1, \dots, q_d]_A$
for it.

In Section~\ref{sec:univ ms}, we defined universal modular symbol
$[q_1, \dots, q_d]$ in the homology of the Tits building.

In Section~\ref{sec:AR ms}, we defined the modular symbol
of Ash-Rudolph to be the image of the universal modular symbol
in the Borel-Moore homology of an arithmetic quotient.
Let us write $[q_1,\dots, q_d]_{AR}$ for it.

The main task of this chapter is to show that they coincide. 
\begin{prop}
\label{PROP:COINCIDENCE MS}
Let the notation be as above.
We have 
\[
[q_1, \dots, q_d]_A=[q_1, \dots, q_d]_{AR}.
\]
In particular, we have $\MS(\Gamma)=\MS(\Gamma)_{AR}$.
\end{prop}

\section{{Proof of Theorem~\ref{lem:apartment}}}
Using this proposition, we are now able to prove our main theorem.
\begin{proof}(Proof of Theorem~\ref{lem:apartment})

Using Proposition~\ref{PROP:COINCIDENCE MS},
the statements follow from
Propositions~\ref{prop:thm13(2) AR}, \ref{prop:thm13(3) AR},
\ref{prop:thm13(4) AR}, \ref{prop:thm13(5) AR} and Corollary \ref{cor:thm13(1) AR}.
\end{proof}

\section{Outline}
\label{sec:CD comparison}
\subsection{}
Let $\alpha$ be a sufficiently large real number.
We defined a subsimplicial complex 
$\cBT_\bullet^{(\alpha)} \subset \cBT_\bullet$ 
in Section~\ref{sec:Tits}.

Given an ordered basis $q_1, \dots, q_d$, 
we defined an injective map 
$\phi=\phi_{q_1,\dots, q_d}: A_\bullet \to \cBT_\bullet$
in Section~\ref{sec:521}.
We regard $A_\bullet$ as a subsimplicial complex
via $\phi$.

We set $A_\bullet^{(\alpha)}=A_\bullet \cap \cBT_\bullet^{(\alpha)}$.

\subsection{}
The $S^{d-2}$ is the $(d-2)$-sphere.

We use the notation from Section~\ref{sec:Tits building},
We let $W=F^{\oplus d}$ and write $T_d=T_W$. 
The notation $\cP'(B)$
appeared in Section~\ref{sec:barycentric subdivision}.
We use $|X|$ to denote the geometric 
realization of a poset $X$.   (See below for precise definition.)

\subsection{}
To prove the proposition, we look at the following diagram 
which appeared in Section~\ref{sec:outline comparison}.
All the coefficients are $\Z$.    
The map (1) is from Lemma~\ref{lem:BM limit isom}.
The map (2) is the limit of the pushfoward by
the canonical projection for each $\alpha$.
The map (3) is the canonical 
projection to one of the 
entries in the limit.  It is an isomorphism
because the maps $(4)(5)$ are isomorphisms.
The map (5) is from Corollary~\ref{cor:Grayson htpy eq}.
The map (4) is the boundary map in the long exact sequence
for a pair of spaces.   It is an isomorphism because $\cBT_\bullet$
is contractible.

The map (6) is the canonical projection.
The map (7) is the boundary map, which
is an isomorphism because $A_\bullet$ 
is contractible.
The map (8) will be constructed below.

The map (9) is the pushforward by 
the map induced by the morphism 
of posets $\cP'(B) \to Q(W)$
associated with the ordered basis.

The map (10) appeared in Section~\ref{sec:def modular symbol}.

The inclusions $A_\bullet \subset \cBT_\bullet$
and $A_\bullet^{(\alpha)} \subset \cBT_\bullet^{(\alpha)}$
give the rest of the horizontal maps.
\begin{equation}
\label{main diagram}
\begin{tikzcd}
H_{d-1}^\BM(\Gamma \backslash \cBT_\bullet)   \arrow[d,"\sim", "(1)"']                                                       &      H_{d-1}^\BM(A_\bullet)  \arrow[dd, "\sim"]  \arrow[l, "(10)"] 
\\
\varprojlim_\alpha H_{d-1}(\Gamma \backslash \cBT_\bullet, \Gamma\backslash \cBT_\bullet^{(\alpha)})    &    \\
\varprojlim_\alpha H_{d-1}(\cBT_\bullet, \cBT_\bullet^{(\alpha)})  \arrow[u, "(2)"']  \arrow[d, "\sim", "(3)"'] &   \varprojlim_\alpha H_{d-1}(A_\bullet, A_\bullet^{(\alpha)})  \arrow[l] \arrow[d, "\sim", "(6)"']
\\
H_{d-1}(\cBT_\bullet, \cBT_\bullet^{(\alpha)}) \arrow[d, "\sim", "(4)"']                                                   &    H_{d-1}(A_\bullet, A_\bullet^{(\alpha)})  \arrow[l]  \arrow[d, "\sim", "(7)"'] 
\\
H_{d-2}(\cBT_\bullet^{(\alpha)})   \arrow[d, "\sim", "(5)"']                                                                     &    H_{d-2}(A_\bullet^{(\alpha)})  \arrow[d, "\sim", "(8)"'],  \arrow[l]                      
\\
H_{d-2}(T_d)                                                                                                                                  &  H_{d-2}(|\cP'(B)|)                    \arrow[l, "(9)"]      
\\
                                                                                                                                                  &  H_{d-2}(S^{d-2})                     \arrow[u,"\sim"']
\end{tikzcd}
\end{equation}
The commutativity of the rectangles except 
the bottom one (with $(5)(8)(9)$) is easy to 
see.

\subsection{}
Recall the definition of 
$[q_1,\dots, q_d]_{AR} \in 
H_{d-1}^\BM(\Gamma \backslash \cBT_\bullet)$
from Section~\ref{sec:AR ms}.
We have 
the universal modular symbol $[q_1, \dots, q_d]$
in $H_{d-2}(T_d)$ of the left bottom corner,
and $[q_1, \dots, q_d]_{AR}$
is the image via the left vertical maps
of $[q_1, \dots, q_d]$.

By construction, $|\cP'(B)|$ is canonically 
homotopy equivalent to $S^{d-2}$.
By the definition (Section~\ref{sec:univ ms}) 
of universal modular symbol,
$[q_1,\dots, q_d]$ is the image of the 
fundamental class of $H_{d-1}(|\cP'(B)|)$
via the map (9).

Our modular symbol $[q_1,\dots, q_d]_{A}$
is the image by (10)
of the fundamental class of the apartment
(Section~\ref{sec:def modular symbol}).

Thus to prove the proposition,
the strategy is to construct the map (8)
such that the bottom square is commutative.

\section{Quillen's lemma}
We put our emphasis on posets in what follows.
For some terminology and results, 
we refer to \cite[p.102--103]{Quillen}
\subsection{}
For a poset $W$, we defined the classifying space 
as the simplicial complex 
$(W, \cP_\tot(W))$ where 
$\cP_\tot(W)$ is the set of finite totally ordered subsets of $W$.
We let $|W|$ denote the 
geometric realization 
of the classifying space. 
We call $W$ contractible if $|W|$ is contractible.
We say that a morphism of posets $W_1 \to W_2$
is a homotopy equivalence if the induced map
$|W_1| \to |W_2|$ of geometric realizations is a homotopy equivalence.

\subsection{}
We use the following lemma of Quillen.
\begin{lem}[Cor 1.8 {\cite{Quillen}}]
\label{lem:Quillen poset}
Let $X, Y$ be posets and $Z \subset X \times Y$
be a closed subset
(i.e., $z' \leq z \in Z$ implies $z' \in Z$).

Let $p_1:Z \to X, p_2:Z \to Y$
be the maps induced by the projections.

We set 
$Z_x=\{ y \in Y \,|\, (x,y) \in Z\}$
and 
$Z_y=\{ x \in X \,|\, (x,y) \in Z\}$.

If $Z_x$ (resp. $Z_y$)
contractible for all $x \in X$ (resp. all $y \in Y$)
then $p_1$ (resp. $p_2$) are homotopy equivalences.
\end{lem}

\section{A key diagram of posets}
\subsection{}
For a simplicial complex $U_\bullet$,
we define $(U_\bullet)$ to be the 
poset of simplicies of $U_\bullet$.
The classifying space of $(U_\bullet)$
is the barycentric subdivision of $U_\bullet$.
(See proof of Lemma 1.9 of \cite{Gra}.)

\subsection{}
We construct the following commutative diagram of posets:
\begin{equation}
\label{key diagram}
\begin{CD}
(\cBT_\bullet^{(\alpha)})    @<<{h_1}<   (A(B)_\bullet^{(\alpha)})    \\
@A{f_1}AA                                        @A{g_1}AA                    \\
U                                  @<<{h_2}<    T                                \\
@V{f_2}VV                                        @V{g_2}VV                    \\
Q(F^{\oplus d})                               @<<{h_3}<     \cP'(B) 
\end{CD}
\end{equation}
We will see that 
the maps $f_1, f_2, g_1, g_2$ are homotopy equivalences.
Upon taking the geometric realizations, 
the diagram gives rise to the bottom square of the diagram \eqref{main diagram}.

\section{The left column}
We construct the left column of the key diagram.   
This is a slight generalization of the construction 
in Grayson \cite{Gra}.

\subsection{}
We recall some notation from Section~\ref{subsec:X_K}.
Let $d \ge 1$ and $\bK=\GL_d(\widehat{A})$.
We defined $\widetilde{X}_{\bK,\bullet}$ 
and $X_{\bK, \bullet}$.    The nonadelic description is
$X_{\bK, \bullet}=\coprod_{j\in J_\bK} \Gamma_j \backslash \cBT_\bullet$.

We consider the maps 
$\phi_j: \cBT_\bullet \to \Gamma_j \backslash \cBT_\bullet$
for $j \in J_{\bK, \bullet}$.
We set 
$\cBT_\bullet^{(\alpha,j)}=\phi_j^{-1}(X_{\bK, \bullet}^{(\alpha)})$.
We note that this does not depend on the choice of $\bK$.
We defined $\cBT_\bullet^{(\alpha)}$ earlier.
We have $\cBT_\bullet^{(\alpha)}=\cBT_\bullet^{(\alpha,e)}$
where $e$ is the identity element. 

\subsection{}
Let us write $F_d=Q(F^{\oplus d})$ for short.
We define a map
\[
S_j:\cBT_0 \to F_d \cup \{\emptyset\}
\]
as follows.   We have the map
$\phi_{j,0}: \cBT_0 \to X_{\bK,0}$
of 0-simplices of $\phi$.
Let $\sigma \in \cBT_0$.
Then by Section~\ref{sec:locfree}, 
$\phi_{j,0}(\sigma) \in X_{\bK,0}$
is represented by a chain of locally free $\cO_C$-submodule of rank $d$
of $\eta_*\eta^*\cO_C^{\oplus d}$.
Take one of the sheaves, say $E$, in the chain and consider 
its Harder-Narasimhan filtration
\[
0 \subsetneqq E_0 \subsetneqq E_1 \subsetneqq \dots
\subsetneqq E_r=E \subset \eta_*F^{\oplus d}
\]
The restriction to the generic fiber 
$\Spec F \subset C$ gives a flag of $F^{\oplus d}$.
We let $S_j(\sigma)$ denote the 
corresponding element in 
$F_d\cup \{\emptyset\}$.

\subsection{}
For sufficiently large $\alpha$, the $0$-simplices of 
$\cBT^{(\alpha)}_\bullet$ consist of unstable ones only.
Thus the map $S$ gives 
$S_j: \cBT_0^{(j, \alpha)} \to F_d$.

Let $f \in F_d$ be a flag.  
We let 
$Z_{f, \bullet}^{(j, \alpha)} \subset 
\cBT_\bullet^{(j, \alpha)}$
be the subsimplicial complex 
consisting of 
simplices $\sigma$
such that
\[
S(v) \ge f
\]
for all vertices $v$ of $\sigma$.

\begin{lem}
\label{lem:Z contractible}
$Z_{f,\bullet}^{(j, \alpha)}$ is contractible.
\end{lem}
\begin{proof}
(We use the proof of Lemma~\ref{7_contract} almost word-for-word.)
We set $X=Z_{f, \bullet}^{(j,\alpha)}$.
We proceed by induction on $d$.
Let $f$ be the flag
\[
0 \subset V_1 \subset \cdots \subset V_{r-1} \subset F^{\oplus d}
\]
with $\dim V_j=i_j$ for $1 \le j \le r-1$
so that $\cD=\{i_1, \dots, i_{r-1}\}$
with $i_1<\cdots <i_{r-1}$.

Set $d'=d-i_1$ and 
$\cD'=
\{i'-i_1 
\,|\,
i'\in \cD \setminus \{i_1\}\}
\subset
\{1,\dots, d'-1\}
$.

We define $f_1 \in \Flag_{\cD'}$
as the image of the flag $f_1$ 
with respect to the projection
\[
F^{\oplus d}
\to
F^{\oplus d}/
V_1
\]
We take an isomorphism $F^{\oplus d}/V_1 \cong F^{\oplus (d-i_1)}$
so that $f_1$ corresponds to the standard flag $f_0'$:
\[
0 \subset F^{\oplus (i_2-i_1)} 
\subset F^{\oplus (i_3-i_1)}
\subset \dots \subset
F^{\oplus (i_{r-1}-i_1)}
\subset F^{\oplus (d-i_1)}
\]

Take a representative $g_j \in \GL_d(\A^\infty)$
of $j \in J_\bK$.
Consider the map 
that sends an $\cO_C$-submodule
$\cF[g_j, L_\infty] \subset \eta_* F^{\oplus d}$
to the $\cO_C$-submodule
$\cF[g_j, L_\infty]/\cF[g_j, L_\infty]_{(i_1)}
\subset \eta_* F^{\oplus d'}$.
There exists $j' \in J_{\GL_{d-i_1}(\widehat{A})}$ 
such that the above map gives a morphism
$h: X \to X'$ where
we set
$X'=Z_{f'_0, \bullet}^{(j', \alpha)}$ if $\cD'\neq \emptyset$,
and 
$X'=\cBT_{\GL_{d'}, \bullet}$ if $\cD'=\emptyset$.
By inductive hypothesis, $|X'|$ is contractible.

Let $\epsilon:\mathrm{Vert}(X) \to \Z$
and $\epsilon':\mathrm{Vert}(X') \to \Z$
denote the maps that send a 
locally free 
$\cO_C$-module $\cF$
to the integer
$[p_\cF(1)/\det(\infty)]$.

We fix an $\cO_C$-submodule
$\cF_0$ of 
$\eta_*F^{\oplus d}$
whose equivalence class belongs to $X$.
By twisting $\cF_0$ by some power
of $\cO_C(\infty)$
if necessary, 
we may assume that $p_{\cF_0}(i_1)-p_{\cF_0}(i_1-1)>\alpha$.
We fix a splitting $\cF_0=\cF_{0,(i_1)}\oplus \cF_0'$.
This splitting induces an isomorphism
$\varphi:\eta_*\eta^*\cF_0' \cong \eta_*F^{\oplus d}$.

Let $h':X' \to X$ denote the morphism that sends an $\cO_C$-submodule
$\cF' \subset 
\eta_*\eta^*F^{\oplus d'}$
to the $\cO_C$-submodule
$\cF_{0,(i_1)}(\epsilon'(\cF')\infty)\oplus \varphi(\cF')$
of $\eta_*F^{\oplus d}$.
For each $n \in \Z$, define a morphism
$G_n:X\to X$ by sending an 
$\cO_C$-submodule
$\cF$ of $\eta_*\eta^*F^{\oplus d}$
to the $\cO_C$-submodule
$\cF_{0,(i_1)}((n+\epsilon(\cF))\infty)+\cF$
of $\eta_*F^{\oplus d}$.

Then the
argument in~\cite[p. 85--86]{Gra} shows that
$f$ and $|h'|\circ|h| \circ f$ are homotopic
for any map $f:Z \to |X|$ from a compact space $Z$
to $|X|$. Since the map $|h'|\circ|h| \circ f$ factors
through the contractible space $|X'|$, $f$ is
null-homotopic. Hence $|X|$ is contractible.

\end{proof}

\subsection{}
We have 
$(Z_{f, \bullet}^{(j, \alpha)})
\subset
(\cBT_\bullet^{(\alpha)})$
for each flag $f \in F_d$ and $j \in J_\bK$.

\subsection{}
We set 
\[
U=\{(\beta,\gamma) \in (\cBT_\bullet^{(\alpha)}) \times F_d
\,|\,
\beta \in (Z_{\gamma, \bullet}^{(j, \alpha)})
\}.
\]
It is a poset.
It is a closed subset of 
$ (\cBT_\bullet^{(\alpha)}) \times F_d$.
(i.e., $u' \leq u \in U$ implies $u' \in U$).

We define each of the maps 
\[
f_1:U \to (\cBT_\bullet^{(\alpha)})\times F_d \to (\cBT_\bullet^{(\alpha)})
\]
\[
f_2:U \to (\cBT_\bullet^{(\alpha)})\times F_d \to F_d
\]
as the inclusion followed by the projection.
Set
\[
\begin{array}{rl}
U_\sigma
&=
\{f \in F_d \,|\,
(\sigma, f) \in U\}
\\
&=\{f \in F_d \,|\,
S(v) \ge f \text{ for all vertices $v$ of $\sigma$}
\}.
\end{array}
\]
\begin{lem}
For sufficiently large $\alpha$,
the poset $U_\sigma$ is contractible.   
\end{lem}
\begin{proof}
We actually prove that if $\alpha>3 \deg(\infty)$ then 
$U_\sigma$ has a maximal element, hence contractible.
Let us show that $\bigcap_v S(v)$ is a maximal element,
where the intersection is taken
by regarding each $S_j(v)$ as a subset of $Q(F^{\oplus d})$.
The maximality is clear.   We need to show that $\bigcap_v S_j(v)$ is nonempty.

For a locally free $\cO_C$-submodule 
$\mathcal{G} \subset \eta_* F^{\oplus d}$ of rank $d$
and a vector subspace $B \subset F^{\oplus d}$ we denote by
$\mathcal{G} \cap \eta_* B$ the intersection 
of $\mathcal{G}$ and $\eta_* B$ in $\eta_* F^{\oplus d}$.
Let us choose a vertex $v_0$ of $\sigma$.
Let $\cF_0 \subset \eta_* F^{\oplus d}$ 
be a locally free $\cO_C$-submodule 
of rank $d$ that represents $v_0$.
Since $v_0$ belongs to $\cBT_\bullet^{(\alpha,j)}$,
there exists a nonzero vector subspace 
$W_0 \subsetneqq F^{\oplus d}$
such that $\cF_0 \cap \eta_* W_0$ is a part of the Harder-Narasimhan
filtration of $\cF_0$ and 
the minimal slope of the
successive quotients of the Harder-Narasimhan filtration of 
$\cF_0 \cap \eta_* W_0$ is greater than
the sum of $3 \deg(\infty)$ and the maximal slope of the
successive quotients of the Harder-Narasimhan filtration of 
$\cF_0/(\cF_0 \cap \eta_* W_0)$.
Let $v$ be a vertex of $\sigma$. Let us choose
a locally free $\cO_C$-submodule 
$\cF_1 \subset \eta_* F^{\oplus d}$
of rank $d$ that represents a vertex of $\sigma$
in such a way that $\cF_0 \subset \cF_1 \subset \cF_0(\infty)$.
Then for any vector subspace $W' \subset F^{\oplus d}$,
we have $\cF_0 \cap \eta_* W' \subset 
\cF_1 \cap \eta_* W' \subset (\cF_0 \cap \eta_* W')(\infty)$.
Hence we have
$$
0 \le \deg(\cF_1 \cap \eta_* W') - \deg(\cF_0 \cap \eta_* W')
\le \dim_F W'.
$$
This implies that
$\cF_1 \cap \eta_* W_0$ is a part of the Harder-Narasimhan
filtration of $\cF_1$.
Hence $\bigcap_v S_j(v)$ contains $W_0$.
In particular $\bigcap_v S_j(v)$ is nonempty.  
This proves the contractiblity.
\end{proof}

This implies that $f_2$ is a homotopy equivalence
by Lemma~\ref{lem:Quillen poset}.

\subsection{}
Set
\[
\begin{array}{rl}
U_f
&=
\{
\sigma \in (\cBT_\bullet^{(\alpha)})
\,|\,
(\sigma, f) \in U
\}
\\
&=\{
\sigma \in (\cBT_\bullet^{(\alpha)})
\,|\,
S(v) \ge f \text{ for all vertices $v$ of $\sigma$}
\}.
\end{array}
\]
This is $Z_{f, \bullet}^{(j, \alpha)}$,
which is contractible.
Hence by Lemma~\ref{lem:Quillen poset}, 
$f_2$ is a homotopy equivalence.

\section{The right column}
\subsection{}
Now, take an ordered basis
$q_1, \dots, q_d$ of $F^{\oplus d}$.

We define $h_3: \cP_\tot(\cP'(B)) \to F(F^{\oplus d})$
to be the map that appeared in Section~\ref{sec:univ ms}.

\subsection{}
We have a map (from Section~\ref{sec:521})
\[
\varphi_{q_1,\dots, q_d}: A(B)_\bullet
=A_\bullet \to \cBT_\bullet.
\]
We identify $A(B)_\bullet$
and $A_\bullet$ as 
subspaces of $\cBT_\bullet$
via this map.

We set
\[
A_\bullet^{(\alpha)}
=\cBT_\bullet^{(\alpha)} \cap A_\bullet
\]
and define
\[
h_1:
(A_\bullet^{(\alpha)}) 
\to
(\cBT_\bullet^{(\alpha)})
\]
to be the morphism of posets
induced by the inclusion
$A_\bullet^{(\alpha)} 
\subset \cBT_\bullet^{(\alpha)}$
of simplicial complexes.


\subsection{}
Note that the image of the restriction map 
\[
S|_{A_0^{(\alpha)}}:
A_0^{(\alpha)}
\to
F_d
\]
is contained in $\cP'(B)$,
which is regarded as a subset of $F_d$
via $h_3$.

\subsection{}
We set 
\[
T=
\{
(\beta, \gamma) \in
(A_\bullet^{(\alpha)})
\times \cP'(B)
\,|\,
S(v) \ge \gamma 
\text{ for all vertices $v$ of $\beta$}
\}.
\]
It is a closed subset of 
$(A_\bullet^{(\alpha)})
\times \cP'(B)$.

\subsection{}
We define each of the maps
\[
g_1: T \subset 
(A_\bullet^{(\alpha)})
\times \cP'(B)
\to
(A_\bullet^{(\alpha)})
\]
\[
g_2:
 T \subset 
(A_\bullet^{(\alpha)})
\times \cP'(B)
\to
\cP'(B)
\]
as the inclusion followed by the projection.

\subsection{}
Let
\[
\begin{array}{rl}
T_f
&=
\{ \sigma \in (A_\bullet^{(\alpha)})
\,|\,
(\sigma,f) \in T
\}
\\
&=
\{ \sigma \in (A_\bullet^{(\alpha)})
\,|\,
S(v) \ge f \text{ for all vertices $v$ of $\sigma$}
\}.
\end{array}
\]

\begin{lem}
$T_f$ is contractible
\end{lem}
\begin{proof}
We use the setup in the proof of Lemma~\ref{lem:Z contractible}.
We proceed by induction on $d$.
We let $A_\bullet^{(\alpha)} \subset X$.
Take $d',i_1, \cD'$ as before.
We set
$Y'=A_\bullet^{(\alpha)} \subset X'$
if $\cD' \neq \emptyset$,
$Y'~A_\bullet \subset X'$ if 
$\cD'=\emptyset$.
The inclusions are those corresponding 
to the ordered basis $q_{i_1+1},\dots, q_d$
of $F^{\oplus d}/(F^{\oplus i_1}\oplus \{0\}^{\oplus d-i_1}$).

By the inductive hypothesis, $|Y'|$ is contractible.
Notice that the restriction of 
the map $h:X \to X'$
to $Y$ has its image inside $Y'$.
We can also check that the image
of the restriction of $h'$ to $Y'$
is in $Y$.
Then 
the same argument as in the proof of 
Lemma~\ref{lem:Z contractible}
proves that $|Y'|$ is contractible.
\end{proof}
This implies $g_1$ is a homotopy equivalence.

\subsection{}
Let 
\[
\begin{array}{rl}
T_\sigma
&=\{f \in T_d
\,|\,
(\sigma, f) \in T\}
\\
&=
\{ f \in T_d
\,|\,
S(v) \ge f \text{ for all vertices $v$ of $\sigma$} 
\}.
\end{array}
\]
This is contractible because it has 
a maximal element.
This implies $g_2$ is a homotopy equivalence.

\section{Proof of Proposition \ref{PROP:COINCIDENCE MS}}
The map (5) in the diagram is defined to be
the isomorphism of homology groups
induced by $f_1 \circ f_2^{-1}$.
The map (8) is defined to be 
the isomorphism of homology groups 
induced by $g_1 \circ g_2^{-1}$.

By the commutativity of the diagram~\eqref{key diagram},
we have the commutativity of the 
bottom square of the diagram~\eqref{main diagram}.
This proves Proposition~\ref{PROP:COINCIDENCE MS}.

\appendix

\backmatter

\printindex

\end{document}